\documentclass[a4paper, geqno, 12pt]{article}

\usepackage{enumerate}
 \usepackage{enumitem}

 \usepackage{amsfonts}
\usepackage{url}
\usepackage[hypertexnames=false]{hyperref}
\usepackage{ulem}
\usepackage{amsmath}
\usepackage{mathabx}
\usepackage[french,english]{babel}
\usepackage[utf8]{inputenc}
\usepackage[T1]{fontenc}
\usepackage{latexsym}
\usepackage{amsthm}
\usepackage{amssymb}
\usepackage{amstext}
\usepackage[all]{xy}
\usepackage{amssymb}
\usepackage{graphicx}
\usepackage{color}
\usepackage{fancyhdr}
\usepackage{lastpage}
\usepackage{enumitem}
\usepackage{graphicx,wrapfig,tikz}
\usepackage{capt-of}

\newtheorem{Df}{\iflanguage{french}{D\'efinition}{Definition}}[section]
\newtheorem{Prop}[Df]{Proposition}

\newtheorem{Th}[Df]{\iflanguage{french}{Théorème}{Theorem}}
\newtheorem{Cor}[Df]{\iflanguage{french}{Corollaire}{Corollary}}

\usepackage{authblk}

\title{Steinitz Classes of Unramified Galois Extensions}
\date{}
\author[1,2]{Ayoub Al Uartassi}
\author[1]{Boucha\"{\i}b Soda\"{\i}gui \thanks{Corresponding author: \texttt{bouchaib.sodaigui@uphf.fr}}}

\author[3]{Mohammed Taous}

\affil[1]{Université Polytechnique Hauts-de-France,  Ceramaths,
    FR CNRS 2037,
    F-59313 Valenciennes, France}
\affil[2]{Université
Moulay Ismaïl, Faculté des Sciences, Meknès, Maroc}
\affil[3]{Université Mohammed Premier, Faculté des Sciences, Département de mathématiques, Oujda, Maroc}

\begin{document}
\maketitle{}
\begin{quote}
\small\textbf{Note on the English translation.}
This is an English translation of the French article by Ayoub Al Uartassi,
Boucha\"{i}b Soda\"{i}gui, and Mohammed Taous,
\emph{Classes de Steinitz d'extensions galoisiennes non ramifiées},
published in \emph{Acta Arithmetica} \textbf{219} (2025), no.~1, 53--80,
\href{https://doi.org/10.4064/aa240611-17-10}{doi:10.4064/aa240611-17-10}.

The English version restricts the unprimed inclusion in Theorem~\ref{t2} to groups of even order with a cyclic Sylow $2$-subgroup. The subgroup conclusion for the primed sets is unchanged.
\end{quote}

\centerline { Mathematics Subject Classification 2020: 11R33}

\centerline{\bf {Abstract}}

Let $K$ be a number field and $Cl_{K}$ its class group. Let $\Gamma$ be a finite group. We denote by $ R_{nr}(K,\Gamma)$  the subset of $Cl_{K}$  formed by the  classes which are realizable as Steinitz classes of   Galois extensions  of $K$,  unramified at the finite places of $K$ and having a Galois group isomorphic to $\Gamma$. If $\Gamma$ is nontrivial and abelian and the narrow class number of $K$ is prime to the order of $\Gamma$, then  $R_{nr}(K,\Gamma)=\emptyset$. In the present article, for any $\Gamma$ we consider the set
 $R'_{nr}(K,\Gamma):=\lbrace 1 \rbrace \cup R_{nr}(K,\Gamma)$.
  We prove that  $R'_{nr}(K, \mathbb{Z}/2\mathbb{Z})$  is a subgroup of  $Cl_{K,2}:=\lbrace c\in Cl_{K}\vert c^{2}=1\rbrace$ ; furthermore, it is equal to $Cl_{K,2}$ under a certain assumption on $K$. Using this result   we show that $R'_{nr}(K,\Gamma)$ is a subgroup of $R'_{nr}(K, \mathbb{Z}/2\mathbb{Z})$ if $\Gamma$ either  has odd order, or  has a noncyclic Sylow $2$-subgroup (for instance $\Gamma$ a nonabelian $2$-group,  $\Gamma=S_n$ or $A_n$, with $n\geq 4$), or has a normal cyclic Sylow $2$-subgroup
  (for instance  $\Gamma$ nilpotent having even order).


\section{Introduction}

Throughout this paper, for a number field $K$, we write $\mathcal{O}_{K}$ for its ring of integers, $Cl_{K}$ for its ideal class group, $h_K$ for its class number, $Cl^{+}_{K}$ for its narrow class group, and $h^{+}_K$ for its narrow class number. For a fractional ideal $I$ of $K$, we write $[I]_K$ for its class in $Cl_K$. We call a finite extension $L/K$ unramified if it is unramified at all finite places. Whenever we require it to be unramified at the infinite places as well, we state this explicitly.

Recall the definition of the Steinitz class. Let $L/K$ be a finite extension of number fields of degree $n$. Since $\mathcal{O}_{L}$ is a torsion-free $\mathcal{O}_{K}$-module of rank $n$, there exists a fractional ideal $I$ of $\mathcal{O}_{K}$ such that $\mathcal{O}_{L}\simeq\mathcal{O}_{K}^{n-1}\oplus I$ as $\mathcal{O}_{K}$-modules (see \cite[Theorem 1.2.19, p. 9]{cohan}). The class $[I]_K$ is called the Steinitz class of $\mathcal{O}_{L}$, or of the extension $L/K$, and is denoted by $cl_{K}(\mathcal{O}_{L})$. For a fixed degree $n$, this class completely determines the $\mathcal{O}_{K}$-module structure of $\mathcal{O}_{L}$. In particular, $\mathcal{O}_{L}$ is free over $\mathcal{O}_{K}$ if and only if $cl_{K}(\mathcal{O}_{L})=1$.

Let $\Gamma$ be a finite group and let $K$ be a number field. A $\Gamma$-extension of $K$ is a Galois extension of $K$ whose Galois group is isomorphic to $\Gamma$. We write $R_m(K,\Gamma)$ (respectively, $R_{nr}(K,\Gamma)$) for the set of classes $c\in Cl_K$ such that there is a tamely ramified (respectively, unramified) $\Gamma$-extension $L/K$ with $c=cl_K(\mathcal{O}_L)$.

It is conjectured that $R_m(K,\Gamma)$ is a subgroup of $Cl_K$ (see \cite[Conjecture 3, \S1, p. 6]{BGS}). When $\Gamma$ is abelian, the conjecture follows from McCulloh's work \cite{Mc}, as explained below. In the nonabelian case it remains open in general. For results establishing it for particular nonabelian groups, see \cite{BrS,BGS,BS,SOCA,bouchaib}.

If $L/K$ is a tamely ramified $\Gamma$-extension, then $\mathcal{O}_L$ determines a class $[\mathcal{O}_L]$ in the locally free class group $Cl(\mathcal{O}_K[\Gamma])$. We write $R_m(\mathcal{O}_K[\Gamma])$ (respectively, $R_{nr}(\mathcal{O}_K[\Gamma])$) for the set of classes $c$ in this group for which there exists a tamely ramified (respectively, unramified) $\Gamma$-extension $L/K$ with $c=[\mathcal{O}_L]$.

For abelian $\Gamma$, McCulloh \cite{Mc} studied the set $\mathfrak{R}_m(\mathcal{O}_K[\Gamma])$ of classes in $Cl(\mathcal{O}_K[\Gamma])$ realized by the integral closures of $\mathcal{O}_K$ in certain tamely ramified Galois algebras over $K$ (see \cite[(1.1), p. 263, the discussion following (2.2), p. 271, and the first paragraph of \S3, p. 276]{Mc}). He also studied the subset $\mathfrak{R}_{nr}(\mathcal{O}_K[\Gamma])$ obtained by requiring these algebras to be unramified at the finite places. Using a Stickelberger map, he proved that both sets are subgroups of $Cl(\mathcal{O}_K[\Gamma])$ (see \cite[Corollary (6.20), p. 291]{Mc}). Furthermore, by \cite[Theorem (6.17), p. 291]{Mc},
\[
\mathfrak{R}_m(\mathcal{O}_K[\Gamma])=R_m(\mathcal{O}_K[\Gamma]).
\]
Thus $R_m(\mathcal{O}_K[\Gamma])$ is a subgroup of $Cl(\mathcal{O}_K[\Gamma])$. Forgetting the $\Gamma$-action then shows that $R_m(K,\Gamma)$ is a subgroup of $Cl_K$.

Agboola and McCulloh refined this result in \cite[Theorem D]{agboola} and substantially generalized it in \cite[Theorem E]{agboola}. Under certain hypotheses on $K$, they showed that $\mathfrak{R}_m(\mathcal{O}_K[\Gamma])$ is a subgroup of $Cl(\mathcal{O}_K[\Gamma])$ for a broad family of groups that includes nonabelian groups. Their work leaves open the corresponding study of $\mathfrak{R}_{nr}(\mathcal{O}_K[\Gamma])$.

McCulloh \cite{Mc} did not compare $\mathfrak{R}_{nr}(\mathcal{O}_K[\Gamma])$ and $R_{nr}(\mathcal{O}_K[\Gamma])$. This comparison was the starting point of the present paper. If these sets were equal, forgetting the $\Gamma$-action would show that $R_{nr}(K,\Gamma)$ is a subgroup of $Cl_K$ for abelian $\Gamma$.

Recall that the Hilbert class field $H_K$ of $K$ is its maximal abelian extension unramified at all finite and infinite places, and $Gal(H_K/K)\simeq Cl_K$. The narrow Hilbert class field $H_K^+$ is its maximal abelian extension unramified at all finite places, and $Gal(H_K^+/K)\simeq Cl_K^+$. If $h_K^+$ is odd and $\Gamma$ is abelian of even order, then both $R_{nr}(\mathcal{O}_K[\Gamma])$ and $R_{nr}(K,\Gamma)$ are empty. In this case $\mathfrak{R}_{nr}(\mathcal{O}_K[\Gamma])\ne R_{nr}(\mathcal{O}_K[\Gamma])$.

For an arbitrary finite group $\Gamma$, we therefore define
\[
R'_{nr}(K,\Gamma):=\{1\}\cup R_{nr}(K,\Gamma).
\]
This set is the main object of study in this paper.

We say that $K$ is of class type if every quadratic extension $K(\sqrt{\varepsilon})/K$ with $\varepsilon$ a unit of $K$ is ramified at some finite place. Otherwise, we say that $K$ is of unit type (see Definition \ref{Dftypeclasse} and the paragraph that follows). If $h_K^+$ is odd, then $K$ is of class type. If $K$ is of unit type, then $h_K^+$ is even.

Recall the relation $h_K^+=h_K\,2^r/[E_K:E_K^+]$, where $r$ is the number of real places of $K$, $E_K$ is its unit group, and $E_K^+$ is the subgroup of totally positive units. In particular, $h_K^+=h_K$ if $h_K^+$ is odd or if $K$ is totally imaginary.

We write $r_2^+(K)$ for the narrow $2$-rank of $K$, namely the $2$-rank of the abelian group $Cl_K^+$. Thus $r_2^+(K)=\dim_{\mathbb{F}_2}(Cl_K^+/(Cl_K^+)^2)$, where $\mathbb{F}_2$ is the field with two elements. See \S2 for further details on the $2$-rank.

The following two theorems are our main results.

\begin{Th}\label{t1}
Let $K$ be a number field and set $Cl_{K,2}:=\{c\in Cl_K\mid c^2=1\}$ (the elementary abelian $2$-subgroup of $Cl_K$ if $h_K$ is even). Then the following hold.

$(1)$ The set $R'_{nr}(K,\mathbb{Z}/2\mathbb{Z})$ is a subgroup of $Cl_{K,2}$.

$(2)$ If $h_K^+$ is odd, then $R_{nr}(K,\mathbb{Z}/2\mathbb{Z})=\emptyset$ and $R'_{nr}(K,\mathbb{Z}/2\mathbb{Z})=\{1\}$.

$(3)$ If $K$ is of class type and $h_K^+$ is even, then
\[
R_{nr}(K,\mathbb{Z}/2\mathbb{Z})=Cl_{K,2}\setminus\{1\}
\qquad\text{and}\qquad R'_{nr}(K,\mathbb{Z}/2\mathbb{Z})=Cl_{K,2}.
\]

$(4)$ If $K$ is of unit type, then $R_{nr}(K,\mathbb{Z}/2\mathbb{Z})=R'_{nr}(K,\mathbb{Z}/2\mathbb{Z})$, and the order of $R_{nr}(K,\mathbb{Z}/2\mathbb{Z})$ is $2^{r_2^+(K)-m}$, where $2^m$ is the degree over $K$ of the compositum of all unramified quadratic extensions $K(\sqrt{\varepsilon})/K$ with $\varepsilon$ a unit of $K$.
\end{Th}

\noindent\textbf{Remarks.}

\textbf{1.} By parts $(2)$ and $(3)$ of Theorem \ref{t1}, if $K$ is of class type, then $1\notin R_{nr}(K,\mathbb{Z}/2\mathbb{Z})$. Hence $R_{nr}(K,\mathbb{Z}/2\mathbb{Z})$ is not a subgroup of $Cl_K$.

\textbf{2.} Consider part $(4)$ of Theorem \ref{t1} for a quadratic extension $K/\mathbb{Q}$ of unit type. If $K$ is imaginary, then $m=1$ because $E_K=\{\pm1\}$. The fields $\mathbb{Q}(i)=\mathbb{Q}(\sqrt{-1})$ and $\mathbb{Q}(j)=\mathbb{Q}(\sqrt{-3})$, where $j$ is a primitive cube root of unity, are of class type since their narrow class number is $1$. If $K$ is real with fundamental unit $\epsilon$, then $E_K=\{\pm\epsilon^n\mid n\in\mathbb{Z}\}$. The compositum in part $(4)$ is therefore contained in the biquadratic extension $K(\sqrt{-1},\sqrt{\epsilon})$, so $m=1$ or $m=2$.

\textbf{3.} In \S3 we classify quadratic fields over $\mathbb{Q}$ according to whether they are of class type or of unit type. We give explicit examples illustrating Theorem \ref{t1}, including fields $K$ of unit type for which $R'_{nr}(K,\mathbb{Z}/2\mathbb{Z})=Cl_{K,2}$ as in part $(3)$ (see Proposition \ref{exempletheoreme1.1unite}). $\blacksquare$

\begin{Th}\label{t2}
Let $K$ be a number field and let $\Gamma$ be a finite group. Suppose that $\Gamma$ has odd order, has a noncyclic Sylow $2$-subgroup, or has a normal cyclic Sylow $2$-subgroup. Then $R'_{nr}(K,\Gamma)$ is a subgroup of $R'_{nr}(K,\mathbb{Z}/2\mathbb{Z})$. If $\Gamma$ has even order and a cyclic Sylow $2$-subgroup, then, even without assuming normality, \[R_{nr}(K,\Gamma)\subseteq R_{nr}(K,\mathbb{Z}/2\mathbb{Z}).\]
\end{Th}

\noindent\textbf{Remark.} By Corollary \ref{cor3}(i) below, if $\Gamma$ has odd order or has a noncyclic Sylow $2$-subgroup, then $R'_{nr}(K,\Gamma)=\{1\}$. $\blacksquare$

\begin{Cor}\label{cort2}
The set $R'_{nr}(K,\Gamma)$ is a subgroup of $R'_{nr}(K,\mathbb{Z}/2\mathbb{Z})$ in the following cases.
\begin{enumerate}
\item $\Gamma$ is a nonabelian $2$-group.
\item $\Gamma=S_n$ (the symmetric group of degree $n$) or $A_n$ (the alternating subgroup of $S_n$), with $n\geq4$.
\item $\Gamma$ is a finite nilpotent group.
\end{enumerate}
\end{Cor}

\begin{proof}
\textbf{1.} The first case follows immediately.

\textbf{2.} Write $(a\ b)$ for a transposition in $S_n$. For $n\geq4$, both $S_n$ and $A_n$ contain the subgroup generated by $(1\ 2)(3\ 4)$ and $(1\ 4)(2\ 3)$. This subgroup is isomorphic to $\mathbb{Z}/2\mathbb{Z}\times\mathbb{Z}/2\mathbb{Z}$. Consequently, a Sylow $2$-subgroup $H$ containing it is noncyclic, and so are all the other Sylow $2$-subgroups, since they are isomorphic to $H$.

By the preceding remark, $R'_{nr}(K,\Gamma)=\{1\}$ in cases \textbf{1} and \textbf{2}.

\textbf{3.} The result follows immediately if $\Gamma$ has odd order or has a noncyclic Sylow $2$-subgroup. If a Sylow $2$-subgroup of $\Gamma$ is cyclic, it is normal because every Sylow subgroup of a finite nilpotent group is normal.
\end{proof}

\noindent\textbf{Remarks.}

\textbf{1.} Let $\Gamma=\mathbb{Z}/2^n\mathbb{Z}\times\Gamma'$, where $n\geq0$ and $\Gamma'$ is a finite group of odd order. Then $\Gamma$ has odd order or has a normal cyclic Sylow $2$-subgroup. It follows from Theorem \ref{t2} that $R'_{nr}(K,\Gamma)$ is a subgroup of $R'_{nr}(K,\mathbb{Z}/2\mathbb{Z})$.

\textbf{2.} If $K$ is of class type and $\Gamma$ has even order with a cyclic Sylow $2$-subgroup, then $R_{nr}(K,\Gamma)$ is not a subgroup of $Cl_{K,2}$. Indeed, the inclusion in Theorem \ref{t2} and the class-type assumption give $1\notin R_{nr}(K,\Gamma)$. In the odd-order and noncyclic cases, Corollary \ref{cor3}(i) instead gives $R_{nr}(K,\Gamma)=\{1\}$ whenever an unramified $\Gamma$-extension exists. $\blacksquare$

\section{Preliminaries}

We collect and prove results used in the proofs of Theorems \ref{t1} and \ref{t2}.

The following proposition, proved by Artin in \cite{artin}, gives a formula for the Steinitz class of a finite extension of number fields. For a generalization and a more recent proof, see \cite[Theorem I.4, p. 9]{martinet}.

\begin{Prop}[Artin's theorem]\label{Artin}
Let $L/K$ be a finite extension of number fields. Then
\[
cl_K(\mathcal O_L)=\left[\left(\frac{\delta_{L/K}}{d}\right)^{1/2}\right]_K,
\]
where $d$ is the discriminant of a $K$-basis of $L$ and $\delta_{L/K}$ is the relative discriminant.
\end{Prop}

\noindent\textbf{Remarks.}

\textbf{1.} We have $cl_K(\mathcal O_L)^2=[\delta_{L/K}]_K$. In particular, if $L/K$ is unramified, then $cl_K(\mathcal O_L)^2=1$.

\textbf{2.} If $L/K$ is unramified and Galois, then $cl_K(\mathcal O_L)$ capitulates in $Cl_L$: if $cl_K(\mathcal O_L)=[I]_K$, then $[I\mathcal O_L]_L=1$. Indeed, $cl_K(\mathcal O_L)=[(d\mathcal O_K)^{-1/2}]_K$, and $d$ is a square in $L$, as shown at the start of the next proof. $\blacksquare$

The next proposition appears in the unpublished doctoral thesis of L. P. Endo, supervised by L. R. McCulloh \cite[Corollary, p. 4]{endo}. It is central to the present paper and may be useful in other work. For the reader's convenience we give a more detailed proof that differs slightly from Endo's.

\begin{Prop}\label{prop1}
Let $L/K$ be a finite Galois extension.

$(i)$ If $Gal(L/K)$ has odd order or a noncyclic Sylow $2$-subgroup, then
\[
cl_K(\mathcal O_L)=[(\delta_{L/K})^{1/2}]_K.
\]

$(ii)$ If $Gal(L/K)$ has even order and a cyclic Sylow $2$-subgroup, then $L/K$ has a unique quadratic subextension $M/K$, and
\[
cl_K(\mathcal O_L)=\left[\left(\frac{\delta_{L/K}}{\alpha}\right)^{1/2}\right]_K,
\]
where $M=K(\sqrt\alpha)$ for $\alpha\in K$.
\end{Prop}

\begin{proof}
Let $n=[L:K]$, write $G=Gal(L/K)=\{\sigma_i\mid1\leq i\leq n\}$, and choose a $K$-basis $\mathcal B=(a_i)_{1\leq i\leq n}$ of $L$. Its discriminant for the trace pairing is $d=\det(M)^2\in K$, where $M=(\sigma_i(a_j))\in M_n(L)$. In particular, $\det(M)\in L$. We distinguish the case where every element of $G$ fixes $\det(M)$, so $K(\sqrt d)=K$, from the case where some element acts nontrivially on $\det(M)$, so $K(\sqrt d)/K$ is quadratic.

\textbf{1.} If $n$ is odd, $L/K$ cannot have a quadratic subextension. Thus $\det(M)\in K$, and Artin's theorem gives $cl_K(\mathcal O_L)=[(\delta_{L/K})^{1/2}]_K$.

\textbf{2.} For $\tau\in G$, the determinant $\tau(\det M)$ is obtained by permuting the rows $L_i=(\sigma_i(a_j))_{1\leq j\leq n}$ of $M$ according to the permutation $\tilde\tau$ induced by left multiplication $\sigma_i\mapsto\tau\sigma_i$. Thus $\tau(\det M)=\epsilon(\tilde\tau)\det M$, where $\epsilon$ is the sign character.

Under the action of $\langle\tilde\tau\rangle$, the orbit of $L_i$ consists of the rows $\tilde\tau^j(L_i)$ for $0\leq j<|\tau|$. These rows are distinct: equality of two rows implies that the corresponding automorphisms agree on the basis $\mathcal B$, so $\tau^j\sigma_i=\tau^{j'}\sigma_i$ and $j=j'$. Hence $\tilde\tau$ decomposes into $n/|\tau|$ disjoint cycles, each of length $|\tau|$ and sign $(-1)^{|\tau|-1}$. Therefore
\[
\epsilon(\tilde\tau)=(-1)^{(n/|\tau|)(|\tau|-1)}.
\]
In particular, an element $\tau$ of odd order induces a permutation of sign $1$.

\textbf{3.} Suppose that $n=2^mk$ with $m\geq1$ and $k$ odd. If $\tau$ has even order $|\tau|=2^{m'}k'$ with $1\leq m'\leq m$ and $k'\mid k$, then $\epsilon(\tilde\tau)=(-1)^{2^{m-m'}}$.

\textbf{3.1.} Suppose the Sylow $2$-subgroups of $G$ are noncyclic. If $m'=m$, then $\langle\tau\rangle$ would contain a cyclic Sylow $2$-subgroup of $G$. Thus $m'<m$ and $\epsilon(\tilde\tau)=1$. It follows that $\det(M)\in K$, and Artin's theorem gives the remaining assertion in $(i)$.

\textbf{3.2.} Suppose that $syl_2(G)$ is cyclic and choose a generator $\tau$. Since $|\tau|=2^m$, we have $\epsilon(\tilde\tau)=-1$, so $K(\sqrt d)/K$ is quadratic. If there were a distinct quadratic subextension $K(\sqrt\alpha)/K$ of $L/K$, their compositum would be biquadratic with Galois group $G/H$ for some normal subgroup $H$ of $G$. Its Sylow $2$-subgroup would be the cyclic group $syl_2(G)H/H$, by Proposition \ref{prp3}(1.2). This is impossible for a biquadratic extension. Thus $K(\sqrt\alpha)=K(\sqrt d)$ whenever $K(\sqrt\alpha)$ is a quadratic subfield of $L$. In that case $d=a^2\alpha$ for some $a\in K^\times$, and substituting in Artin's formula proves $(ii)$.
\end{proof}

\begin{Cor}\label{cor3}
Let $L/K$ be a finite unramified Galois extension.
\begin{itemize}
\item[(i)] If $Gal(L/K)$ has odd order or a noncyclic Sylow $2$-subgroup, then the Steinitz class is trivial.
\item[(ii)] If $Gal(L/K)$ has even order and a cyclic Sylow $2$-subgroup, then the Steinitz class is $cl_K(\mathcal O_M)$, where $M/K$ is the unique quadratic subextension of $L/K$.
\end{itemize}
\end{Cor}

Let $K$ be a number field. We recall the ray congruence $\bmod^*$ of class field theory for a modulus $\mathfrak M$ given by a nonzero proper integral ideal of $\mathcal O_K$. We need only the finite part of the modulus. Write $K^\times=K\setminus\{0\}$ and let $v_{\mathfrak p}$ denote the $\mathfrak p$-adic valuation. The set
\[
K^\times_{\mathfrak M}=
\{x\in K^\times\mid v_{\mathfrak p}(x-1)\geq v_{\mathfrak p}(\mathfrak M)
\text{ for every }\mathfrak p\mid\mathfrak M\}
\]
is a subgroup of $K^\times$. Its cosets define a congruence compatible with multiplication. In the usual notation, $x\equiv y\bmod^*\mathfrak M$ means $xy^{-1}\in K^\times_{\mathfrak M}$. In particular,
\begin{equation}\label{1}
x\equiv y\bmod^*\mathfrak M
\quad\Longleftrightarrow\quad
x\equiv y\bmod^*\mathfrak p^{v_{\mathfrak p}(\mathfrak M)}
\text{ for every }\mathfrak p\mid\mathfrak M,
\end{equation}
and
\begin{equation}\label{2}
x\equiv y\bmod^*\mathfrak M
\quad\Longrightarrow\quad
v_{\mathfrak p}(x)=v_{\mathfrak p}(y)
\text{ for every }\mathfrak p\mid\mathfrak M.
\end{equation}

For $x,y\in\mathcal O_K$, the ordinary congruence $x\equiv y\bmod\mathfrak M$ means $x-y\in\mathfrak M$, equivalently $v_{\mathfrak p}(x-y)\geq v_{\mathfrak p}(\mathfrak M)$ for every $\mathfrak p\mid\mathfrak M$. The units of $\mathcal O_K/\mathfrak M$ are precisely the classes of the elements $x$ with $v_{\mathfrak p}(x)=0$ for all such primes. If $x$ and $y$ are units modulo $\mathfrak M$, then
\begin{equation}\label{3}
x\equiv y\bmod^*\mathfrak M\quad\Longleftrightarrow\quad
x\equiv y\bmod\mathfrak M,
\end{equation}
since $v_{\mathfrak p}(x/y-1)=v_{\mathfrak p}(x-y)-v_{\mathfrak p}(y)=v_{\mathfrak p}(x-y)$.

Let $x\ne0$ belong to the localization $\mathcal O_{K,\mathfrak p}$. For every positive integer $n$,
\begin{equation}\label{4}
x\equiv a\bmod^*\mathfrak p^n\quad\text{for some }a\in\mathcal O_K.
\end{equation}
Indeed, write $x=b/c$ with $b,c\in\mathcal O_K$ and $c\notin\mathfrak p$. Let $d$ represent the inverse of $c$ modulo $\mathfrak p^n$ and take $a=bd$. Then $a/x-1=cd-1\in\mathfrak p^n$, which gives the congruence.

\begin{Prop}\label{unitemodulocarre}
Let $K$ be a number field, $\varepsilon$ a unit of $K$, and $\mathfrak M$ a finite modulus of $K$. The following are equivalent.

(i) There exists $x\in K^\times$ such that $\varepsilon\equiv x^2\bmod^*\mathfrak M$.

(ii) There exists $a\in\mathcal O_K$ such that $\varepsilon\equiv a^2\bmod\mathfrak M$.
\end{Prop}

\begin{proof}
By (\ref{3}), in (ii) one may replace $\bmod$ by $\bmod^*$, since the classes of $\varepsilon$ and $a^2$ are units modulo $\mathfrak M$. This gives (i).

Conversely, assume (i). For each $\mathfrak p\mid\mathfrak M$, (\ref{2}) gives $v_{\mathfrak p}(x^2)=v_{\mathfrak p}(\varepsilon)=0$. Thus $x\in\mathcal O_{K,\mathfrak p}$, and (\ref{4}) gives $a_{\mathfrak p}\in\mathcal O_K$ with $x\equiv a_{\mathfrak p}\bmod^*\mathfrak p^{v_{\mathfrak p}(\mathfrak M)}$. The Chinese remainder theorem gives $a\in\mathcal O_K$ with $a\equiv a_{\mathfrak p}\bmod\mathfrak p^{v_{\mathfrak p}(\mathfrak M)}$ for each $\mathfrak p\mid\mathfrak M$. Since $a_{\mathfrak p}$ is a unit at $\mathfrak p$, (\ref{3}) gives the same congruence with $\bmod^*$. By (\ref{1}), $x\equiv a\bmod^*\mathfrak M$, so $\varepsilon\equiv a^2\bmod^*\mathfrak M$. Applying (\ref{3}) again gives (ii).
\end{proof}

The following proposition follows from Kummer theory \cite[\S39, pp. 132--139]{hecke} and Artin's theorem. It is the case $p=2$ of \cite[Proposition 3.2]{BrS}.

\begin{Prop}\label{th2}
Let $L/K$ be quadratic and write $L=K(\sqrt\alpha)$ with $\alpha\in K^\times$.
\begin{enumerate}
\item The extension $L/K$ is tamely ramified if and only if there exists $x\in K^\times$ such that $\alpha\equiv x^2\bmod^*4\mathcal O_K$.
\item Suppose $L/K$ is tamely ramified. Then
\begin{enumerate}
\item There is a unique factorization $\alpha\mathcal O_K=I(\alpha)^2J$, with $I(\alpha)$ a fractional ideal and $J$ a squarefree integral ideal of $\mathcal O_K$.
\item $\delta_{L/K}=J$.
\item $cl_K(\mathcal O_L)=[I(\alpha)]_K^{-1}$.
\end{enumerate}
\end{enumerate}
\end{Prop}

\begin{proof}
We briefly prove assertion 1. The quadratic extension $L/K$ is tamely ramified if and only if the primes above $2$ are unramified. If $\alpha$ has valuation zero at all primes over $2$, the assertion is Hecke's Theorem 119 \cite[p. 136]{hecke}. For a general $\alpha$, if $\alpha\equiv x^2\bmod^*4\mathcal O_K$, then $\alpha'=\alpha x^{-2}$ has valuation zero at these primes and satisfies $\alpha'\equiv1\bmod^*4\mathcal O_K$. Apply the theorem to $\alpha'$, which defines the same extension. Conversely, if the primes above $2$ are unramified, the valuations of $\alpha$ there are even, since $\sqrt\alpha$ lies in an extension of ramification index $1$. By weak approximation, choose $y\in K^\times$ with $v_{\mathfrak p}(y)=v_{\mathfrak p}(\alpha)/2$ at these primes. The theorem applied to $\alpha y^{-2}$ gives the required congruence for $\alpha$.
\end{proof}

Since every unramified extension is tame, we obtain the following corollary.

\begin{Cor}\label{cor1}
Let $L=K(\sqrt\alpha)$ be quadratic over $K$. Then $L/K$ is unramified if and only if there exist $x\in K^\times$ and a fractional ideal $I$ of $K$ such that $\alpha\equiv x^2\bmod^*4\mathcal O_K$ and $\alpha\mathcal O_K=I^2$. In that case $cl_K(\mathcal O_L)=[I]_K^{-1}=[I]_K$.
\end{Cor}

\noindent\textbf{Remark.} The first assertion is the quadratic case of Hecke's Theorem 120 \cite[p. 137]{hecke}. $\blacksquare$

\begin{Cor}\label{cor4}
Let $N/K$ be an unramified biquadratic extension. Write $L_1,L_2,L_3$ for its quadratic subfields over $K$. Then
\[
cl_K(\mathcal O_{L_3})=cl_K(\mathcal O_{L_1})cl_K(\mathcal O_{L_2}).
\]
\end{Cor}

\begin{proof}
This follows immediately from \cite[Proposition 3.2(1)]{bouchaib}.
\end{proof}

For a prime $p$ and a finite group $G$, write $Syl_p(G)$ for the set of its Sylow $p$-subgroups and $syl_p(G)$ for a chosen member of this set. Thus $Syl_p(G)=\{g\,syl_p(G)g^{-1}\mid g\in G\}$. When convenient we use $syl_p(G)$ for any conjugate of the chosen subgroup. We write $H\leq G$ for a subgroup and $H\trianglelefteq G$ for a normal subgroup. The order of a group $G$ or an element $g$ is denoted by $|G|$ or $|g|$, respectively.

We will use the following standard facts.

\begin{Prop}\label{prp3}
Let $G$ and $G'$ be finite groups.

$(1)$ If $H\trianglelefteq G$, then

$(1.1)$ $Syl_p(H)=\{P\cap H\mid P\in Syl_p(G)\}$.

$(1.2)$ $Syl_p(G/H)=\{PH/H\mid P\in Syl_p(G)\}$.

$(2)$ $Syl_p(G\times G')=Syl_p(G)\times Syl_p(G')$, under the natural identification of pairs with product subgroups.

$(3)$ If $\gcd(|G|,|G'|)=1$, then $K\leq G\times G'$ if and only if $K=H\times H'$ for subgroups $H\leq G$ and $H'\leq G'$.

$(4)$ If $H,K\trianglelefteq G$, then $G$ is the internal direct product of $H$ and $K$ if and only if it is their internal semidirect product.
\end{Prop}

\begin{proof}
For (1.1), (1.2), and (2), see \cite[Corollary 6.4, Theorem 6.6, and Theorem 7.1]{CO}. For (3), see \cite[Theorem 2]{AN}. In (4) the two normal factors commute since their commutators lie in their trivial intersection.
\end{proof}

\begin{Prop}\label{Scott}
Let $G$ have even order $2^nk$, where $k$ is odd. If $G$ has a cyclic Sylow $2$-subgroup, then it has a unique normal subgroup $H$ of order $k$. The subgroup $H$ consists of the elements of odd order in $G$.
\end{Prop}

\begin{proof}
See \cite[Theorem 6.8, p. 14]{CO} or \cite[6.2.11, p. 138]{SC}.
\end{proof}

\noindent\textbf{Remarks.}

\textbf{1.} With the notation of Proposition \ref{Scott}, $G=H\rtimes syl_2(G)$, so $G/H\simeq syl_2(G)$. If $syl_2(G)$ is normal, Proposition \ref{prp3}(4) makes this a direct product. Consequently, $syl_2(G)$ is normal if and only if $G\simeq\mathbb Z/2^n\mathbb Z\times G'$ for some group $G'$ of order $k$.

\textbf{2.} A finite group $G$ is $p$-nilpotent if it has a normal subgroup $N$ such that $G=N\rtimes syl_p(G)$. By the preceding remark, the group in Proposition \ref{Scott} is $2$-nilpotent, a fact first proved by Burnside. It is also solvable, since $H$ is solvable by the Feit--Thompson theorem and $G/H$ is abelian. By a theorem of Shafarevich, $G$ can be realized over any number field $K$ as the Galois group of a tamely ramified extension. $\blacksquare$
\smallskip

We use the definition of $p$-rank in \cite[Chap. IX, pp. 232--233]{CAFR}. It extends the usual definition for abelian groups and is convenient for Galois theory.

\begin{Df}\label{rang}
Let $G$ be a finite group. Its $p$-rank $r_p(G)$ is the largest integer $n\geq0$ for which $G$ has a normal subgroup $H$ with $G/H$ elementary abelian of order $p^n$. In this case $n=\dim_{\mathbb F_p}(G/H)$, where $\mathbb F_p$ is the field of $p$ elements.
\end{Df}

\noindent\textbf{Remark.} In group theory, $p$-rank can also mean the largest $n$ for which $G$ has an elementary abelian subgroup of order $p^n$. For abelian $G$, the two definitions agree, and $r_p(G)$ is the rank of $G_p=\{g\in G\mid g^p=1\}$. $\blacksquare$

\begin{Prop}\label{2rang}
Let $G$ be a finite group, write $G^2=\{g^2\mid g\in G\}$, and let $\langle G^2\rangle$ be the subgroup it generates.
 
(1) The subgroup $\langle G^2\rangle$ is normal, $G/\langle G^2\rangle$ is elementary abelian, and $r_2(G)=\dim_{\mathbb F_2}(G/\langle G^2\rangle)$. If $G$ is abelian, then $\langle G^2\rangle=G^2$.

(2) The number of subgroups of $G$ of index $2$ is $2^{r_2(G)}-1$.

(3) For any finite group $G'$, we have $r_2(G\times G')=r_2(G)+r_2(G')$.

(4) If $G$ has even order and a cyclic Sylow $2$-subgroup, then $r_2(G)=1$.

(5) If $H\trianglelefteq G$, then $r_2(G/H)\leq r_2(G)$.
\end{Prop}

\begin{proof}
(1) Since $gh^2g^{-1}=(ghg^{-1})^2$, the subgroup generated by the squares is normal. Every element of the quotient has order $1$ or $2$, so the quotient is elementary abelian. If $G/H$ is any elementary abelian $2$-quotient, every square belongs to $H$, and therefore $G/H$ is a quotient of $G/\langle G^2\rangle$. This proves the dimension formula.

(2) Every subgroup of index $2$ contains $\langle G^2\rangle$ and corresponds to an index-$2$ subgroup of $G/\langle G^2\rangle$. These are the hyperplanes of an $\mathbb F_2$-space of dimension $r_2(G)$, of which there are $2^{r_2(G)}-1$.

(3) The squares in $G\times G'$ are the pairs in $G^2\times (G')^2$. They include the squares supported in either factor, so
\[
\langle(G\times G')^2\rangle
=\langle G^2\rangle\times\langle(G')^2\rangle.
\]
The corresponding quotient is a direct product, so the dimensions add.

(4) Every elementary abelian $2$-quotient of $G$ is cyclic by Proposition \ref{prp3}(1.2), so $r_2(G)\leq1$. In the notation of Proposition \ref{Scott}, $G/H$ is cyclic of order $2^n$ with $n\geq1$, and hence has a quotient of order $2$. Thus $r_2(G)=1$.

(5) Any elementary abelian quotient of $G/H$ is also a quotient of $G$, so its rank is at most $r_2(G)$.
\end{proof}

Throughout we write $r_2(K)=r_2(Cl_K)$ and $r_2^+(K)=r_2(Cl_K^+)$. These are the ordinary and narrow $2$-ranks of $K$. Set $Cl^+_{K,2}=\{c\in Cl_K^+\mid c^2=1\}$ and recall that $Cl_{K,2}=\{c\in Cl_K\mid c^2=1\}$. The remark after Definition \ref{rang} gives $r_2(K)=r_2(Cl_{K,2})$ and $r_2^+(K)=r_2(Cl^+_{K,2})$.

Let $k$ be a number field and $K/k$ a quadratic extension. Write $u$ and $u_\infty$ for the numbers of finite and infinite places of $k$ that ramify in $K$. Let $E_k$ and $E_k^+$ be the units and totally positive units of $k$, respectively. Define $e$ and $e^+$ by
\[
2^e=[E_k:E_k\cap N_{K/k}(K^\times)],\qquad
2^{e^+}=[E_k^+:E_k^+\cap N_{K/k}(K^\times)],
\]
where $N_{K/k}$ is the norm. Write $Am(K/k)$ and $Am^+(K/k)$ for the ambiguous class groups, namely the subgroups of $Cl_K$ and $Cl_K^+$ fixed by $Gal(K/k)$, respectively.

\begin{Prop}\label{K/kquadratique}
With the preceding notation:

(i) If $h_k$ is odd, then $syl_2(Am(K/k))=Cl_{K,2}$ and $r_2(K)=u+u_\infty-e-1$.

(ii) If $h_k^+$ is odd, then $syl_2(Am^+(K/k))=Cl^+_{K,2}$ and $r_2^+(K)=u-e^+-1$.

(iii) If $h_k=1$, then $Am(K/k)=Cl_{K,2}$. If $h_k^+=1$, then $Am^+(K/k)=Cl^+_{K,2}$.
\end{Prop}

\begin{proof}
(i) If $c\in Cl_K$ has $2$-power order, then $N_{K/k}(c)$ also has $2$-power order. Since $h_k$ is odd, this norm is trivial. Extending ideals from $\mathcal O_k$ to $\mathcal O_K$ gives $c\sigma(c)=1$, where $\sigma$ generates $Gal(K/k)$. Thus $\sigma(c)=c$ if and only if $c^2=1$, giving $syl_2(Am(K/k))=Cl_{K,2}$. Chevalley's ambiguous class number formula \cite[a, p. 406]{CH} gives $|Am(K/k)|=h_k2^{u+u_\infty-e-1}$. Taking the $2$-part gives the stated rank formula (see also \cite[(1.1), p. 27]{SORO}).

(ii) The argument for $syl_2(Am^+(K/k))=Cl^+_{K,2}$ is analogous. Gras's narrow ambiguous class number formula \cite[Theorem 4.1, p. 26]{gras} gives $|Am^+(K/k)|=h_k^+2^{u-e^+-1}$. Taking the $2$-part yields the formula for $r_2^+(K)$.

(iii) In either class group, the norm of every class is trivial under the corresponding class-number hypothesis. Thus $\sigma(c)=c^{-1}$ for every class $c$, and being fixed is equivalent to $c^2=1$.
\end{proof}

Taking $k=\mathbb Q$ gives the following classical results of Hilbert and Hasse.

\begin{Cor}\label{corK/kquadratique}
Let $K/\mathbb Q$ be quadratic.

(1) $r_2^+(K)=u-1$.

(2) If $K$ is imaginary, then $r_2(K)=u-1$.

(3) Suppose $K$ is real and let $\epsilon$ be its fundamental unit.

(a) If $N_{K/\mathbb Q}(\epsilon)=-1$, then $r_2(K)=u-1$.

(b) Suppose $N_{K/\mathbb Q}(\epsilon)=1$.

(b.1) If $-1$ is a norm from $K$, then $r_2(K)=u-1$.

(b.2) If $-1$ is not a norm from $K$, then $r_2(K)=u-2$.
\end{Cor}

\noindent\textbf{Remark.} Since $h_K^+$ is odd exactly when $r_2^+(K)=0$, and $h_K$ is odd exactly when $r_2(K)=0$, the corollary identifies the quadratic fields with odd ordinary or narrow class number. $\blacksquare$
\smallskip

Hasse's norm theorem makes the conditions in (b.1) and (b.2) explicit. The element $-1$ is a norm from $K$ if and only if it is a local norm at every finite and infinite place. If $K=\mathbb Q(\sqrt d)$ with $d$ positive and squarefree, this is equivalent to $\left(\frac{-1}{p}\right)=1$ for every odd prime $p\mid d$. Here $\left(\frac{\cdot}{p}\right)$ is the Legendre symbol. Since $\left(\frac{-1}{p}\right)=(-1)^{(p-1)/2}$, this means that all odd prime divisors of $d$ are congruent to $1$ modulo $4$. Equivalently, $d$ is a sum of two squares, or $K$ can be embedded in a cyclic quartic extension of $\mathbb Q$.

By \cite[I, p. 94]{FRMA}, for a number field with $r$ real embeddings, $r_2^+(K)-r_2(K)\leq\lfloor r/2\rfloor$. For a real quadratic field this difference is $0$ or $1$, and the preceding corollary shows that it is $1$ precisely in case (b.2). The element $-1$ fails to be a norm exactly when some prime divisor $p$ of $d$ satisfies $p\equiv3\bmod4$. This condition also implies $N_{K/\mathbb Q}(\epsilon)=1$ by \cite[(1.7), p. 179]{FT}.

\begin{Cor}\label{r^+_2(K)= r_2(K)+1}
Let $K=\mathbb Q(\sqrt d)$ be real quadratic with $d$ positive and squarefree. Then $r_2^+(K)=r_2(K)+1$ if and only if $d$ has a prime divisor $p\equiv3\bmod4$.
\end{Cor}

\noindent\textbf{Remark.} Another proof is given in \cite[Remark 2]{Dummitquebec}. $\blacksquare$
\smallskip

We now describe $Cl_{K,2}$ explicitly for a quadratic field $K$. Following Hilbert \cite[\S75, p. 136]{HI}, an integral ideal is ambiguous in Hilbert's sense if it is fixed by $Gal(K/\mathbb Q)$ and divisible by no rational integer other than $\pm1$. By Hilbert's Theorem 105 \cite[p. 136]{HI}, there are exactly $u$ ambiguous prime ideals, one above each ramified rational prime. Their squarefree products give exactly $2^u$ ambiguous ideals. Let $Am_0(K/\mathbb Q)$ be the subgroup of $Am(K/\mathbb Q)$ consisting of classes represented by such ideals. Hilbert's Theorem 106 \cite[p. 136]{HI} gives
\[
|Am_0(K/\mathbb Q)|=
\begin{cases}
2^{u-1},&\text{if }K\text{ is imaginary or real with }N_{K/\mathbb Q}(\epsilon)=-1,\\
2^{u-2},&\text{if }K\text{ is real with }N_{K/\mathbb Q}(\epsilon)=1.
\end{cases}
\]

\noindent\textbf{Remarks.}

\textbf{1.} The subgroup $Am_0(K/\mathbb Q)$ is generated by the classes of the prime ideals above the ramified rational primes.

\textbf{2.} In the general sense, an ambiguous fractional ideal is any fractional ideal fixed by $Gal(K/\mathbb Q)$. The classes of all such ideals still give $Am_0(K/\mathbb Q)$. Indeed, every such class has an ambiguous integral representative $\mathfrak a$, which by \cite[(1.16), p. 182]{FT} can be written $r\mathfrak q$ with $r\in\mathbb Q$ and $\mathfrak q$ ambiguous in Hilbert's sense. $\blacksquare$
\smallskip

Combining Hilbert's Theorem 106 with Corollary \ref{corK/kquadratique} gives the following description.

\begin{Cor}\label{descclasseK2}
Let $K/\mathbb Q$ be quadratic.

If $K$ is imaginary or real with $N_{K/\mathbb Q}(\epsilon)=-1$, then
\[
Cl_{K,2}=Am_0(K/\mathbb Q),\qquad |Cl_{K,2}|=2^{u-1}.
\]
Suppose that $K$ is real with $N_{K/\mathbb Q}(\epsilon)=1$. If $-1$ is not a norm from $K$, then
\[
Cl_{K,2}=Am_0(K/\mathbb Q),\qquad |Cl_{K,2}|=2^{u-2}.
\]
If $-1$ is a norm from $K$, then
\[
Cl_{K,2}=Am(K/\mathbb Q),\qquad |Cl_{K,2}|=2^{u-1},
\qquad [Am(K/\mathbb Q):Am_0(K/\mathbb Q)]=2.
\]
\end{Cor}

\noindent\textbf{Remark.} Corollary 1 of \cite[p. 181]{FT} requires a small correction: replace $C_{K,2}$ by $Am_0(K/\mathbb Q)$ in the notation and statement. The symbol $C_{K,2}$ is undefined there, but the definition of $C^+_{K,2}$ in \cite[Theorem 39, p. 181]{FT} suggests that it denotes our $Cl_{K,2}$. With this correction, the corollary is Hilbert's Theorem 106. Also, in the definition of $\nu_K$ on p. 181, replace $C_K$ by $C_K^+$. $\blacksquare$

\section{Steinitz Classes of Unramified Quadratic Extensions with Galois Group $\Gamma\simeq\mathbb{Z}/2\mathbb{Z}$}

Let $l$ be a prime and let $K$ be a number field containing a primitive $l$th root of unity. Set
\[
V_l(K)=\{\alpha\mid I^l=\alpha\mathcal{O}_K
\text{ for some fractional ideal } I \text{ of } K\}.
\]
In \cite[Definition 5.2.4(3), p. 231]{cohan}, Cohen defines the $l$-Selmer group of $K$, denoted by $S_l(K)$, by analogy with the Selmer groups of elliptic curves:
\[
S_l(K)=V_l(K)/(K^\times)^l.
\]
The elements of $V_l(K)$ are called virtual units, which explains the notation $V$ (see \cite[Remark, p. 233]{cohan}). For work relating $S_2(K)$ to the narrow class group rank of a number field and to quadratic reciprocity, see \cite{Dummitlondon,LE}.

\noindent\textbf{Remark.} By Kummer theory, more precisely \cite[Theorem 118, p. 134]{hecke}, the extensions of degree $l$ unramified at all finite places away from $l$ correspond to the one-dimensional subspaces of $S_l(K)$, regarded as an $\mathbb{F}_l$-vector space. The line generated by $\bar\alpha$ corresponds to $K(\sqrt[l]{\alpha})$. In particular, for $l=2$ this gives a bijection with $S_2(K)\setminus\{\bar1\}$. $\blacksquare$
\smallskip

Let $K$ be any number field and let $R_K$ be the subgroup of $K^\times$ consisting of the elements $\alpha$ satisfying both of the following conditions.
\begin{enumerate}
\item There exists $x\in K^\times$ such that $\alpha\equiv x^2\bmod^*4\mathcal{O}_K$.
\item There exists a fractional ideal $I$ of $K$ such that $I^2=\alpha\mathcal{O}_K$.
\end{enumerate}
Since $R_K$ is a subgroup of $V_2(K)$ containing $(K^\times)^2$, the quotient $\bar R_K:=R_K/(K^\times)^2$ is a subgroup of the $2$-Selmer group $S_2(K)$.

If $L/K$ is an unramified quadratic extension, Proposition \ref{th2} and Corollary \ref{cor1} give a unique $\bar\alpha\in\bar R_K\setminus\{\bar1\}$ with $L=K(\sqrt\alpha)$. This gives a bijection between $\bar R_K\setminus\{\bar1\}$ and the unramified quadratic extensions of $K$. The study of their Steinitz classes thus provides another way to investigate the $2$-Selmer group of a number field.

Write $U_K=E_K\cap R_K$, where $E_K$ is the unit group of $K$. Then
\[
U_K=\{\varepsilon\in E_K\mid
\varepsilon\equiv x^2\bmod^*4\mathcal{O}_K\text{ for some }x\in K^\times\}.
\]
Proposition \ref{unitemodulocarre} gives the simpler description
\begin{equation}\label{U_Kextension2}
U_K=\{\varepsilon\in E_K\mid
\varepsilon\equiv a^2\bmod4\mathcal{O}_K\text{ for some }a\in\mathcal{O}_K\}.
\end{equation}
We have $E_K^2\subset U_K$. For $\varepsilon\in E_K$,
\begin{equation}\label{U_Kramifiee}
K(\sqrt\varepsilon)/K\text{ is quadratic and unramified}
\quad\Longleftrightarrow\quad\varepsilon\in U_K\setminus E_K^2.
\end{equation}

\noindent\textbf{Remark.} In the notation of \cite[Definition, p. 226]{Dummitquebec}, the assertion following the definition of $E_{K,4}$ identifies it, in effect, with our $U_K\setminus E_K^2$, using \cite[Lemma 1, p. 226]{Dummitquebec}. Justifying that identification requires (\ref{U_Kextension2}), even if one invokes \cite[Proposition 4.8]{Dummitlondon}. The $\bmod^*$ congruence also gives an immediate answer to the questions raised in \cite[Remark 4, p. 226]{Dummitquebec}. $\blacksquare$

\begin{Df}\label{Dftypeclasse}
Following Hasse \cite[\S22]{hass}, we say that $K$ is of class type if $U_K=E_K^2$, and of unit type otherwise.
\end{Df}

By (\ref{U_Kramifiee}), $K$ is of class type if and only if every quadratic extension $K(\sqrt\varepsilon)/K$, with $\varepsilon$ a unit of $K$, ramifies at some prime ideal $\mathfrak p$. We used this characterization as the definition in \S1 to simplify the presentation. Such a $\mathfrak p$ necessarily lies above $2$: the basis $\mathcal B=(1,\sqrt\varepsilon)$ of the $K$-vector space $K(\sqrt\varepsilon)$ consists of integral elements and has discriminant ideal $4\mathcal{O}_K$, which is divisible by the relative discriminant $\delta_{K(\sqrt\varepsilon)/K}$.
\smallskip

\noindent\textbf{Remark.} If $h_K^+$ is odd, then $h_K^+=h_K$ and $K$ is of class type. Examples include $\mathbb{Q}$ and every totally imaginary number field with odd class number. For the latter, recall that a Galois extension of $\mathbb{Q}$ is either totally imaginary or totally real. $\blacksquare$
\smallskip

Let $K$ have $r\geq1$ real embeddings $\sigma_i$. The signature map $s:K^\times\to\{1,-1\}^r$ is given by \[s(x)=(\operatorname{sgn}(\sigma_1(x)),\ldots,\operatorname{sgn}(\sigma_r(x))),\] where $\operatorname{sgn}(\sigma_i(x))=\sigma_i(x)/|\sigma_i(x)|$.

\begin{Prop}\label{totalement reel}
Let $K$ be totally real. Then $K$ is of class type if and only if the restriction of $s$ to $E_K$ is surjective.
\end{Prop}

\begin{proof}
Consider the reduction map
\[
\pi:E_K/E_K^2\longrightarrow
(\mathcal O_K/4\mathcal O_K)^\times/
\bigl((\mathcal O_K/4\mathcal O_K)^\times\bigr)^2.
\]
By \cite[Theorem 3.5, $(c)\Leftrightarrow(a)$]{muller}, $\pi$ is injective if and only if the restriction of $s$ to $E_K$ is surjective. The kernel of $\pi$ consists of the classes of units that are squares modulo $4\mathcal O_K$. By (\ref{U_Kextension2}) it equals $U_K/E_K^2$, so injectivity is equivalent to $K$ being of class type.
\end{proof}

\noindent\textbf{Remark.} This proposition is Theorem 3.1.1 of \cite[p. 23]{Derhem}. We proved Proposition \ref{unitemodulocarre} to supply a step missing from that proof. $\blacksquare$
\smallskip

The next two propositions characterize the quadratic extensions of $\mathbb{Q}$ of class type and those of unit type.

Let $K/\mathbb{Q}$ be imaginary quadratic. The fields $\mathbb{Q}(i)=\mathbb{Q}(\sqrt{-1})$ and $\mathbb{Q}(j)=\mathbb{Q}(\sqrt{-3})$, where $j$ is a primitive cube root of unity, are of class type since their class number is $1$. For every other imaginary quadratic field, $E_K=\{1,-1\}$, so $K(\sqrt{-1})/K$ is the only quadratic extension obtained by adjoining the square root of a unit. Thus $K$ is of class type if and only if this extension is ramified.

\begin{Prop}\label{imaginaire}
Let $K=\mathbb{Q}(\sqrt d)$ be imaginary quadratic, with $d$ a squarefree integer.

(i) The field $K$ is of unit type if and only if $d\ne-1$ and $d\equiv3\bmod4$.

(ii) The field $K$ is of class type if and only if $d=-1$ or $d\equiv1$ or $2\bmod4$.
\end{Prop}

\begin{proof}
The two types are complementary, so it suffices to prove the reverse implications. We have already treated $d=-1$ and assume henceforth that $d\ne-1$.

(i) If $d\equiv3\bmod4$, then $(\sqrt d)^2=d\equiv-1\bmod4\mathcal{O}_K$. By (\ref{U_Kextension2}) and (\ref{U_Kramifiee}), $K(\sqrt{-1})/K$ is unramified, so $K$ is of unit type.

(ii) Consider the biquadratic compositum $K(i)$ of $K$ and $\mathbb{Q}(i)$. Since $-1$ is a unit, only primes above $2$ can ramify in $K(i)/K$.

If $d\equiv1\bmod4$, then $2$ is unramified in $K/\mathbb{Q}$ and ramified in $\mathbb{Q}(i)/\mathbb{Q}$. It therefore has ramification index $2$ in $K(i)/\mathbb{Q}$, since total ramification would force ramification in $K/\mathbb{Q}$. Thus a prime above $2$ ramifies in $K(i)/K$, and $K$ is of class type.

If $d\equiv2\bmod4$, then \cite[Theorem 3(iii)]{Wil} gives $\delta_{K(i)/\mathbb{Q}}=64d^2\mathbb{Z}$. In the notation of that theorem, $m=d$, $n=-1$, $l=\gcd(m,n)=1$, $m_1=d$, and $m_2=-1$, with $m\equiv2\bmod4$ and $-1\equiv3\bmod4$. The discriminant formula in a tower gives
\[
\delta_{K(i)/\mathbb{Q}}
=\delta_{K/\mathbb{Q}}^2 N_{K/\mathbb{Q}}(\delta_{K(i)/K}).
\]
Since $\delta_{K/\mathbb{Q}}=4d\mathbb{Z}$, we obtain $N_{K/\mathbb{Q}}(\delta_{K(i)/K})=4\mathbb{Z}$. Hence $K(i)/K$ is ramified and $K$ is of class type.
\end{proof}

\begin{Prop}\label{reel}
Let $K=\mathbb{Q}(\sqrt d)$ be real quadratic with $d$ a positive squarefree integer, and let $\epsilon$ be its fundamental unit.

(i) The field $K$ is of unit type if and only if $N_{K/\mathbb{Q}}(\epsilon)=1$.

(ii) The field $K$ is of class type if and only if $N_{K/\mathbb{Q}}(\epsilon)=-1$.
\end{Prop}

\begin{proof}
Let $s_{E_K}$ be the restriction of the signature map $s:K^\times\to\{1,-1\}^2$ to $E_K$. If $N_{K/\mathbb{Q}}(\epsilon)=-1$, then the image of $s_{E_K}$ is $\{1,-1\}^2$. If the norm is $1$, the image is $\{(1,1),(-1,-1)\}$. Since $K$ is totally real, Proposition \ref{totalement reel} gives the result.
\end{proof}

\begin{Prop}\label{unité}
Let $L=K(\sqrt\alpha)$ be an unramified quadratic extension of $K$. There exists a unit $\varepsilon$ of $K$ with $L=K(\sqrt\varepsilon)$ if and only if $cl_K(\mathcal O_L)=1$.
\end{Prop}

\begin{proof}
If $L=K(\sqrt\varepsilon)$ with $\varepsilon$ a unit, Artin's theorem gives \[cl_K(\mathcal O_L)=[(\mathcal O_K/(4\varepsilon))^{1/2}]_K=1.\] Conversely, suppose $cl_K(\mathcal O_L)=1$. By Corollary \ref{cor1}, there is a fractional ideal $I$ of $K$ with $I^2=\alpha\mathcal O_K$ and $[I]_K=cl_K(\mathcal O_L)$. Thus $I$ is principal, say $I=\beta\mathcal O_K$ with $\beta\in K^\times$, and $\alpha=\beta^2\varepsilon$ for some unit $\varepsilon$. It follows that $K(\sqrt\alpha)=K(\sqrt\varepsilon)$.
\end{proof}

\noindent\textbf{Remark.} This result is also \cite[Proposition 4]{Dummitquebec}, where the proof gives an explicit basis of the $\mathcal O_K$-module $\mathcal O_L$. $\blacksquare$
\smallskip

The preceding proposition gives a characterization of fields of class type in terms of Steinitz classes.

\begin{Cor}\label{cor2}
A number field $K$ is of class type if and only if every unramified quadratic extension $L/K$ has $cl_K(\mathcal O_L)\ne1$. Equivalently, $\mathcal O_L$ is nonfree as an $\mathcal O_K$-module, or $cl_K(\mathcal O_L)$ has order $2$.
\end{Cor}

\noindent\textbf{Remark.} The equivalence with nonfreeness follows from the definition of the Steinitz class. The equivalence with having order $2$ follows from $cl_K(\mathcal O_L)^2=1$. $\blacksquare$

\begin{proof}[\textbf{Proof of Theorem \ref{t1}}]
(1) By the first remark following Proposition \ref{Artin}, $R'_{nr}(K,\mathbb{Z}/2\mathbb{Z})\subset Cl_{K,2}$. Since this set is finite and nonempty, it suffices to prove closure under multiplication. Let $c_1,c_2$ belong to the set. If either is trivial or if $c_1=c_2$, their product belongs to the set. Otherwise, choose unramified quadratic extensions $L_i/K$ with $cl_K(\mathcal O_{L_i})=c_i$ for $i=1,2$. Their compositum is an unramified biquadratic extension. Its third quadratic subextension $L_3/K$ satisfies $cl_K(\mathcal O_{L_3})=c_1c_2$ by Corollary \ref{cor4}. This proves closure.

(2) If $h_K^+$ is odd, we already know that \[R_{nr}(K,\mathbb{Z}/2\mathbb{Z})=\emptyset\quad\text{and}\quad R'_{nr}(K,\mathbb{Z}/2\mathbb{Z})=\{1\}.\]

(3) Suppose $h_K^+$ is even and $K$ is of class type. Since $R'_{nr}(K,\mathbb{Z}/2\mathbb{Z})\subset Cl_{K,2}$ and $|Cl_{K,2}|=2^{r_2(K)}$, it suffices to prove that these groups have the same order. Proposition \ref{2rang}(2) and Galois theory give $n=2^{r_2^+(K)}-1$ unramified quadratic extensions, denoted by $K_i/K$ for $1\leq i\leq n$. For distinct $i,j$, the third quadratic subextension of the unramified biquadratic compositum $K_iK_j$ is some $K_m$ with $m\notin\{i,j\}$. By Corollary \ref{cor4},
\[
cl_K(\mathcal O_{K_i})cl_K(\mathcal O_{K_j})=cl_K(\mathcal O_{K_m}).
\]
By Corollary \ref{cor2}, the right-hand side has order $2$, so the two classes on the left are distinct. Thus $R_{nr}(K,\mathbb{Z}/2\mathbb{Z})$ has cardinality $2^{r_2^+(K)}-1$, and $R'_{nr}(K,\mathbb{Z}/2\mathbb{Z})$ has order $2^{r_2^+(K)}$. The inclusion in $Cl_{K,2}$ gives $r_2^+(K)\leq r_2(K)$. Since $Cl_K$ is a quotient of $Cl_K^+$, Proposition \ref{2rang}(5) gives $r_2(K)\leq r_2^+(K)$. Hence $r_2(K)=r_2^+(K)$ and the two assertions follow.

(4) Suppose $K$ is of unit type. By the definition of $m$, there are $2^m-1$ unramified quadratic extensions $M_i=K(\sqrt{\varepsilon_i})$ with trivial Steinitz class. Thus the number with nontrivial Steinitz class is
\[
a=(2^{r_2^+(K)}-1)-(2^m-1)=2^{r_2^+(K)}-2^m.
\]
If $L=K(\sqrt{\alpha})$ is one of these, exactly $2^m$ quadratic extensions have its Steinitz class. They are $L$ and the third quadratic subextensions $K(\sqrt{\alpha\varepsilon_i})$ of the biquadratic composita $LM_i$. Hence there are $a/2^m=2^{r_2^+(K)-m}-1$ nontrivial Steinitz classes. Adding the trivial class completes the proof.
\end{proof}

\noindent\textbf{Remarks.}

\textbf{1.} If $K$ is of class type and $h_K^+$ is even, the proof gives $r_2(K)=r_2^+(K)$. Thus every quadratic extension of $K$ unramified at the finite places is also unramified at the infinite places.

\textbf{2.} If $r_2(K)\ne r_2^+(K)$, then $K$ is of unit type, by the first remark and the fact that odd $h_K^+$ gives $r_2(K)=r_2^+(K)=0$. $\blacksquare$

\medskip
\noindent\textbf{Explicit examples for Theorem \ref{t1} in the case of quadratic extensions $K/\mathbb{Q}$.}

We use the notation of Proposition \ref{K/kquadratique} and Corollary \ref{corK/kquadratique}.

For an abelian extension $L/\mathbb{Q}$, its narrow genus field $L^+_{gen}$ is the largest abelian extension of $\mathbb Q$ containing $L$ and unramified over $L$ at all finite places. Recall also the following class field theory result. If $L'/L$ is a quadratic extension of number fields and $L''/L'$ is an unramified abelian extension of degree coprime to $h_L^+$, then $L''/L$ is Galois.

The next proposition identifies the field containing all the unramified quadratic extensions of $K$.

\begin{Prop}\label{genre}
Let $K/\mathbb Q$ be quadratic and let $M$ be the compositum of all unramified quadratic extensions of $K$. Then $M=K^+_{gen}$.
\end{Prop}

\begin{proof}
For the description of $K^+_{gen}$, see \cite[(1.4), p. 4]{Ichida}. In particular, $[K^+_{gen}:K]=2^{u-1}$ by \cite[(1.5), p. 5]{Ichida}. Let $L/K$ be an unramified quadratic extension. Since $[K:\mathbb Q]=2$ and $\gcd(h_{\mathbb Q}^+,[L:K])=1$, the preceding class field theory result shows that $L/\mathbb Q$ is Galois. Its degree is $4$, so it is abelian. Hence $M\subseteq K^+_{gen}$. Since $[M:K]=2^{r_2^+(K)}$ and $r_2^+(K)=u-1$ by Corollary \ref{corK/kquadratique}(1), equality follows.
\end{proof}

Recall the explicit description of $K^+_{gen}$ for $K=\mathbb Q(\sqrt d)$ with $d$ a squarefree integer \cite{Ichida}. Let $p_1,\ldots,p_u$ be the primes dividing $\delta_{K/\mathbb Q}$. If $d\equiv1\bmod4$, set $p_i^*=\left(\frac{-1}{p_i}\right)p_i=(-1)^{(p_i-1)/2}p_i$. If $d\equiv2$ or $3\bmod4$, take $p_1=2$. Set $p_1^*=-4$ when $d\equiv3\bmod4$, $p_1^*=8$ when $d\equiv2\bmod8$, and $p_1^*=-8$ when $d\equiv-2\bmod8$. For $2\leq i\leq u$, set $p_i^*=(-1)^{(p_i-1)/2}p_i$. Then
\[
K^+_{gen}=K(\sqrt{p_i^*}\mid 2\leq i\leq u)
=\mathbb Q(\sqrt d,\sqrt{p_i^*}\mid 2\leq i\leq u).
\]
Write $\mathfrak p_i$ for the prime ideal of $\mathcal O_K$ above $p_i$.

\begin{Cor}\label{corgenre}
With the preceding notation, $R'_{nr}(K,\mathbb{Z}/2\mathbb{Z})$ is the subgroup of $Am_0(K/\mathbb Q)$ generated by $[\mathfrak p_i]_K$ for $2\leq i\leq u$. Moreover, $R'_{nr}(K,\mathbb{Z}/2\mathbb{Z})=Am_0(K/\mathbb Q)$ if $d\equiv1$ or $2\bmod4$.
\end{Cor}

\begin{proof}
The first assertion follows from Proposition \ref{genre} and the description of $K^+_{gen}$. If $d\equiv1$ or $2\bmod4$, then $\prod_{i=1}^u\mathfrak p_i=\sqrt d\,\mathcal O_K$. Thus $[\mathfrak p_1]_K\in R'_{nr}(K,\mathbb{Z}/2\mathbb{Z})$, which gives the second assertion.
\end{proof}

\textbf{1. Odd $h_K^+$.}

In this case, $K$ is of class type, $R_{nr}(K,\mathbb{Z}/2\mathbb{Z})=\emptyset$, and $R'_{nr}(K,\mathbb{Z}/2\mathbb{Z})=\{1\}$.

\begin{Prop}
The narrow class number $h_K^+$ is odd if and only if $K$ is one of the following fields: $\mathbb Q(\sqrt{-1})$, $\mathbb Q(\sqrt{\pm2})$, $\mathbb Q(\sqrt p)$ for a prime $p\equiv1\bmod4$, or $\mathbb Q(\sqrt{-p})$ for a prime $p\equiv3\bmod4$.
\end{Prop}

\begin{proof}
The condition that $h_K^+$ be odd is equivalent to $r_2^+(K)=0$, or equivalently $u=1$ by Corollary \ref{corK/kquadratique}(1). The listed fields are precisely the quadratic fields with only one ramified finite prime, as follows from the quadratic discriminant formula.
\end{proof}

\textbf{2. Even $h_K^+$.}

\textbf{2.1. Fields of class type.}

We have $R'_{nr}(K,\mathbb{Z}/2\mathbb{Z})=Cl_{K,2}$. Moreover, Propositions \ref{imaginaire}(ii) and \ref{reel}(ii), together with Corollary \ref{descclasseK2}, give $R'_{nr}(K,\mathbb{Z}/2\mathbb{Z})=Am_0(K/\mathbb Q)$.

\textbf{2.2. Fields of unit type.}

Recall that $|R_{nr}(K,\mathbb{Z}/2\mathbb{Z})|=2^{r_2^+(K)-m}$ with $m\geq1$.

\textbf{The case $r_2^+(K)=r_2(K)$.}
Here $K=\mathbb Q(\sqrt d)$ is either imaginary with $d\ne-1$ and $d\equiv3\bmod4$, or real with $N_{K/\mathbb Q}(\epsilon)=1$ and $-1$ a norm from $K$. Since $m\geq1$, the subgroup $R_{nr}(K,\mathbb{Z}/2\mathbb{Z})$ is proper in $Cl_{K,2}$. In the imaginary case, it has order $2^{r_2(K)-1}$.

In the real case, suppose $d$ is odd. Then $d\equiv1\bmod4$ by the discussion of Hasse's theorem following Corollary \ref{corK/kquadratique}. Since $2$ is unramified in $K/\mathbb Q$, it ramifies in $K(\sqrt{-1})/K$. Hence $m=1$ and $|R_{nr}(K,\mathbb{Z}/2\mathbb{Z})|=2^{r_2(K)-1}$. If $d\equiv1\bmod8$, then by \cite[Theorem 1.1(i)]{LA} the total positivity of $\epsilon$ implies $\epsilon\equiv1\bmod4\mathcal O_K$. Thus $K(\sqrt{\epsilon})/K$ is unramified and $K(\sqrt{-\epsilon})/K$ is ramified.

\textbf{The case $r_2^+(K)=r_2(K)+1$.}
Here $K$ is real, $N_{K/\mathbb Q}(\epsilon)=1$, and $-1$ is not a norm from $K$ (Corollary \ref{r^+_2(K)= r_2(K)+1}). In this case $Cl_{K,2}=Am_0(K/\mathbb Q)$ has order $2^{u-2}$ and $|R_{nr}(K,\mathbb{Z}/2\mathbb{Z})|=2^{r_2(K)+1-m}$. The next proposition gives infinitely many examples with $R_{nr}(K,\mathbb{Z}/2\mathbb{Z})=Cl_{K,2}$, equivalently $m=1$.

\begin{Prop}\label{exempletheoreme1.1unite}
Let $K=\mathbb Q(\sqrt d)$ be real quadratic, with $d$ a positive squarefree integer. Suppose that all prime divisors of $d$ are congruent to $3$ modulo $4$ and that their number is positive and even. Then $K$ is of unit type and $R_{nr}(K,\mathbb{Z}/2\mathbb{Z})=Cl_{K,2}=Am_0(K/\mathbb Q)$.
\end{Prop}

\begin{proof}
Some prime $p\equiv3\bmod4$ divides the discriminant. Hence $N_{K/\mathbb Q}(\epsilon)=1$ (see \cite[(1.7), p. 179]{FT}). By Proposition \ref{reel}(i), $K$ is of unit type. Since $\left(\frac{-1}{p}\right)=-1$, the element $-1$ is not a norm from $K$, so $r_2^+(K)=r_2(K)+1$. The conditions on the prime divisors imply $d\equiv1\bmod4$, so $\delta_{K/\mathbb Q}=d\mathbb Z$. Thus $2$ is unramified in $K/\mathbb Q$ and $K(\sqrt{-1})/K$ is ramified. Remark 2 following Theorem \ref{t1} then gives $m=1$. The set of realizable classes therefore has order $2^{r_2(K)}$ and equals $Cl_{K,2}$.
\end{proof}

\section{Proof of Theorem 1.2}

\begin{proof}[\textbf{Proof of Theorem \ref{t2}}]
Let $\Gamma$ be a finite group. If $\Gamma$ has odd order or has a noncyclic Sylow $2$-subgroup, Corollary \ref{cor3}(i) gives $R'_{nr}(K,\Gamma)=\{1\}$, which proves the subgroup assertion in this case.

Henceforth assume that $\Gamma$ has even order and that its Sylow $2$-subgroups are cyclic. To keep the argument as general as possible, we will invoke normality, and hence uniqueness, of the Sylow $2$-subgroup only where needed.

Corollary \ref{cor3}(ii) gives $R_{nr}(K,\Gamma)\subseteq R_{nr}(K,\mathbb{Z}/2\mathbb{Z})$ by taking the unique quadratic subextension of each unramified $\Gamma$-extension. Adjoining the trivial class gives the inclusion for the sets $R'_{nr}$. Since $R'_{nr}(K,\Gamma)$ is finite and nonempty, it suffices to show that it is closed under multiplication in $Cl_K$.

Let $c_1,c_2\in R'_{nr}(K,\Gamma)$. If either class is trivial, or if $c_1=c_2$, then $c_1c_2\in R'_{nr}(K,\Gamma)$. We may therefore assume that $c_1$ and $c_2$ are distinct and nontrivial.

Let $N_1/K$ and $N_2/K$ be unramified $\Gamma$-extensions with $cl_K(\mathcal{O}_{N_1})=c_1$ and $cl_K(\mathcal{O}_{N_2})=c_2$. For $i=1,2$, Proposition \ref{prop1}(ii) gives a unique quadratic subextension $L_i/K$ of $N_i/K$. Corollary \ref{cor3}(ii) gives $c_i=cl_K(\mathcal{O}_{L_i})$, so
\[
c_1c_2=cl_K(\mathcal{O}_{L_1})cl_K(\mathcal{O}_{L_2}).
\]

Consider the following diagram.
\begin{center}

  \begin{tikzpicture}[node distance = 1cm, auto]
  \node (Q) {$K$};
  \node (E) [above of=Q] {$N_{1}\cap N_{2}$};
  \node (E1) [above of=E, left of=E] {$N_{1}$};
  \node (F) [above of=E, right of=E] {$N_{2}$};
  \node (K) [above of=E, node distance =2cm] {$N_{1}N_{2}$};
   \draw (Q) --(E);
  \draw (E) -- (F);
  \draw (E) --  (E1);
  \draw (F) --  (K);
  \draw (E1) -- (K);
  \end{tikzpicture}

  \end{center}
Put $G=Gal(N_1N_2/K)$, $G'=Gal(N_1N_2/N_1\cap N_2)$, $G_i=Gal(N_i/K)\simeq\Gamma$, $H_i=Gal(N_i/N_1\cap N_2)$, $H'_1=Gal(N_1N_2/N_2)$, and $H'_2=Gal(N_1N_2/N_1)$. We distinguish two cases according to the parity of $[N_1\cap N_2:K]$.

\textbf{Case 1.} The degree $[N_1\cap N_2:K]$ is even.

By Galois theory, $Gal(N_1\cap N_2/K)\simeq G_1/H_1$. Proposition \ref{prp3}(1.2) gives $syl_2(G_1/H_1)=syl_2(G_1)H_1/H_1$, which is cyclic because $syl_2(G_1)$ is cyclic. By Proposition \ref{prop1}(ii), $N_1\cap N_2$ has a unique quadratic subfield over $K$. Thus $L_1=L_2$ and $c_1c_2=1\in R'_{nr}(K,\Gamma)$.

\noindent\textbf{Remark.} This case requires only that the Sylow $2$-subgroups of $\Gamma$ be cyclic. $\blacksquare$

\textbf{Case 2.} The degree $[N_1\cap N_2:K]$ is odd.

Let $L_3/K$ be the third quadratic subextension of the biquadratic extension $L_1L_2/K$. By Corollary \ref{cor4}, $c_1c_2=cl_K(\mathcal{O}_{L_3})$. It suffices to embed $L_3/K$ in an unramified $\Gamma$-extension of $K$.

Write $|\Gamma|=2^nk$ with $n\geq1$ and $k$ odd. If $M/K$ is a Galois subextension of $N_1N_2/K$ containing $N_1\cap N_2$ and $H=Gal(M/N_1\cap N_2)$, then a Sylow $2$-subgroup of $H$ is also a Sylow $2$-subgroup of $Gal(M/K)$. Indeed, the index of $H$ is $[N_1\cap N_2:K]$, which is odd.

Applying this observation to $M=N_1N_2$ and $M=N_i$ gives $syl_2(G)=syl_2(G')$ and $syl_2(G_i)=syl_2(H_i)$ for compatible choices of Sylow subgroups. Since $N_1$ and $N_2$ are linearly disjoint over $N_1\cap N_2$, we have $G'\simeq H_1\times H_2$. Proposition \ref{prp3}(2) now gives
\[
syl_2(G)\simeq syl_2(H_1)\times syl_2(H_2)
\simeq syl_2(\Gamma)\times syl_2(\Gamma).
\]
The restriction isomorphisms $H'_i\simeq H_i$ also give $syl_2(H'_i)\simeq syl_2(H_i)\simeq syl_2(\Gamma)$. Both observations hold for arbitrary Sylow $2$-subgroups of $\Gamma$.

Cyclicity now gives $syl_2(G)\simeq(\mathbb{Z}/2^n\mathbb{Z})^2$, so $r_2(syl_2(G))=2$. Choose generators $\sigma$ of $syl_2(H'_2)$ and $\tau$ of $syl_2(H'_1)$ in the two factors of $syl_2(G)$. Then $syl_2(G)=\langle\sigma,\tau\rangle$.

For the rest of the proof assume that $syl_2(\Gamma)$ is normal and hence unique. Let $I$ be the unique normal subgroup of $\Gamma$ of order $k$ (Proposition \ref{Scott}). Since $syl_2(\Gamma)$ and $I$ are normal and have coprime orders, $\Gamma$ is their internal direct product by Proposition \ref{prp3}(4). The same follows from the remark after Proposition \ref{Scott}. Thus $\Gamma\simeq syl_2(\Gamma)\times I$.

By Galois theory, $G$ embeds in
\[
\Gamma\times\Gamma\simeq syl_2(\Gamma)\times syl_2(\Gamma)\times I\times I.
\]
The two factors $syl_2(\Gamma)\times syl_2(\Gamma)$ and $I\times I$ have coprime orders. By Proposition \ref{prp3}(3), $G$ is a direct product of a subgroup $R$ of the first factor and a subgroup $S$ of the second. Since $syl_2(G)\simeq syl_2(\Gamma)\times syl_2(\Gamma)$, we obtain $G\simeq syl_2(G)\times S$. As in the proof of Proposition \ref{2rang}(3),
\[
G/\langle G^2\rangle\simeq
\bigl(syl_2(G)/\langle syl_2(G)^2\rangle\bigr)
\times\bigl(S/\langle S^2\rangle\bigr).
\]
Every element of the odd-order group $S$ is a square. Thus \[r_2(G)=r_2(syl_2(G))=2.\] By Proposition \ref{2rang}(2), the subgroups of index $2$ are
\[
\Gamma_1=\langle\sigma,\tau^2\rangle S,\quad
\Gamma_2=\langle\tau,\sigma^2\rangle S,\quad
\Gamma_3=\langle\sigma\tau,\tau^2\rangle S,
\]
where $S$ is identified with a subgroup of $G$. Since $\sigma$ fixes $N_1$ and $\tau$ fixes $N_2$, the subgroups $\Gamma_i$ correspond to the fields $L_i$ for $i=1,2,3$.

Since $\gcd(|syl_2(G)|,|S|)=1$, Proposition \ref{prp3}(3) expresses $H'_1$ as a direct product of a subgroup of $syl_2(G)$ and a subgroup $S_1$ of $S$. Both factors are normal in their respective ambient groups because $H'_1$ is normal in $G$. As $syl_2(H'_1)=\langle\tau\rangle$, we have $H'_1=\langle\tau\rangle S_1$.

Consider $\Gamma_4=\langle\sigma\tau\rangle S_1$. We use Proposition \ref{prp3} throughout the rest of the argument. The two factors of $\Gamma_4$ are normal in $syl_2(G)$ and $S$, respectively, so $\Gamma_4$ is normal in $G$. Let $N_3$ be the fixed field of $\Gamma_4$ in $N_1N_2$. Then
\[
Gal(N_3/K)\simeq G/\Gamma_4
\simeq \bigl(\langle\sigma,\tau\rangle/\langle\sigma\tau\rangle\bigr)\times(S/S_1)
\simeq \mathbb{Z}/2^n\mathbb{Z}\times S/S_1.
\]
Similarly,
\[
\Gamma\simeq Gal(N_2/K)\simeq G/H'_1
\simeq \bigl(\langle\sigma,\tau\rangle/\langle\tau\rangle\bigr)\times(S/S_1)
\simeq \mathbb{Z}/2^n\mathbb{Z}\times S/S_1.
\]
Hence $Gal(N_3/K)\simeq\Gamma$. Since $\Gamma_4\subseteq\Gamma_3$, we have $L_3\subseteq N_3$. Moreover, $N_3/K$ is unramified because it is a subextension of the compositum of two unramified extensions. This completes the proof.
\end{proof}

\noindent\textbf{Remark.} The same argument writes $H'_2=\langle\sigma\rangle S_2$ with $S_2$ normal in $S$. Taking $\Gamma_5=\langle\sigma\tau\rangle S_2$ and letting $N_4$ be its fixed field in $N_1N_2$, we likewise obtain an unramified $\Gamma$-extension $N_4/K$ containing $L_3/K$. $\blacksquare$

The following proposition illustrates Theorem \ref{t2} for $\Gamma\simeq\mathbb{Z}/4\mathbb{Z}$.

Let $H^+_{K,2}$ be the narrow Hilbert $2$-class field of $K$, that is, the maximal abelian $2$-extension of $K$ unramified at all finite places. Recall that $Gal(H^+_{K,2}/K)\simeq syl_2(Cl_K^+)$. Every unramified $\mathbb{Z}/4\mathbb{Z}$-extension $L/K$ satisfies $L\subseteq H^+_{K,2}$. Write $H_0=Gal(H^+_{K,2}/K)$.

\begin{Prop}\label{exempletheor1.2}
Let $K$ be a number field of class type whose narrow $2$-class group has type $(2,4)$, so \[syl_2(Cl_K^+)\simeq\mathbb{Z}/2\mathbb{Z}\times\mathbb{Z}/4\mathbb{Z}.\] Then $\left|R'_{nr}(K,\mathbb{Z}/4\mathbb{Z})\right|=2$.
\end{Prop}

\begin{proof}
Choose elements $\tau$ and $\sigma$ of $H_0$ of orders $2$ and $4$, respectively, with
\[
H_0=\langle\tau,\sigma\rangle=\{\tau^i\sigma^j\mid 0\leq i\leq1,\ 0\leq j\leq3\}.
\]
The three subgroups of order $4$ and index $2$ are $H_1=\langle\sigma\rangle$, $H_2=\langle\tau\sigma\rangle$, and $H_3=\langle\tau,\sigma^2\rangle$. The three subgroups of order $2$ and index $4$ are $H_4=\langle\tau\rangle$, $H_5=\langle\tau\sigma^2\rangle$, and $H_6=\langle\sigma^2\rangle$.

Set $H_7=\{1\}$, $L_0=K$, and $L_7=H^+_{K,2}$. The Galois correspondence associates to each $H_i$, for $0\leq i\leq7$, a unique unramified abelian extension $L_i/K$ with $Gal(L_i/K)\simeq H_0/H_i$.

The only unramified cyclic extensions of degree $4$ over $K$ are $L_4$ and $L_5$. By Corollary \ref{cor3}(ii),
\[
cl_K(\mathcal{O}_{L_4})=cl_K(\mathcal{O}_{L_5})=cl_K(\mathcal{O}_{L_3}).
\]
Thus $R'_{nr}(K,\mathbb{Z}/4\mathbb{Z})=\{1,cl_K(\mathcal{O}_{L_3})\}$. Since $K$ is of class type, $cl_K(\mathcal{O}_{L_3})\ne1$, and the set has order $2$.
\end{proof}

\noindent\textbf{Remark.} Under the hypotheses of the preceding proposition,
\[R'_{nr}(K,\mathbb{Z}/2\mathbb{Z})=\{1,cl_K(\mathcal{O}_{L_i})\mid 1\leq i\leq3\}=Cl_{K,2}.\qquad\blacksquare\]
\medskip

\textbf{Example.}
Let $K=\mathbb{Q}(\sqrt{p_0p_1p_2})$ be a real quadratic field, where the primes $p_i$, for $0\leq i\leq2$, are congruent to $1$ modulo $4$ and satisfy
\[
\left(\frac{p_1}{p_0}\right)=\left(\frac{p_2}{p_0}\right)=1,\qquad
\left(\frac{p_1}{p_2}\right)=-1,
\]
and either $\left(\frac{p_0}{p_1p_2}\right)_4=-1$ or $\left(\frac{p_1p_2}{p_0}\right)_4=-1$, where $\left(\frac{a}{b}\right)_4$ denotes the quartic residue symbol. By \cite[Proposition B1, p. 326]{Ka76} (see also \cite[Proposition B, p. 327]{Ka76}), $syl_2(Cl_K^+)$ has type $(2,4)$. Proposition \ref{reel} shows that $K$ is of class type if and only if its fundamental unit $\epsilon$ has norm $-1$.

Take $p_0=5$, $p_1=13$, and $p_2=101$, so $K=\mathbb{Q}(\sqrt{6565})$. PARI/GP \cite{Pari} gives $syl_2(Cl_K^+)$ of type $(2,4)$ and $\epsilon=(81+\sqrt{6565})/2$, of norm $-1$. This field therefore satisfies the hypotheses of the preceding proposition. $\blacksquare$

We first proved Theorem \ref{t2} for abelian groups and then for nilpotent groups. In both cases normality of the Sylow $2$-subgroup was the property that led to the more general statement.

We then asked whether our argument in Case 2, where $[N_1\cap N_2:K]$ is odd, also works when the Sylow $2$-subgroup is nonnormal. Case 1 already holds without assuming normality. Consider the dihedral group $\Gamma=D_n$ of order $2n$, where $n\geq3$ is odd. If $[N_1\cap N_2:K]$ is odd, then $[N_1\cap N_2:K]=1$. Thus $N_1$ and $N_2$ are linearly disjoint over $K$, and $Gal(N_1N_2/K)\simeq D_n\times D_n$.

The group $D_n$ has presentation $D_n=\langle\sigma,\tau\mid \sigma^n=\tau^2=1,\ \tau\sigma\tau=\sigma^{-1}\rangle$. We have $r_2(D_n\times D_n)=2$, and the three subgroups of index $2$ are $D_n\times\langle\sigma\rangle$, $\langle\sigma\rangle\times D_n$, and
\[
(\langle\sigma\rangle\times\langle\sigma\rangle)
\cup(\langle\sigma\rangle\tau\times\tau\langle\sigma\rangle).
\]
Using SageMath \cite{Sage}, we checked that for $n\leq100$ the only normal subgroups $J$ of $D_n\times D_n$ for which $(D_n\times D_n)/J\simeq D_n$ are $D_n\times\{1\}$ and $\{1\}\times D_n$. Under the Galois correspondence these give $N_1/K$ and $N_2/K$. Thus, for these values of $n$, $L_3/K$ cannot be embedded in any subextension of $N_1N_2/K$ with Galois group isomorphic to $D_n$. An embedding into an unramified $D_n$-extension outside $N_1N_2$ could provide a way to generalize Theorem \ref{t2}.

\clearpage
\selectlanguage{french}
\setcounter{section}{0}
\setcounter{equation}{0}
\setcounter{Df}{0}
\setcounter{Rem}{0}
\setcounter{Conj}{0}
\pdfbookmark[0]{Version originale en français}{original-fr}
\begin{center}
{\LARGE Classes de Steinitz d'extensions galoisiennes non ramifiées\par}
\bigskip
{\large Ayoub Al Uartassi, Bouchaïb Sodaïgui et Mohammed Taous\par}
\medskip
\emph{Version originale en français}\par
\emph{Acta Arithmetica} \textbf{219} (2025), no~1, 53--80.
\end{center}

\centerline { Mathematics Subject Classification 2020: 11R33}

\centerline{\bf {Abstract}}

Let $K$ be a number field and $Cl_{K}$ its class group. Let $\Gamma$ be a finite group. We denote by $ R_{nr}(K,\Gamma)$  the subset of $Cl_{K}$  formed by the  classes which are realizable as Steinitz classes of   Galois extensions  of $K$,  unramified at the finite places of $K$ and having a Galois group isomorphic to $\Gamma$. If $\Gamma$ is abelian and the  narrow class number of $K$ is prime to the order of $\Gamma$, then  $R_{nr}(K,\Gamma)=\emptyset$. In the present article, for any $\Gamma$ we consider the set
 $R'_{nr}(K,\Gamma):=\lbrace 1 \rbrace \cup R_{nr}(K,\Gamma)$.
  We prove that  $R'_{nr}(K, \mathbb{Z}/2\mathbb{Z})$  is a subgroup of  $Cl_{K,2}:=\lbrace c\in Cl_{K}\vert c^{2}=1\rbrace$ ; furthermore, it is equal to $Cl_{K,2}$ under a certain assumption on $K$. Using this result   we show that $R'_{nr}(K,\Gamma)$ is a subgroup of $R'_{nr}(K, \mathbb{Z}/2\mathbb{Z})$ if $\Gamma$ either  has odd order, or  has a noncyclic  $2$-Sylow subgroup (for instance $\Gamma$ a nonabelian $2$-group,  $\Gamma=S_n$ or $A_n$, with $n\geq 4$), or has a normal cyclic  $2$-Sylow subgroup
  (for instance  $\Gamma$ nilpotent having even order).

\textbf{Keywords} :   Rings of integers; Steinitz classes;  Unramified extensions; Embedding problem;  Realizable Galois module classes;  Class group; Class number;   $2$-rank of a group.

\section{Introduction}

 \ \ Dans tout cet article, si $K$ est un corps de nombres,  $\mathcal{O}_{K}$ désigne son anneau d'entiers, $Cl_{K}$ son groupe de classes, $h_K$ son nombre de classes,  $Cl^{+}_{K}$ son groupe de classes au sens restreint, $h^{+}_K$ son nombre de classes au sens restreint et si $I$ est un idéal fractionnaire de $K$ on note  $[I]_K$ sa classe dans $Cl_{K}$.  Si une  extension finie $L/K$ est  non ramifiée aux places finies, on convient de  dire simplement qu'elle est non ramifiée et quand il s'agira de la non ramification aux places infinies on le précisera.

Rappelons la définition de la classe de Steinitz. Soit $L/K$ une extension finie de corps de nombres de degré $n$. L'anneau $\mathcal{O}_{L}$ étant un $\mathcal{O}_{K}$-module sans torsion de rang $n$,  il existe un idéal fractionnaire $I$ de $\mathcal{O}_{K}$ tel que $ \mathcal{O}_{L}\simeq \mathcal{O}_{K}^{n-1}\oplus I$ en tant que $ \mathcal{O}_{K}$-module (pour plus de détails voir \cite [Théorème 1.2.19, p. 9]{fr-cohan}) ; $[I]_K$ est appelée la classe de Steinitz de l'anneau $\mathcal{O}_{L}$  ou de l'extension $L/K$ et on la note  $cl_{K}(\mathcal{O}_{L})$. La structure de $\mathcal{O}_{K}$-module de $\mathcal{O}_{L}$ est complètement  déterminée par $cl_{K}(\mathcal{O}_{L})$ ;  en particulier, $\mathcal{O}_{L}$ est un $\mathcal{O}_{K}$-module libre si et seulement si $cl_{K}(\mathcal{O}_{L})=1$.

   Soient $\Gamma$ un groupe fini et $K$ un corps de nombres. On appelle $\Gamma$-extension de $K$  une extension galoisienne de $K$ dont le  groupe de Galois est  isomorphe à $\Gamma$.
    On note $R_{m}(K,\Gamma)$ (resp. $R_{nr}(K,\Gamma)$)   l'ensemble des
  $ c\in Cl_{K}$ telles qu'il existe  une $\Gamma$-extension $L/K$ modérée (resp. non ramifiée) satisfaisant  $c= cl_{K}(\mathcal{O}_{L})$.

   Il est conjecturé que $R_{m}(K,\Gamma)$ est  un sous-groupe de $Cl_{K}$ (voir \cite[Conjecture 3, \S1, p. 6]{fr-BGS}). Cette conjecture est vraie lorsque $\Gamma$ est abélien d'après \cite{fr-Mc} comme on va le voir ci-dessous. Mais elle est  toujours  ouverte dans le cas non  abélien ;  on peut consulter par exemple  les articles \cite{fr-BrS,fr-BGS,fr-BS,fr-SOCA,fr-bouchaib} qui démontrent qu'elle est vraie pour certains groupes non abéliens.

Lorsque $L/K$ est une $\Gamma$-extension  moderée, on peut associer à $\mathcal{O}_{L}$ une classe notée $[\mathcal{O}_{L}]$ dans le groupe de classes $Cl(\mathcal{O}_{K}[\Gamma])$ des $\mathcal{O}_{K}[\Gamma]$-modules localement libres. On désigne par $R_{m}(\mathcal{O}_{K}[\Gamma])$ (resp. $R_{nr}(\mathcal{O}_{K}[\Gamma])$) l'ensemble des classes $c$ de $Cl(\mathcal{O}_{K}[\Gamma])$ telle qu'il existe une $\Gamma$-extension $L/K$ modérément ramifiée (resp. non ramifiée),  avec $c=[\mathcal{O}_{L}]$.

Dans le cas où $\Gamma$  est abélien, dans \cite{fr-Mc} McCulloh a étudié  $\mathfrak{R}_{m}(\mathcal{O}_{K}[\Gamma])$ l'ensemble des classes de $Cl(\mathcal{O}_{K}[\Gamma])$  réalisables par les classes des clôtures intégrales de $\mathcal{O}_K$ dans  certaines algèbres galoisiennes modérément ramifiées sur $K$ (pour plus de précisions : voir \cite[ (1.1), p. 263 ; p. 271 après (2.2) ; le premier paragraphe de la section  \S3, p. 276.]{fr-Mc}). Il a étudié aussi le sous-ensemble $\mathfrak{R}_{nr}(\mathcal{O}_{K}[\Gamma])$ de $\mathfrak{R}_{m}(\mathcal{O}_{K}[\Gamma])$, où $nr$ désigne non ramifié aux places finies.  En utilisant une correspondance de Stickelberger, il a montré que ces deux ensembles sont des sous-groupes de $Cl(\mathcal{O}_{K}[\Gamma])$ (Voir \cite[ (6.20) Corollary, p. 291]{fr-Mc}). De plus, d'après \cite[ (6.17) Theorem, p. 291]{fr-Mc}),
$\mathfrak{R}_{m}(\mathcal{O}_{K}[\Gamma])=R_{m}(\mathcal{O}_{K}[\Gamma])$, par suite ce dernier
   est un sous-groupe de $Cl(\mathcal{O}_{K}[\Gamma])$, et donc, en oubliant l'action de $\Gamma$,   $R_{m}(K,\Gamma)$ est un sous-groupe de  $Cl_{K}$.

Signalons que dans \cite{fr-agboola} les auteurs ont raffiné (voir \cite[Théorème D]{fr-agboola}) et   considérablement  généralisé (voir \cite[Théorème E]{fr-agboola})  le résultat de McCulloh disant que $\mathfrak{R}_{m}(\mathcal{O}_{K}[\Gamma])$ est un sous-groupe de $Cl(\mathcal{O}_{K}[\Gamma])$  pour une large famille de  groupes pouvant être non abéliens sous certaines conditions sur $K$.  Cependant l'étude de $\mathfrak{R}_{nr}(\mathcal{O}_{K}[\Gamma])$ n'a pas été considérée.

   Dans \cite{fr-Mc}, rien n'est dit sur la comparaison de $\mathfrak{R}_{nr}(\mathcal{O}_{K}[\Gamma])$ et $R_{nr}(\mathcal{O}_{K}[\Gamma])$. C'était le point de départ du présent article. S'il y avait égalité entre ces deux ensembles, on en déduirait (en oubliant l'action de $\Gamma$) que $R_{nr}(K,\Gamma)$ serait un sous-groupe de $Cl(K)$ lorsque $\Gamma$ est abélien.

   Rappelons que  le corps de classes de Hilbert de $K$, noté $H_K$, est l'extension abélienne  non ramifiée aux places finies et infinies maximale de $K$ et
 $Gal(H_K/K)\simeq Cl_{K}$.
Rappelons aussi que  le corps de classes de Hilbert de $K$  au sens restreint, noté $H^+_K$, est l'extension abélienne  non ramifiée (aux places finies) maximale de $K$ et
 $Gal(H^+_K/K)\simeq Cl^{+}_{K}$.
   Nous avons tout de suite constaté le fait suivant : si $h^{+}_K$   est impair et $\Gamma$ est abélien d'ordre pair, alors $R_{nr}(\mathcal{O}_{K}[\Gamma])$ et $R_{nr}(K,\Gamma)$  sont égaux à l'ensemble vide.
   Clairement, dans ce cas on a $\mathfrak{R}_{nr}(\mathcal{O}_{K}[\Gamma])\not =R_{nr}(\mathcal{O}_{K}[\Gamma])$.

   Nous avons alors défini l'ensemble suivant pour un groupe fini $\Gamma$ quelconque : $$R'_{nr}(K,\Gamma):=\lbrace 1\rbrace \cup R_{nr}(K,\Gamma).$$ L'objectif principal du présent article est l'étude de cet ensemble.

   On dit que $K$ est de type classe si toute extension quadratique de $K$  de la forme $K(\sqrt{\varepsilon})$, où $\varepsilon$ est une unité de $K$, est ramifiée en au moins une place finie ; sinon on dit qu'il est de type unité (pour plus de précisions voir Définition \ref{fr-Dftypeclasse} et le paragraphe qui la suit). Noter que si $h^+_K$ est impair, alors $K$ est de type classe, et  si $K$ est de type unité, alors  $h^{+}_K$ est pair.

 On rappelle le lien bien connu entre $h^+_K$ et $h_K$ : $h^+_K=h_K \frac{2^r}{[E_K:E^+_K]}$, où $r$ est le nombre de places réelles de $K$, $E_K$ le groupe des unités de $K$ et $E^+_K$ son sous-groupe des unités totalement positives. En particulier, $h^+_K=h_K$ si $h^+_K$ est impair ou $K$ est totalement imaginaire.

On note  $r^+_2(K)$ le $2$-rang de $K$ au sens restreint ; rappelons que  $r^+_2(K)$ est le $2$-rang du groupe abélien  $Cl^+_{K}$ et $r^+_2(K)=\dim_{\mathbb{F}_2}(Cl^+_{K}/Cl^+_{K}{^2})$, où $\mathbb{F}_2$ est le corps à deux éléments. (Voir \S2 pour plus de détails sur le $2$-rang.)

Les principaux résultats  sont les deux théorèmes suivants.

\begin{Th}\label{fr-t1} Soient $K$ un corps de nombres et $Cl_{K,2}:=\lbrace c\in Cl_{K}\vert c^{2}=1\rbrace$ (le $2$-sous-groupe élémentaire de $Cl_{K}$ si $h_K$ est pair).  Alors :

$(1)$   $ R'_{nr}(K,\mathbb{Z}/2\mathbb{Z})$ est un sous-groupe de $Cl_{K, 2}$.

$(2)$ Si $h^+_K$ est impair, alors $R_{nr}(K,\mathbb{Z}/2\mathbb{Z})=\emptyset$ et  $R'_{nr}(K,\mathbb{Z}/2\mathbb{Z})=\lbrace 1\rbrace$.

 $(3)$ Si $K$ est de type classe et $h^+_K$ est pair, alors
$$R_{nr}(K,\mathbb{Z}/2\mathbb{Z})= Cl_{K,2}\setminus \{1\} \ \ \text{et } \ \ R'_{nr}(K,\mathbb{Z}/2\mathbb{Z})=Cl_{K,2}.$$

$(4)$ Si $K$ est de type unité, alors
$R_{nr}(K,\mathbb{Z}/2\mathbb{Z})=  R'_{nr}(K,\mathbb{Z}/2\mathbb{Z})$, et l'ordre de $R_{nr}(K,\mathbb{Z}/2\mathbb{Z})$ est égal à
$2^{r^+_2(K)-m}$, où  $2^m$ est le degré de la composée de toutes les extensions quadratiques non ramifiées de la forme
$K(\sqrt{\varepsilon})/K$, où $\varepsilon$ est une unité de $K$.

\end{Th}

\noindent \textbf{Remarques.}

\textbf{1.} D'après $(2)$ et $(3)$ du Théorème \ref{fr-t1},  si $K$ est de type classe, alors $1\not \in R_{nr}(K,\mathbb{Z}/2\mathbb{Z})$, d'où $R_{nr}(K,\mathbb{Z}/2\mathbb{Z})$ n'est pas un sous-groupe de $Cl_K$.

\textbf{2.} On se place dans la situation (4) du Théorème \ref{fr-t1} en considérant une extension quadratique $K/\mathbb{Q}$, avec $K$ de type unité. Si $K$ est imaginaire, alors $m=1$ car $E_K=\{\pm 1\}$ (noter que $ \mathbb{Q}(i)= \mathbb{Q}(\sqrt{-1})$ et $\mathbb{Q}(j)= \mathbb{Q}(\sqrt{-3})$, où $j$ est une racine primitive 3ème de l'unité, sont de type classe car leur nombre de classes au sens restreint ($=1$) est impair). Si $K$ est réel et $\epsilon$ est son unité fondamental, comme  $E_K=\{\pm \epsilon^n, n\in \mathbb{Z}\}$, il est immédiat que  la composée dans (4) est contenue dans  l'extension biquadratique $K(\sqrt{-1}, \sqrt{\epsilon})$, et donc $m=1$ ou bien
$m=2$.

\textbf{3.} Dans la section \S3 nous classifierons les corps quadratiques sur $\mathbb{Q}$ de type classe ou de type unité et nous donnerons des exemples assez explicites pour le Théorème \ref{fr-t1}, en particulier des exemples de $K$ de type unité vérifiant aussi comme dans $(3)$ l'égalité $R'_{nr}(K,\mathbb{Z}/2\mathbb{Z}) =Cl_{K,2}$ (voir Proposition \ref{fr-exempletheoreme1.1unite}). $\blacksquare$

\begin{Th}\label{fr-t2} Soient $K$ un corps de nombres et $\Gamma$ un groupe fini. Si l'ordre de $\Gamma$ est impair, ou bien il a un $2$-sous-groupe de Sylow non cyclique,
ou bien  il a un $2$-sous-groupe de Sylow
 cyclique et normal, alors $ R_{nr}(K,\Gamma)\subset R_{nr}(K,\mathbb{Z}/2\mathbb{Z})$ et
 $ R'_{nr}(K,\Gamma)$ est un sous-groupe de $ R'_{nr}(K,\mathbb{Z}/2\mathbb{Z})$.
\end{Th}

\noindent \textbf{Remarque.} Une précision : d'après le Corollaire \ref{fr-cor3} (i) ci-dessous, si l'ordre de $\Gamma$ est impair ou bien il a un $2$-sous-groupe de Sylow non cyclique, alors $ R'_{nr}(K,\Gamma)=\{1\}$. $\blacksquare$.

  \begin{Cor} \label{fr-cort2} L'ensemble  $R'_{nr}(K,\Gamma)$ est un sous-groupe de
  $ R'_{nr}(K,\mathbb{Z}/2\mathbb{Z})$ dans les cas suivants :
   \begin{enumerate}
     \item $\Gamma$ est un $2$-groupe non abélien ;
     \item $\Gamma=S_n$ (le groupe symétrique de degré $n$) ou $A_n$ (le sous-groupe alterné de $S_n$), avec $n\geq 4$ ;
     \item $\Gamma$ est un groupe nilpotent fini.
   \end{enumerate}

\end{Cor}

\begin{proof}

 \textbf{1.} Le premier cas est évident.

\textbf{2.} On note (a \ b) une transposition de $S_n$.  Les groupes $S_n$ et $A_n$, avec $n\geq 4$,   contiennent le sous-groupe engendré par les deux  produits de transpositions suivants : $(1 \ 2)(3 \ 4)$ et $(1 \ 4)(2 \ 3)$. Clairement ce sous-groupe est isomorphe à $\mathbb{Z}/2\mathbb{Z}\times \mathbb{Z}/2\mathbb{Z}$. Par conséquent un $2$-sous-groupe de Sylow $H$ qui le contient ne peut être cyclique (les autres non plus car ils sont tous isomorphes à $H$).

(Noter que $ R'_{nr}(K,\Gamma)=\{1\}$ dans les cas  \textbf{1} et \textbf{2} d'après la remarque précédente.)

\textbf{3.} C'est clair si l'ordre de $\Gamma$ est impair ou  il a un $2$-sous-groupe de Sylow non cyclique. Si  $\Gamma$ a un $2$-sous-groupe de Sylow  cyclique alors il est normal car les  $p$-sous-groupes de Sylow d'un groupe nilpotent fini  sont normaux.
 \end{proof}

\noindent \textbf{Remarques.}

\textbf{1.} Soit $\Gamma= \mathbb{Z}/2^n\mathbb{Z }\times \Gamma' $, où $n\geq 0$ et $\Gamma'$ est un groupe fini d'ordre impair. Il est clair que l'ordre de $\Gamma$ est impair, ou bien  $\Gamma$  a un $2$-sous-groupe de Sylow
 cyclique et normal. Une autre conséquence du Théorème \ref{fr-t2} est : $ R'_{nr}(K,\Gamma)$ est un
  sous-groupe de $ R'_{nr}(K,\mathbb{Z}/2\mathbb{Z})$.

\textbf{2.}    Si $K$ est de type classe, alors $R_{nr}(K,\Gamma)$ n'est pas un sous-groupe de $Cl_{K, 2}$. En effet : $ R_{nr}(K,\Gamma)\subset R_{nr}(K,\mathbb{Z}/2\mathbb{Z})$ et  $1\not \in R_{nr}(K,\mathbb{Z}/2\mathbb{Z})$. $\blacksquare$

\section{Préliminaires}

Dans cette section, on donne ou on prouve des résultats qui seront utiles pour la démonstration des  Théorèmes \ref{fr-t1} et \ref{fr-t2}.

La proposition suivante est démontrée par Artin dans \cite{fr-artin} ;  pour   sa généralisation et  une preuve plus récente    on peut consulter \cite[Théorème I.4, p. 9]{fr-martinet}.  Elle permet de calculer la classe de Steinitz  d'une extension finie de corps de nombres.

\begin{Prop} (Théorème d'Artin) Soit $L/K$ une extension finie de corps de nombres.\label{fr-Artin}
Alors $$cl_{K}(\mathcal{O}_{L})=\left[\left(\frac{\delta_{L/K}}{d}\right)^{1/2}\right]_{K},$$
où $d$ est le discriminant d'une base du $K$-espace vectoriel $L$ et $\delta_{L/K}$ est le discriminant  de $L/K$.
\end{Prop}

\noindent \textbf{Remarques.}

\textbf{1.} On a $cl_{K}(\mathcal{O}_{L})^2=[\delta_{L/K}]_{K}$. En particulier, si $L/K$ est non ramifiée, alors $cl_{K}(\mathcal{O}_{L})^2=1$.

\textbf{2.} Si $L/K$ est galoisienne non ramifiée, alors $cl_{K}(\mathcal{O}_{L})$ capitule dans $Cl_L$ (i.e., si $cl_{K}(\mathcal{O}_{L})=[I]_{K}$, alors
$[I\mathcal{O}_L]_{L}=1$). En effet : $cl_{K}(\mathcal{O}_{L})=$ \break $[(d\mathcal{O}_K)^{-1/2}]_{K}$, et  $L/K$ est galoisienne implique que $d$ est un carré d'un élément de $L$ (voir le début de la démonstration de la proposition suivante).$\blacksquare$
\smallskip

La proposition suivante (voir \cite[Corollary, p. 4]{fr-endo}) se trouve dans la thèse de doctorat ("PhD Thesis") non publiée de L. P. Endo dont le directeur était L. R. McCulloh  . Elle est cruciale pour le présent article mais elle est intéressante en elle-même et pourra servir dans d'autres travaux, en particulier nos futures recherches. Pour la convenance du lecteur, nous donnons une démonstration plus précise et légèrement différente de celle se trouvant dans \cite{fr-endo}.

\begin{Prop}\label{fr-prop1}
Soient $L/K$ une extension galoisienne finie et $Gal(L/K)$ son groupe de Galois.

  $(i)$ Si  $Gal(L/K)$ est d'ordre impair ou il a un $2$-sous-groupe de Sylow non cyclique, alors
$$cl_{K}(\mathcal{O}_{L})=[(\delta_{L/K})^{1/2}]_{K}.$$

 $(ii)$  Si  $Gal(L/K)$ est d'ordre pair et a un  $2$-sous-groupe de Sylow  cyclique, alors $L/K$ admet une unique sous-extension quadratique $M/K$ et
$$cl_{K}(\mathcal{O}_{L})=\left[\left(\frac{\delta_{L/K}}{\alpha}\right)^{1/2}\right]_{K},$$
où $\alpha$ est un élément de $K$ tel que $M=K(\sqrt{\alpha})$.

\end{Prop}

\begin{proof} Soient $L/K$ une extension galoisienne de degré $n$.  Posons $Gal(L/K)=\{\sigma_i, 1\leq i \leq n\}$. Soient $\mathcal{B}=(a_i)_ {1\leq i \leq n}$ une base du $K$-espace vectoriel $L$ et $d$ son discriminant par rapport à la forme trace de $L/K$. Soit la matrice carrée  $M=(\sigma_i(a_j)) \in M_n(L)$. On sait que
$d= \det(M)^2 \in K$, où $\det$ désigne le déterminant. L'extension $L/K$ étant galoisienne, on a $\det(M)\in L$. L'idée de la démonstration est d'utiliser un raisonnement par disjonction des cas selon que $\det(M)\in K$ ou non, ce qui est équivalent à $K(\sqrt{d})=K$ ou bien $K(\sqrt{d} )/K$ est quadratique ;  équivalent aussi à tout $\tau \in Gal(L/K)$, $\tau(\det(M))=\det(M)$, ou bien il existe  $\tau \in Gal(L/K)$, $\tau(\det(M))\not =\det(M)$.
Nous divisons la preuve en 3 parties.

\textbf{1)} Si $n$ est impair, alors $\det(M)\in K$ ; sinon, $L/K$ aurait une sous-extension quadratique sur $K$, ce qui est impossible puisque $n$ est impair. Par le théorème d'Artin on obtient $cl_{K}(\mathcal{O}_{L})=[(\delta_{L/K})^{1/2}]_{K}$, d'où la première partie de $(i)$.

\textbf{2)} Soit $\tau \in Gal(L/K)$. Il est clair que  $\tau(\det(M))=\det(\tau\sigma_i(a_j))$. Notons $E_n$ l'ensemble des lignes $L_i:= (\sigma_i(a_j))_{ 1\leq j \leq n}$ de $M$. Puisque l'application de $G$ dans $G$ qui à $\sigma_i$ associe  $\tau\sigma_i$ est une bijection, il est immédiat que l'application $\tilde{\tau}$ de $E_n$ dans $E_n$ qui à $L_i$ fait correspondre $\tilde{\tau}(L_i):=(\tau\sigma_i(a_j))_{ 1\leq j \leq n}$ est une permutation de $E_n$ (autrement dit $\tilde{\tau} \in S(E_n)$ le groupe symétrique de $E_n$). Puisque $\tilde{\tau}$ consiste en une permutation des lignes de $M$,  immédiatement
$\tau(\det(M))=\epsilon(\tilde{\tau})\det(M),$
où $\epsilon$ désigne la signature.

 Pour calculer $\epsilon(\tilde{\tau})$,   nous allons maintenant décomposer $\tilde{\tau}$ en un produit de cycles de  $S(E_n)$.
Le sous-groupe $\langle\tilde{\tau}\rangle $ de $S(E_n)$ agit de façon naturelle sur $E_n$ de la manière suivante : $\tilde{\tau}^j.L_i=\tilde{\tau}^j(L_i)$ ($=(\tau^j\sigma_i(a_j))_{ 1\leq j \leq n}$) (noter que  $\widetilde{\tau^j}=\tilde{\tau}^j$).  Sous cette action l'orbite $O(L_i)$ de $L_i$ est égal à  $\{ \tilde{\tau}^j(L_i),  1\leq j \leq \vert \tau \vert -1   \}$.
  Le cardinal de $O(L_i)$ est $\vert \tau \vert$ l'ordre de $\tau$ ; en effet : soient $j$ et $j'$, avec $1\leq j, j' \leq \vert \tau \vert -1$. On a    $\tilde{\tau}^j(L_i)=\tilde{\tau}^{j'}(L_i)$ est équivalent à tout $l$,  $1\leq l \leq n$,
 $\tau^j\sigma_i(a_l)=\tau^{j'}\sigma_i(a_l)$ ; on en déduit que $\tau^j\sigma_i=\tau^{j'}\sigma_i$ car $\mathcal{B}$ est une base, et donc $j=j'$.
 Soit
 $\{L_{i_{k}}, 1\leq k \leq r \}$ un système de représentants des orbites. Alors $\tilde{\tau}=\prod_{k=1}^{r}c_k,$ où $c_k$ est le cycle de $S(E_n)$ de longueur
 $\vert \tau \vert$  suivant : $c_k=(\tilde{\tau}^j(L_{i_{k}}))_{ 1\leq j \leq \vert \tau \vert -1 } $.

  D'une part,  $\epsilon(c_k)=(-1)^{\vert \tau \vert -1}$. D'autre part $n=r\vert \tau \vert$ (les orbites $O(L_{i_{k}})$  forment une partition de $E_n$) . Il découle de $\epsilon(\prod_{k=1}^{r}c_k)=\prod_{k=1}^{r}\epsilon(c_k)$ que
  $\epsilon(\tilde{\tau})=(-1)^{\frac{n}{\vert \tau \vert}(\vert \tau \vert-1)}.$
  Notons tout de suite que si $\vert \tau \vert$ est impair, alors $\epsilon(\tilde{\tau})=1$.

  \textbf{3)}  Supposons que l'ordre de $Gal(L/K)$, noté $\vert Gal(L/K) \vert$, est pair et posons $\vert Gal(L/K) \vert=2^mk$ ($=n$), où $m\geq 1$ et $k$ est impair.
   Soit $\tau \in G$ d'ordre pair. On peut écrire  $\vert \tau \vert=2^{m'}{k'}$, où
   $1\leq m' \leq m$ et $k'$ est un diviseur de $k$, de sorte que  $\epsilon(\tilde{\tau})=(-1)^{2^{m-m'}}$.

  \textbf{3.1.} Supposons qu'un $2$-sous-groupe de Sylow de $Gal(L/K)$ soit non cyclique ;  les autres le sont aussi puisqu'ils sont ses  conjugués. Alors $m \not =m'$ ; sinon $\langle\tau\rangle $ aurait un $2$-sous-groupe cyclique d'ordre $2^m$, lequel serait un $2$-sous-groupe de Sylow cyclique de $Gal(L/K)$. Donc $\epsilon(\tilde{\tau})=1$. On conclut que $\det(M)\in K$. Par le théorème d'Artin, $cl_{K}(\mathcal{O}_{L})=[(\delta_{L/K})^{1/2}]_{K}$, d'où la deuxième partie de $(i)$.

  \textbf{3.2.} Supposons qu'un $2$-sous-groupe de Sylow de $Gal(L/K)$, noté\break $syl_2(Gal(L/K))$, soit cyclique (les autres le   sont aussi) et soit $\tau$ son générateur.
 Comme $\vert \tau \vert=2^{m}$, on a
  $\epsilon(\tilde{\tau})=-1$. Par suite $\det(M) \not \in K$ et $K(\sqrt{d})/K$ est quadratique. Montrons que c'est l'unique sous-extension quadratique de $L/K$.
  Si $K(\sqrt{\alpha})/K$ est quadratique et  différente de $K(\sqrt{d})/K$,   alors $K(\sqrt{d}, \sqrt{\alpha})/K$ est une extension biquadratique galoisienne, donc, par la théorie de Galois, il existe un sous-groupe  normal  $H$  de $Gal(L/K)$ tel que  $Gal(K(\sqrt{d}, \sqrt{\alpha})/K)\simeq Gal(L/K)/H$. Or le 2-sous-groupe de Sylow de $Gal(L/K)/H$ est $syl_2(Gal(L/K))H/H$ (voir ci-dessous :  Proposition \ref{fr-prp3} (1.2)), donc il serait cyclique, d'où une contradiction.

 L'égalité  $K(\sqrt{\alpha})=K(\sqrt{d})$ implique l'existence de $a\in K$ tel que $d=a^2\alpha$. On termine la démonstration de $(ii)$ en remplaçant $d$ par $a^2\alpha$
  dans le théorème d'Artin.
\end{proof}


\begin{Cor}\label{fr-cor3}
Soit $L/K$  une extension galoisienne finie non ramifiée.

\begin{itemize}
  \item [(i)] Si  $Gal(L/K)$ est d'ordre impair ou il a un $2$-sous-groupe de Sylow non cyclique, alors
$cl_{K}(\mathcal{O}_{L})=1$.
\item [(ii)]
   Si  $Gal(L/K)$ est d'ordre pair et a un  $2$-sous-groupe de Sylow  cyclique, alors $cl_{K}(\mathcal{O}_{L})=cl_{K}(\mathcal{O}_{M})$, où $M/K$  est l'unique sous-extension quadratique de $L/K$.
\end{itemize}
\end{Cor}

 Soit $K$ un corps de nombres. Nous allons  rappeler  la congruence  $\mod^*$ de la théorie du corps de classes pour un  module (ou cycle) $\mathfrak{M}$ qui est un idéal entier  de $\mathcal{O}_K$ différent de $\mathcal{O}_K$. (Dans le présent article nous n'avons pas besoin du cycle  $\mathfrak{M}_\infty$ défini par les places réelles infinies). Dans toute la suite $K^{\times}:=K\setminus \{0\}$. Si $\mathfrak{p}$ est un idéal premier de $\mathcal{O}_K$ divisant $\mathfrak{M}$, on note $\mathfrak{p}|\mathfrak{M}$.

  Soit $K^{\times}_\mathfrak{M}=\{x \in K^{\times}, \forall \  \mathfrak{p}| \mathfrak{M},
v_{\mathfrak{p}}(x-1)\geq v_{\mathfrak{p}}(\mathfrak{M})  \},$ où $v_{\mathfrak{p}}$ est la valuation $\mathfrak{p}$-adique.
On vérifie facilement que $K^{\times}_\mathfrak{M}$ est un sous-groupe du groupe $(K^{\times}, \times)$.
D'où une relation d'équivalence sur $K^{\times}$ compatible avec la multiplication :
$x\equiv y \mod K^{\times}_\mathfrak{M} \Leftrightarrow xy^{-1}\in K^{\times}_\mathfrak{M}.$ La notation usuelle de cette congruence dans la théorie du corps de classes est
$x\equiv y \mod^* \mathfrak{M}.$
 On vérifie sans difficulté qu'on a :
\begin{equation}\label{fr-1}
  x\equiv y \bmod^* \mathfrak{M} \Leftrightarrow \forall \    \mathfrak{p}|\mathfrak{M},  \ x\equiv y \bmod^* \mathfrak{p}^{ v_{\mathfrak{p}}(\mathfrak{M})},
\end{equation}
 \begin{equation}\label{fr-2}
 x\equiv y \mod^* \mathfrak{M}\Rightarrow \forall \   \mathfrak{p} | \mathfrak{M}, v_{\mathfrak{p}}(x)= v_{\mathfrak{p}}(y).
 \end{equation}

Rappelons la congruence usuelle modulo $\mathfrak{M}$ qui permet de définir l'anneau quotient $\mathcal{O}_K/\mathfrak{M}$ : soit $(x,y)\in \mathcal{O}_K^2$,    $$x\equiv y \mod   \mathfrak{M} \Leftrightarrow x-y \in \mathfrak{M} \Leftrightarrow \forall \    \mathfrak{p}|\mathfrak{M}, v_{\mathfrak{p}}(x-y)\geq v_{\mathfrak{p}}(\mathfrak{M}).$$

Soit $(\mathcal{O}_K/\mathfrak{M})^\times$ le groupe des éléments inversibles de l'anneau quotient $\mathcal{O}_K/\mathfrak{M}$ ;
 il est clair que $(\mathcal{O}_K/\mathfrak{M})^\times=\{\overline{x}\in \mathcal{O}_K/\mathfrak{M}, \forall \    \mathfrak{p}|\mathfrak{M}, v_{\mathfrak{p}}(x)=0  \}$.

 Soient $x$ et $y$  des éléments de $\mathcal{O}_K$ inversibles modulo $\mathfrak{M}$. Alors
 \begin{equation}\label{fr-3}
  x\equiv y \mod^* \mathfrak{M} \Leftrightarrow x\equiv y \mod   \mathfrak{M},
 \end{equation}
car si $\mathfrak{p}|\mathfrak{M}$, alors $v_{\mathfrak{p}}(\frac{x}{y}-1)=v_{\mathfrak{p}}(x-y)$, puisque
$v_{\mathfrak{p}}(\frac{x-y}{y})=v_{\mathfrak{p}}(x-y)-v_{\mathfrak{p}}(y)$ et $v_{\mathfrak{p}}(y)=0$.

Soit $x\in \mathcal{O}_{K, \mathfrak{p}}$ le localisé de $\mathcal{O}_K$ en $\mathfrak{p}$. Alors, pour  tout  entier naturel non nul $n$
\begin{equation}\label{fr-4}
 \exists \ a \in \mathcal{O}_K, \ \  x\equiv a \bmod^* \mathfrak{p}^{ n}.
 \end{equation}
En effet : écrivons $x=\frac{b}{c}$, où $b \in  \mathcal{O}_K$ et $c \in \mathcal{O}_K \setminus \mathfrak{p}$ ; $v_{\mathfrak{p}}(c)=0$ implique $\overline{c}\in (\mathcal{O}_K/\mathfrak{p}^{ n})^\times$. Soit $\overline{d}$ l'inverse de $\overline{c}$. Posons $a=bd$ ; on a $\frac{a}{x}-1 =cd-1 \in  \mathfrak{p}^{ n}$, d'où  $a\equiv x \bmod^* \mathfrak{p}^{ n}$ et donc $x\equiv a \bmod^* \mathfrak{p}^{ n}$ (par la symétrie de $\equiv$).

\begin{Prop}\label{fr-unitemodulocarre} Soient $K$ un corps de nombres,  $\varepsilon$  une unité de $K$ et $\mathfrak{M}$ un module de $K$. Alors les deux assertions suivantes sont équivalentes.

(i) Il existe $x\in K ^\times $ tel que $\varepsilon \equiv x^2 \bmod^{*} \mathfrak{M}$.

 (ii)  Il existe $a\in \mathcal{O}_K$ tel que $\varepsilon \equiv a^2 \bmod \mathfrak{M}$.

 \end{Prop}
 \begin{proof} Remarquons que d'après  (\ref{fr-3}), $(ii) \Leftrightarrow$  Il existe $a\in \mathcal{O}_K$ tel que $\varepsilon \equiv a^2 \bmod^* \mathfrak{M}$
  (puisque $\overline{\varepsilon}$ et $ \overline{a^2}\in (\mathcal{O}_K/\mathfrak{M})^\times$).

(1) On a  (ii) implique (i) par la remarque précédente.

(2) Supposons (i). Soit $\mathfrak{p}$ un idéal premier tel que $\mathfrak{p}|\mathfrak{M}$. Il suit de  (\ref{fr-2}) que $v_{\mathfrak{p}}(\varepsilon)= v_{\mathfrak{p}}(x^2)=0$. D'où $v_{\mathfrak{p}}(x)=0$ et donc $x\in \mathcal{O}_{K, \mathfrak{p}}$. L'assertion (\ref{fr-4}) nous donne :
$\exists \ a_\mathfrak{p} \in \mathcal{O}_K, \ \  x\equiv a_\mathfrak{p} \bmod^* \mathfrak{p}^{v_{\mathfrak{p}}(\mathfrak{M})}$. Le théorème des restes chinois nous permet d'affirmer : $\exists \ a \in \mathcal{O}_K, \ \  a\equiv a_\mathfrak{p} \bmod \mathfrak{p}^{v_{\mathfrak{p}}(\mathfrak{M})}$ pour tout $\mathfrak{p}|\mathfrak{M}$ ; on a $a\equiv a_\mathfrak{p} \bmod^* \mathfrak{p}^{v_{\mathfrak{p}}(\mathfrak{M})}$ par  (\ref{fr-3}), car $v_{\mathfrak{p}}(x)=v_{\mathfrak{p}}(a_\mathfrak{p})=0$, et donc $x\equiv a_\mathfrak{p} \bmod^* \mathfrak{p}^{v_{\mathfrak{p}}(\mathfrak{M})}$. D'après (\ref{fr-1}), $x\equiv a \bmod^* \mathfrak{M}$. Par conséquent $\varepsilon \equiv a^2 \bmod^* \mathfrak{M}$, ce qui termine la démonstration.
 \end{proof}

En  utilisant la théorie de Kummer (voir \cite[\S39, pp. 132-139]{fr-hecke}) et le théorème d'Artin, on démontre facilement la proposition suivante
(signalons, par exemple, que c'est le cas $p=2$ dans \cite[Proposition 3.2]{fr-BrS}).
\begin{Prop}\label{fr-th2}
  Soient $L/K$ une extension quadratique et $\alpha\in K^{*}$ tel que $L=K(\sqrt{\alpha})$.
  \begin{enumerate}
   \item L'extension L/K  est modérément ramifiée si et seulement s'il existe $x\in K^{\times}$ tel que  $\alpha \equiv x^2 \mod^{*}  4\mathcal{O}_{K}.$
   \item Supposons  $L/K$  modérément ramifiée. Alors
\begin{enumerate}
\item On peut écrire d'une manière unique $\alpha \mathcal{O}_{K}=I(\alpha)^{2}J$, où $I(\alpha)$ est un idéal fractionnaire de $\mathcal{O}_{K}$, et $J$ est un idéal entier de $\mathcal{O}_{K}$ sans facteur carré ;
\item  $\delta_{L/K}=J$ ;
\item  $cl_{K}(\mathcal{O}_{L})=[I(\alpha)]_{K}^{-1}$.
 \end{enumerate}
\end{enumerate}
\end{Prop}

\begin{proof}Pour la convenance du lecteur on démontre l'assertion 1 brièvement.
L'extension $L/K$  est modérément ramifiée si et seulement si tout idéal premier $\mathfrak{p}$ de $\mathcal{O}_K$ au dessus de   $2$ est non ramifié. Comme $L/K$ est une extension de Kummer de degré $2$, si $v_{\mathfrak{p}}(\alpha)=0$ pour tout $\mathfrak{p}$  au dessus   de $2$,  le Théorème 119 de Hecke (voir \cite[p. 136]{fr-hecke}) nous donne 1. Sinon, la condition $\alpha \equiv x^2 \mod^{*}  4\mathcal{O}_{K}$ est équivalente à $\alpha  x^{-2} \equiv 1 \mod^{*}  4\mathcal{O}_{K}$. Posons $\alpha'=\alpha  x^{-2}$, alors $L=K(\sqrt{\alpha'})$ et pour tout $\mathfrak{p}$  au dessus   de $2$ on a $v_{\mathfrak{p}}(\alpha')=0$ ; appliquer le Théorème 119 avec $\alpha'$ nous permet de conclure.
\end{proof}

Comme toute extension de corps de nombres non ramifiée  est modérée, on obtient :

\begin{Cor}\label{fr-cor1}
 Soit $L=K(\sqrt{\alpha})$  une extension quadratique de $K$. Alors,  $L/K$ est  non ramifiée si et seulement s'il existe $x\in K^{\times}$ tel que  $\alpha \equiv x^2 \mod^{*}  4\mathcal{O}_{K}$ et il existe un idéal fractionnaire $I$ de $K$ tel que $\alpha \mathcal{O}_{K}=I^{2}$ ;
dans ce cas   $cl_{K}(\mathcal{O}_{L})=[I]_{K}^{-1}=[I]_{K}$.
\end{Cor}

\noindent \textbf{Remarque.} En fait, la première partie de ce corollaire  est le Théorème 120 dans \cite[p. 137]{fr-hecke} dans le cas d'une extension quadratique. $\blacksquare$

\begin{Cor}\label{fr-cor4}
 Soit $N/K$  une extension biquadratique non ramifiée. Soient  $L_i/K$, $1\leq i \leq 3$, ses trois sous-extensions quadratiques. Alors $$cl_{K}(\mathcal{O}_{L_3})=cl_{K}(\mathcal{O}_{L_1})cl_{K}(\mathcal{O}_{L_2}).$$
\end{Cor}

\begin{proof}
  C'est une conséquence immédiate de \cite[Proposition 3.2 (1)]{fr-bouchaib}.
\end{proof}

On fixe les notations suivantes pour toute la suite :   $p$ est un nombre premier et  si  $G$ est un groupe fini   on  désigne par $Syl_p(G)$ l'ensemble de ses $p$-sous-groupes de Sylow. On note $syl_p(G)$  l'un des éléments de ce dernier (unique à conjugaison près)  de sorte que les autres sont ses conjugués : $Syl_p(G)=\{gsyl_p(G)g^{-1}, g\in G\}$.
Parfois en vue de simplification,  nous ferons un abus de notation en notant $syl_p(G)$ un conjugué de $syl_p(G)$. Si $H$ est un sous-groupe (resp. sous-groupe normal) de $G$, on notera $H\leq G$ (resp. $H\trianglelefteq G$). L'ordre de $G$ (resp. un élément $g$ de $G$) sera noté $\vert G \vert$ (resp. $\vert g\vert$).

Nous aurons besoin des résultats bien connus suivants :
\begin{Prop}\label{fr-prp3}
Soient $G$ et $G'$ des groupes finis.

$(1)$ Si $H\trianglelefteq G$, alors

$(1.1)$ $Syl_p(H)=\{P\cap H, P\in Syl_p(G) \}$ ;

$(1.2)$ $Syl_p(G/H)=\{PH/H, P\in Syl_p(G) \}$.

$(2)$ $Syl_p(G\times G')=Syl_p(G)\times Syl_p(G')$.

$(3)$ Supposons  $PGCD(\vert G \vert, \vert G' \vert)=1$. Alors,
$K\leq G\times G'$ si et seulement s'il existe $H\leq G$ et $H'\leq G'$ tel que $K=H\times H'$.

$(4)$ Soient $H\trianglelefteq G$ et $K\trianglelefteq G$. Alors, $G$ est un produit direct de $H$ et $K$  si et seulement si c'est un produit semi-direct de $H$ et $K$.

\end{Prop}

\begin{proof} Pour (1.1), (1.2) et (2) on pourrait consulter  \cite[ Corollary 6.4, Theorem 6.6,  Theorem 7.1]{fr-CO}.
 Pour (3) voir \cite[Theorem 2]{fr-AN}. La preuve de (4) ne présente aucune difficulté.
\end{proof}

\begin{Prop}\label{fr-Scott} Soit $G$ un groupe d'ordre pair $2^nk$, où $k$ est impair. Si $G$ a un $2$-sous-groupe de Sylow cyclique, alors $G$ admet un unique sous-groupe normal H d'ordre $k$.  De plus, $H$ est l'ensemble des éléments de $G$ d'ordres impairs.
\end{Prop}

\begin{proof}
Voir \cite[ Theorem 6.8, p. 14]{fr-CO} ou \cite[6.2.11, p. 138]{fr-SC}.

\end{proof}

\noindent \textbf{Remarques}.

\textbf{1.} Sous les notations et les hypothèses de la  Proposition \ref{fr-Scott}, il est immédiat que $G= H \rtimes syl_2(G)$ et donc $G/H \simeq syl_2(G)$.
Si $syl_2(G)$ est normal, alors $G$ est un produit direct de $H$ et $syl_2(G)$ par la Proposition \ref{fr-prp3} (4). On en déduit : $G$ admet un $syl_2(G)$ normal si et seulement s'il existe un  groupe fini $G'$ d'ordre $k$ tel que $G\simeq \mathbb{Z}/2^n\mathbb{Z }\times G'$.

\textbf{2.} Soit $p$ un nombre premier ;  dans la théorie des groupes, un groupe fini $\Gamma$ est dit $p$-nilpotent s'il admet un sous-groupe normal $N$ et un $p$-sous-groupe de Sylow  tels que $G=N\rtimes syl_p(G)$. D'après la remarque \textbf{1} précédente, le groupe $G$ de la  Proposition \ref{fr-Scott} est $2$-nilpotent (ce fait a été démontré pour la première fois par Burnside). De plus, $G$ est résoluble car $H$ et $G/H$ le sont,   en effet : $|H|$ étant impair, d'après le théorème de Feit-Thompson $H$ est résoluble ; $G/H$ est résoluble car il est abélien. Donc,  si $K$ est un corps de nombres, d'après un  théorème de Shafarevich $G$ est réalisable comme groupe de Galois sur $K$ (et même d'une façon modérément ramifiée).
$\blacksquare$
\smallskip

Pour la   définition du $p$-rang d'un groupe fini, nous avons suivi celle donnée dans \cite[Chap. IX, pp. 232-233]{fr-CAFR}, elle généralise la définition dans la théorie des groupes du $p$-rang d'un groupe abélien et elle est  plus pratique pour la théorie de Galois comme on va le voir ci-dessous.

\begin{Df}  \label{fr-rang}
Soit $G$ un groupe fini. Soit $A$ l'ensemble des $n\in \mathbb{N}$   tels qu'il existe un sous-groupe normal $H$ de $G$ tel que $G/H$ est un $p$-groupe élémentaire abélien d'ordre $p^n$ (donc $n=\dim_{\mathbb{F}_p}(G/H)$, où $\mathbb{F}_{p}$ est le corps fini à $p$ éléments). Le $p$-rang de $G$, noté $r_{p}(G)$, est le plus grand élément de $A$.
\end{Df}

\noindent {\bf Remarque.}  Dans la théorie des groupes, le $p$-rang de $G$ est défini comme étant le plus grand élément de l'ensemble des $n\in \mathbb{N}$   tels qu'il existe un sous-groupe  H de $G$ satisfaisant $H$ est un $p$-groupe élémentaire abélien d'ordre $p^n$. On vérifie sans peine que cette définition et  Définition  \ref{fr-rang} sont équivalentes  lorsque $G$ est abélien ;  donc $r_{p}(G)$ est égal au rang du $p$-sous-groupe élémentaire $G_p:=\{\ g\in G, g^p=1\}$ de $G$. $\blacksquare$

 \begin{Prop}\label{fr-2rang} Soient  $G$ un groupe fini, $G^2=\{g^2, g\in G\}$ et $\langle G^2\rangle $ le sous-groupe de $G$ engendré par $G^2$. Alors

  (1)  $\langle G^2\rangle \trianglelefteq G$, le groupe  quotient $G/\langle G^2\rangle $ est un $2$-groupe élémentaire (abélien) et  $r_{2}(G)=\dim_{\mathbb{F}_2}(G/\langle G^2\rangle )$. (Noter que $\langle G^2\rangle =G^2$ si $G$ est abélien.)

  (2) Le nombre de sous-groupes de $G$ d'indice $2$ est égal à $2^{r_{2}(G)}-1$.

   (3) Si $G'$ est un groupe fini, alors $r_{2}(G\times G')=r_{2}(G)+r_{2}(G')$.

  (4) Si $G$ est  d'ordre pair ayant un  $2$-sous-groupe de Sylow cyclique, alors $r_{2}(G)=1$.

   (5) Soit $H \trianglelefteq G$. Alors  $r_{2}(G/H)\leq  r_{2}(G)$.

 \end{Prop}
 \begin{proof}
 (1) Pour tout $g\in G$ et  $h^2\in G^2$, on a  $gh^2g^{-1}=(ghg^{-1})^2 \in G^2$ ; on en déduit que $\langle G^2\rangle $ est normal.

 Tout élément de $G/\langle G^2\rangle $ est d'ordre $1$ ou $2$, par suite $G/\langle G^2\rangle $ est un $2$-groupe élémentaire abélien, disons d'ordre $2^m$.

Soit $H$  un sous-groupe normal de $G$ satisfaisant $G/H$ est un $2$-groupe élémentaire abélien d'ordre $2^n$. Pour tout classe $\overline{g}$ de $G/H$, $\overline{g}^2=\overline{1}$, d'où $g^2\in H$ ; par conséquent   $\langle G^2\rangle \subset H$, et donc $n\leq m$. On conclut que $r_{2}(G)=m$. Il découle de   $G/\langle G^2\rangle $ est un
 $\mathbb{F}_2$-espace vectoriel de cardinal $2^m$ que  $\dim_{\mathbb{F}_2}(G/\langle G^2\rangle )=m$, ce qui termine la preuve de (1).

 (2) D'une part, le nombre des sous-groupes de $G/\langle G^2\rangle $ d'indice $2$ est égal au nombre d'hyperplans du $\mathbb{F}_2$-espace vectoriel $G/\langle G^2\rangle $. Il est bien connu que ce nombre est $2^{\dim_{\mathbb{F}_2}(G/\langle G^2\rangle )}-1$, donc égal à $2^{r_{2}(G)}-1$.

  D'autre part il y a une bijection entre l'ensemble des  sous-groupes normaux $H$ de  $G$ contenant $\langle G^2\rangle $ et celui des  sous-groupes normaux de\break  $G/\langle G^2\rangle $ : elle fait correspondre $H/\langle G^2\rangle $ à $H$. On a  $(G/\langle G^2\rangle )/(H/\langle G^2\rangle )\simeq (G/H)$, d'où $H$ est d'indice $2$ dans $G$ si et seulement si $H/\langle G^2\rangle $ est d'indice $2$ dans $G/\langle G^2\rangle $. Donc le nombre de sous-groupes de $G$ d'indice $2$ est égal à celui des sous-groupes de $G/\langle G^2\rangle $ d'indice $2$.

(3) On a $(G\times G')^2=G^2\times G'^2$. Immédiatement $\langle G^2\times G'^2\rangle  \subset \langle G^2\rangle \times \langle G'^\rangle $.
Puisque $1\in G^2$, on a $\langle\{1\}\times G'^2\rangle  \subset \langle G^2\times G'^2\rangle $ ; facilement on a $\langle\{1\}\times G'^2\rangle =\{1\}\times \langle G'^\rangle $ et  de même $\langle G^2\rangle \times \{1\}\subset \langle G^2\times G'^2\rangle $, par conséquent $\langle G^2\rangle \times \langle G'^\rangle  \subset \langle G^2\times G'^2\rangle $. On conclut que $\langle G^2\times G'^2\rangle  =\langle G^2\rangle \times \langle G'^\rangle =\langle(G\times G')^2\rangle $, ce qui implique que  $(G\times G')/\langle(G\times G')^2\rangle  \simeq (G/ \langle G^2\rangle )\times (G'/ \langle G'^\rangle )$. Un calcul simple de dimension sur $\mathbb{F}_2$ (en utilisant (1)) nous donne $(3)$.

(4) Soit $K\unlhd G$ tel que $G/K \simeq (\mathbb{Z}/2\mathbb{Z})^m$. Comme $syl_2(G)$ est cyclique, il découle  de la Proposition \ref{fr-prp3} (1.2) que $syl_2(G/K)$ est cyclique, d'où $m\leq 1$ et donc $r_{2}(G)\leq 1$.
Dans les notations de la Proposition \ref{fr-Scott} et d'après la remarque qui la suit, $G/H$ est cyclique d'ordre $2^n$, où $n\geq 2$, donc admet un sous-groupe $K'/H$ d'indice $2$. Comme $(G/H)/(K'/H)
\simeq (G/K')$, on a  $K'$ est un sous-groupe de $G$ d'indice 2 et donc $G/K' \simeq (\mathbb{Z}/2\mathbb{Z})$, d'où $r_{2}(G)\geq 1$. On conclut que $r_{2}(G)=1$.

(5) Soit $H_1$ un sous-groupe normal  de $G$ tel que $G/H_1\simeq (\mathbb{Z}/2\mathbb{Z})^{r_{2}(G)}$. Soit $H_2$ un sous-groupe normal  de $G$ contenant $H$ tel que $(G/H)/(H_2/H)\simeq (\mathbb{Z}/2\mathbb{Z})^{r_{2}(G/H)}$. Comme $(G/H)/(H_2/H)\simeq G/H_2$, on obtient $r_{2}(G/H)\leq  r_{2}(G)$ par maximalité de $r_{2}(G)$.
\end{proof}

 Nous fixons la notation suivante pour toute la suite de l'article : $r_2(K)=r_2(Cl_{K})$,   et nous rappelons que $r^+_2(K)=r_2(Cl^+_{K})$ (introduite dans $\S1$) ; dans la littérature, ils sont appelés respectivement le $2$-rang de $K$ et le $2$-rang de $K$ au sens restreint.

 On note $Cl^+_{K,2}$ l'ensemble des $c\in Cl^+_K$ telles que $c^2=1$. Rappelons que $Cl_{K,2}$ est l'ensemble des $c\in Cl_K$ telles que $c^2=1$ (introduit dans $\S1$).
 Signalons que d'après la remarque qui suit la Définition \ref{fr-rang}, on a  $$r_2(K)=r_2(Cl_{K,2}), \ \   r^+_2(K)=r_2(Cl^+_{K,2}).$$

Soient  $k$ un corps de nombres et  $K/k$ une extension quadratique. Soit $u$ (resp. $u_{\infty}$) le nombre des idéaux premiers (resp. places infinies) de $k$ ramifiés (resp. ramifiées) dans $K/k$. Soit $E_k$ (resp. $E^+_k$) le groupe des unités (resp. unités totalement positives) de $k$. Soient $e$ et  $e^+$  les entiers satisfaisants  $$2^e=[E_k : E_k\cap N_{K/k}(K^\times)], \   2^{e^+}=[E^+_k : E^+_k\cap N_{K/k}(K^\times)],$$ où  $N_{K/k}$ est la norme dans $K/k$ .

On désigne par $Am(K/k)$ (resp. $Am^+(K/k)$) le sous-groupe de $Cl_K$ (resp. $Cl^+_K$) des classes ambiges de l'extension $K/k$ (i.e., les éléments de $Cl_K$ (resp. $Cl^+_K$)  fixes par $Gal(K/k)$).

\begin{Prop} \label{fr-K/kquadratique} Sous les notation précédentes, on a les assertions suivantes.

  (i) Supposons $h_k$ impair. Alors $syl_{2}(Am(K/k))=Cl_{K,2}$ et $r_2(K)=u+u_{\infty} - e -1$.

  (ii) Supposons $h^+_k$ impair. Alors $syl_{2}(Am^+(K/k))=Cl^+_{K,2}$ et $r^+_2(K)=u - e^+ -1$.

  (iii) Si  $h_k=1$ (resp. $h^+_k=1$),  alors $Am(K/k)=Cl_{K,2}$ (resp. $Am^+(K/k)$\break $=Cl^+_{K,2}$).
\end{Prop}

\begin{proof}
  (i) Soit  $c\in Cl_K$ ; si l'ordre de  $c$ est une puissance de $2$, il en est de même pour  $N_{K/k}(c)$ dans $Cl_k$, ce qui nous donne  $ N_{K/k}(c)=1$ dans $Cl_k$ puisque  $h_k$  est impair; en étendant les idéaux de $\mathcal{O}_k$ représentants de classes  dans $Cl_k$ à $\mathcal{O}_K$  on  en déduit   $c\sigma(c)=1$ dans $Cl_K$, d'où $\sigma(c)=c^{-1}$. Donc
    $c^2=1 \Leftrightarrow \sigma(c)=c$, d'où $syl_{2}(Am(K/k))=Cl_{K,2}$.

   D'après la formule bien connue des classes ambiges de C. Chevalley (voir \cite[a, p. 406]{fr-CH}) on a
   $$| Am(K/k)| =h_k2^{u+u_{\infty}-e-1}.$$ Si $h_k$ est impair, alors $|syl_{2}(Am(K/k))|=$ $2^{u+u_{\infty}-e-1}$, d'où la formule (connue, voir par exemple \cite[(1.1), p. 27]{fr-SORO}) de $r_2(K)$.

   (ii) La démonstration de $syl_{2}(Am^+(K/k))=Cl^+_{K,2}$ est analogue à la   première partie de  celle de (i).

      La formule des classes ambiges au sens restreint de   G. Gras (voir \cite[Théorème 4.1, p. 26 ]{fr-gras}) donne  $$| Am^+(K/k)| =h^+_k 2^{u-e^+-1}.$$ Si $h^+_k$ est impair, comme dans (i) on obtient $r^+_2(K)$.

     (iii) Soit $c\in Cl_K$ (resp. $Cl^+_{K}$). On a $ N_{K/k}(c)=1$, d'où $\sigma(c)=c^{-1}$ et donc $c^2=1 \Leftrightarrow \sigma(c)=c$.
     \end{proof}

 Le corollaire suivant est une conséquence immédiate de la proposition précédente en prenant $k=\mathbb{Q}$ ; il résume des résultats bien connus depuis l'époque de Hilbert et Hasse.

 \begin{Cor}\label{fr-corK/kquadratique}
   Soit $K/\mathbb{Q}$ une extension quadratique. Sous les notations précédentes, on a :

   (1) $r^+_2(K)=u-1$.

   (2) Si $K$ est imaginaire, alors $r_2(K)=u-1$.

   (3) Supposons $K$ réel et soit $\epsilon$ son unité fondamental.

    (a) Si $N_{K/\mathbb{Q}}(\epsilon)=-1$, alors $r_2(K)=u-1$.

    (b) Supposons $N_{K/\mathbb{Q}}(\epsilon)=1$.

   (b.1) Si $-1$ est une  norme dans $K/\mathbb{Q}$, alors $r_2(K)=u-1$.

   (b.2)  Si $-1$ n'est pas une norme dans $K/\mathbb{Q}$, alors $r_2(K)=u-2$.
   \end{Cor}

\noindent \textbf{Remarque.} Sachant que $h^+_K$ (resp. $h_K$) est impair si et seulement si $r^+_2(K)=0$ (resp. $r_2(K)=0$), ce corollaire permet d'expliciter  toutes les extensions quadratiques $K/\mathbb{Q}$ ayant $h^+_K$ (resp. $h_K$) impair. $\blacksquare$
 \smallskip

Concernant les conditions dans (b.1) et (b.2), on peut les expliciter grâce à  un théorème de Hasse (local-global) sur les normes disant que $-1$ est une  norme dans $K/\mathbb{Q}$ si et seulement si $-1$ est une   norme locale pour toute place (finie ou infinie). Dans le cas $K=\mathbb{Q}(\sqrt{d})$ réel, où $d$ est un entier sans facteur carré,  ceci est équivalent à tout $p\not =2$ premier divisant $d$,  on a $(\dfrac{-1}{p})=1$ (où $(\dfrac{a}{p})$ est le symbole de Legendre), ce qui est équivalent à $p \equiv 1 \mod 4$ (car $(\dfrac{-1}{p})=(-1)^{\frac{p-1}{2}}$), équivalent aussi à $d$ est somme de deux carrés ($\Leftrightarrow$ $K/\mathbb{Q}$ est plongeable dans une extension cyclique de degré 4 sur $\mathbb{Q}$).

   D'après \cite[I, p. 94]{fr-FRMA}, si $K$ est un corps de nombres ayant $r$ plongements réels, alors $r^+_2(K)-r_2(K)\leq E(\frac{r}{2})$, où $E$ désigne la partie entière.  Donc, si $K/\mathbb{Q}$ est quadratique réelle, alors $$r^+_2(K)= r_2(K)   \text{ ou bien} \  r^+_2(K)= r_2(K)+1.$$
    Le corollaire précédent nous dit que c'est  seulement dans  (b.2) qu'on a $r^+_2(K)= r_2(K)+1$. On a $ \ -1 \  \text{n'est pas une norme dans}\  K/\mathbb{Q}$ si et seulement s'il existe un diviseur premier $p$ de $d$ tel que $p \equiv 3 \mod 4$. Cette dernière condition donne  $N_{K/\mathbb{Q}}(\epsilon)=1$ grâce à  \cite[(1.7), p. 179]{fr-FT}.
   On a donc  la caractérisation suivante.

  \begin{Cor}\label{fr-r^+_2(K)= r_2(K)+1} Soit $K/\mathbb{Q}$ une extension quadratique réelle. Alors,
  $r^+_2(K)= r_2(K)+1  \Leftrightarrow$ Il existe un diviseur premier $p$ de $d$ tel que $p \equiv 3 \mod 4.$
  \end{Cor}

 \noindent \textbf{Remarque} On peut trouver une autre démonstration de ce corollaire (un peu plus compliquée) dans \cite[Remarque 2 ]{fr-Dummitquebec}. $\blacksquare$ \smallskip

  Soit $K/\mathbb{Q}$ une extension quadratique. Le but de ce qui suit est de fournir  une description explicite de  $Cl_{K,2}$.

   En suivant Hilbert (voir \cite[\S 75,  p. 136]{fr-HI}) un idéal de $\mathcal{O}_K$ est dit ambige s'il est fixe par $Gal(K/\mathbb{Q})$ et ses seuls diviseurs dans $\mathbb{Z}$ sont $\pm1$. Le Théorème 105 de Hilbert (voir \cite[ p. 136]{fr-HI}) nous dit qu'il y  a exactement $u$ idéal premiers ambiges qui sont les idéaux premiers au dessus des premiers ramifiés dans $K/\mathbb{Q}$, en les multipliant entre eux ils donnent exactement $2^u$ idéaux ambiges.  Notons $Am_0(K/\mathbb{Q})$ le sous-groupe de $Am(K/\mathbb{Q})$ formé par les éléments contenant un idéal ambige. Le Théorème 106 de Hilbert (voir \cite[ p. 136]{fr-HI}) affirme  :

  - Si $K$ est imaginaire ou réel avec  $N_{K/\mathbb{Q}}(\epsilon)=-1$, alors  $$|Am_0(K/\mathbb{Q})|=2^{u-1}.$$

 - si  $K$ est réel avec $N_{K/\mathbb{Q}}(\epsilon)=1$, alors $$|Am_0(K/\mathbb{Q})|=2^{u-2}.$$

\noindent \textbf{Remarques.}

\textbf{1.} Le groupe $Am_0(K/\mathbb{Q})$ est engendré par les classes des idéaux premiers de $\mathcal{O}_K$ au dessus des  premiers ramifiés dans $K/\mathbb{Q}$.

\textbf{2.} La définition générale d'un idéal fractionnaire ambige $I$ de $K$ est :  $I$ est  fixe par $Gal(K/\mathbb{Q})$ ;  en fait
$Am_0(K/\mathbb{Q})$ est  le sous-groupe de $Am(K/\mathbb{Q})$ formé par les classes de tous les idéaux ambiges. En effet : une telle classe contient un idéal entier ambige $\mathfrak{a}$ ; d'après \cite[(1.16), p. 182]{fr-FT}, il existe $r\in \mathbb{Q}$ et $\mathfrak{q}$ un idéal ambige au sens de Hilbert tel que $\mathfrak{a}=r\mathfrak{q}$. $\blacksquare$
\smallskip

 La combinaison du  Théorème 106 de Hilbert  et le Corollaire \ref{fr-corK/kquadratique} aboutit à la description de $Cl_{K,2}$ suivante.

\begin{Cor} \label{fr-descclasseK2} Soit $K/\mathbb{Q}$ une extension quadratique.

  - Si $K$ est imaginaire ou réel avec  $N_{K/\mathbb{Q}}(\epsilon)=-1$, alors
 $$Cl_{K,2}=Am_0(K/\mathbb{Q}), \ \ |Cl_{K,2}|=2^{u-1}.$$

 Supposons $K$   réel avec  $N_{K/\mathbb{Q}}(\epsilon)=1$.

 -  Si $-1$ n'est pas norme dans $K/\mathbb{Q}$ , alors
 $$Cl_{K,2}=Am_0(K/\mathbb{Q}), \ \ |Cl_{K,2}|=2^{u-2}.$$

 - Si $-1$ est  une norme dans $K/\mathbb{Q}$, alors
$$ Cl_{K,2}=Am(K/\mathbb{Q}), \ \ |Cl_{K,2}|=2^{u-1},  \ \ [Am(K/\mathbb{Q}): Am_0(K/\mathbb{Q})]=2.$$
\end{Cor}

  \noindent \textbf{Remarque.} Une légère  correction :  le Corollaire 1 dans \cite[Corrolaire 1, p. 181]{fr-FT} est vrai seulement si l'on remplace dans son énoncé et ses notations $C_{K,2}$ (qui n'est pas défini, mais vu la définition du $C^+_{K,2}$ dans \cite[Theorem 39, p. 181]{fr-FT}, il semble que c'est notre  $C_{K,2}$)   par $Am_0(K/\mathbb{Q})$. En fait ce corollaire est le Théorème 106 de Hilbert. Aussi, dans p. 181, à remplacer dans la définition de  $\nu_K$, $C_K$ par $C^+_K$. $\blacksquare$

\section{Classes de Steinitz des extensions quadratiques non ramifiées à groupe de Galois $\Gamma\simeq\mathbb{Z}/2\mathbb{Z}$}

Soient $l$ un nombre premier et $K$ un corps de nombres contenant une racine primitive $l$ème de l'unité. Soit
$$V_l(K)=\{\alpha, \text {\textit{Il existe un idéal fractionnaire}} \ I \ \text {\textit{de}} \  K \  \text {\textit{tel que} } \ I^l=\alpha \mathcal{O}_K \}.$$
 Dans \cite[Définition 5.2.4 (3), p. 231]{fr-cohan}, l'auteur a introduit la définition du $l$-Selmer groupe de $K$, noté $S_l(K)$,  par analogie aux groupes de Selmer des courbes elliptiques (voir \cite[Remarque, p. 233]{fr-cohan} ; selon l'auteur, les éléments de $V_l(K)$ sont appelés des unités virtuelles, d'où la notation $V$) :
$$S_l(K)=V_l(K)/(K^\times)^l.$$
On pourrait consulter \cite{fr-Dummitlondon,fr-LE} pour des investigations concernant $S_2(K)$ en lien avec le rang au sens restreint d'un corps de nombres et les lois de réciprocité quadratiques.\smallskip

\noindent \textbf{Remarque.} D'après la théorie de Kummer, plus précisément \cite[Théorème 118, p. 134]{fr-hecke}, il y a une bijection entre   $ S_l(K)\setminus \{\bar {1}\}$ et l'ensemble  des extensions de Kummer
 $L/K$ de degré $l$, non ramifiées à toutes les places finies en dehors de celles  au dessus de $l$ ; elle fait correspondre à la classe $\bar{\alpha}$ l'extension  $L=K(\sqrt[l]{\alpha})$. $\blacksquare$
\smallskip

 Soient $K$ un corps de nombres quelconque  et $R_{K}$ le sous-groupe de  $K^\times$ formé par les éléments $\alpha$
 vérifiant :
\begin{enumerate}
\item \textbf{ il existe $x\in K^\times$ tel que  $\alpha\equiv x^2 \bmod^{*} 4\mathcal{O}_{K}$} ;
\item \textbf{il existe un idéal fractionnaire $I$ de $K$ tel que $I^2=\alpha \mathcal{O}_K$.}
\end{enumerate}
Il est clair que $R_{K}$ est un  sous-groupe de $V_2(K)$ contenant  $(K^\times)^2$ et donc
 $\bar{R}_{K}: =R_{K}/(K^\times)^2$ est un sous-groupe du $2$-Selmer groupe $S_2(K)$.

  Soit $L/K$ une extension
quadratique non ramifiée ;    d'après la  Proposition \ref{fr-th2} et le  Corollaire \ref{fr-cor1}, il existe un unique élément
$\bar{\alpha}\in \bar{R}_{K}$ tel que $\bar{\alpha}\not = \bar {1}$ et  $L=K(\sqrt \alpha)$.  On en déduit une bijection entre   $ \bar{R}_{K}\setminus \{\bar {1}\}$ et l'ensemble  des extensions quadratiques non ramifiées de $K$.  L'étude des classes de Steinitz des extensions quadratiques non ramifiées pourrait donc être considérée comme une nouvelle investigation du $2$-Selmer groupe d'un corps de nombres.

On note  $U_{K}=E_{K}\cap R_{K}$  (rappelons que $E_{K}$ est  le groupe des
unités de $K$). On a   $$U_{K}=\{\varepsilon \in E_K,  \exists \ x\in K^\times,  \varepsilon \equiv x^2 \bmod^* 4\mathcal{O}_K\}.$$
La Proposition \ref{fr-unitemodulocarre}   permet une description plus simple de  $U_K$ :
\begin{equation}\label{fr-U_Kextension2}
  U_{K}=\{\varepsilon \in E_K,  \exists \ a\in \mathcal{O}_K,  \varepsilon \equiv a^2 \bmod 4\mathcal{O}_K\}.
\end{equation}

On a  $E_{K}^2\subset U_{K}$. Soit  $\varepsilon \in E_K$, alors
\begin{equation}\label{fr-U_Kramifiee}
 K(\sqrt \varepsilon)/K \  \text {est quadratique et  non ramifiée} \Leftrightarrow \varepsilon \in U_K \setminus E_{K}^2
\end{equation}

\noindent \textbf{Remarque.} Dans les notations de \cite[Définition, p. 226]{fr-Dummitquebec},  l'affirmation sur  $E_{K, 4}$ (disant en gros que $E_{K, 4}$ - décrit à l'aide du \cite[Lemme 1, p. 226]{fr-Dummitquebec} - est notre $U_K \setminus E_{K}^2$)  qui suit cette définition n'est pas bien justifiée, même en utilisant la référence \cite[Proposition 4.8]{fr-Dummitlondon}, si on ne prouve pas notre assertion (\ref{fr-U_Kextension2}).
De plus, les questions posées dans Remarque 4 de \cite[p. 226]{fr-Dummitquebec} ont une réponse immédiate grâce à l'utilisation  de la congruence $\mod ^*$, laquelle est inexistante dans \cite{fr-Dummitquebec}. $\blacksquare$

\begin{Df}\label{fr-Dftypeclasse}
Nous dirons, d'après Hasse (voir \cite[\S22]{fr-hass}), que le corps $K$ est de type classe si
$U_K=E_{K}^2$, sinon  $K$ est dit de type unité.
\end{Df}

     D'après (\ref{fr-U_Kramifiee}) : pour que  le corps $K$ soit de type classe, il faut et il suffit que toute extension quadratique de $K$  de la forme $K(\sqrt{\varepsilon})$, où $\varepsilon$ est une unité de $K$, soit ramifiée en au moins un idéal premier  $\mathfrak{p}$. Noter qu'on a utilisé cette caractérisation de corps de type classe  comme définition dans \S1 pour simplifier la présentation des résultats. Noter aussi qu'un tel $\mathfrak{p}$ est nécessairement au dessus de $2$, car
 $\mathcal{B}=(1, \sqrt{\varepsilon} )$ est une base entière du $K$-espace vectoriel $K(\sqrt{\varepsilon})$ de discriminant  $d=4\mathcal{O}_{K}$ et
 $\delta_{ K(\sqrt{\varepsilon})/ K}$ divise $d$.
\smallskip

 \noindent \noindent \textbf{Remarque.}
Si $h^+_K$ est impair (donc $h^+_K=h_K$), alors $K$ est de type classe. Par exemple $K=\mathbb{Q}$, ou si $K/\mathbb{Q}$ est totalement imaginaire (donc $h^+_K=h_K$) et $h_K$ est impair ; pour construire des exemples dans le dernier cas, rappelons que si $K/\mathbb{Q}$ est galoisienne, alors elle est totalement imaginaire ou bien totalement réelle. $\blacksquare$ \smallskip

 Soit $K$ un corps de nombres ayant $r\geq 1$ plongements réels $\sigma_i$. Notons $s$ le morphisme de signature de $K$ : il est défini sur  $K^\times$ et à valeur dans    $ \{1, -1\}^r$,   et à $x$ fait correspondre $(signe(\sigma_1(x)), signe(\sigma_2(x)), ..., signe(\sigma_r(x)))$, où\break $signe(\sigma_i(y))=\frac{\sigma_i(y)}{|\sigma_i(y)|}$ lorsque $y\in  K^\times$.

\begin{Prop}\label{fr-totalement reel} Soit $K$ un corps de nombres totalement réel. Alors, $K$ est de type classe si et seulement si la restriction de la signature  $s$ de $K^\times$ à $E_K$  est surjective.
  \end{Prop}
\begin{proof}  Soit $\pi : E_K/E_{K}^2 \rightarrow (\mathcal{O}_K/4\mathcal{O}_K)^\times/((\mathcal{O}_K/4\mathcal{O}_K)^\times)^2$, $\overline{\varepsilon} \mapsto \overline{\varepsilon}$.
D'après \cite[Théorème 3.5,  $(c)\Leftrightarrow (a)$]{fr-muller}, on a : $\pi$ est injective ((c)) est équivalent à la restriction de   $s$  à $E_K$  est surjective ((a)).  Un calcul simple nous donne : $\pi$ est injective si  et  seulement si $E_{K}^2=\{\varepsilon,  \exists \ a \in \mathcal{O}_K,  \varepsilon \equiv a^2 \bmod 4\mathcal{O}_K\}$, ce qui est équivalent à $E_{K}^2=U_K$ par (\ref{fr-U_Kextension2}), équivalent aussi à $K$ est de type classe.
\end{proof}

\noindent \textbf{Remarque.} Cette proposition est le Théorème 3.1.1 dans \cite[p. 23]{fr-Derhem}. Nous avons constaté que sa démonstration était incomplète, ce qui nous a obligé de démontrer la Proposition \ref{fr-unitemodulocarre} pour la compléter. $\blacksquare$ \smallskip

Le but des deux propositions suivantes est de donner des exemples d'extensions quadratiques de $\mathbb{Q}$   de type classe et d'autres de type unité.
En fait, nous allons caractériser celles qui sont de type classe et celles qui sont de type unité.

  Soit $K/\mathbb{Q}$  une extension quadratique imaginaire. Si  $K = \mathbb{Q}(i)= \mathbb{Q}(\sqrt{-1})$ ou $K = \mathbb{Q}(j)= \mathbb{Q}(\sqrt{-3})$, où $j$ est une racine primitive 3ème de l'unité, alors $h_K=1$ (impair) implique que $K$ est de type classe. Si $K \not= \mathbb{Q}(\sqrt{-1})$ et $K \not= \mathbb{Q}(\sqrt{-3})$, on a : $K$ est de type classe si et seulement si $K(\sqrt{-1})/K$ est ramifiée, en effet :  $E_{K}=\{1, -1\}$ et donc $K(\sqrt{-1})$ est l'unique extension quadratique de $K$ de la forme  $K(\sqrt{\varepsilon})$, où $\varepsilon$ est une unité de $K$.

\begin{Prop}\label{fr-imaginaire}
Soit  $K=\mathbb{Q}(\sqrt{d})/\mathbb{Q}$ une  extension quadratique imaginaire, où $d$ est un élément de $\mathbb{Z}$   sans facteur carré. Alors :

(i)  Le corps $K$ est de type unité $\Longleftrightarrow$  $d\not=-1$ et   $d\equiv  3 \mod 4$.

(ii) Le corps $K$ est de type classe $\Longleftrightarrow$    $d=-1$ ou $d\equiv  1$ (resp. $2$) $\mod 4$.
\end{Prop}

\begin{proof}

Comme la négation de l'assertion  type classe est type unité, pour avoir l'équivalence dans (i) et (ii), il suffit de montrer dans  (i) et (ii) les implications $\Longleftarrow$. On sait que si $d=-1$, alors $K$ est de type classe ; dans la suite on suppose $d\not=-1$.

(i) On a $ (\sqrt{d})^2  = d \equiv -1 \mod 4\mathcal{O}_K$. D'après (\ref{fr-U_Kextension2}) et (\ref{fr-U_Kramifiee}),   $K(\sqrt{-1})/K$ est non  ramifiée, par suite $K$ est de type unité.

$(ii)$ 
On considère  l'extension biquadratique $K(\sqrt{-1})=K(i)$ qui est la composée de $K$ et $\mathbb{Q}(i)$.  Les seuls premiers qui peuvent se ramifier dans $K(i)/K$ sont ceux au dessus de $2$ (car $-1$ est une unité).

Supposons $d\equiv 1 \mod 4$.  Le premier $2$ est non ramifié dans $\mathbb{Q}(\sqrt{d})$ puisque $d\equiv 1 \mod 4$ ; comme il est ramifié dans $\mathbb{Q}(i)/\mathbb{Q}$,
son indice de ramification dans $K(i)/\mathbb{Q}$ est donc $2$ (sinon il serait totalement ramifié). Maintenant, $2$ n'est pas  ramifié dans $K$, par suite il est  ramifié dans $K(i)/K$. On conclut que $K$ est de type classe.

Supposons $d \equiv 2 \mod 4$. L'extension biquadratique $K(i)$ est une composée de $\mathbb{Q}(\sqrt{d})$ et $\mathbb{Q}(\sqrt{-1}) $, d'après \cite[Theorem 3, (iii)]{fr-Wil}, $\delta_{ K(i)/ \mathbb{Q}}=64d^2\mathbb{Z}$ (dans les notations de \cite[Theorem 3, (iii)]{fr-Wil}, $m=d$, $n=-1$, $l=PGCD(m,n)$\break $=1$, $m_1=(m/l)=d$, $m_2=(-1/l)=-1$,  $m\equiv2 \mod 4$, $-1 \equiv 3 \mod  4$).  La formule de transitivité du discriminant dans une tour de corps de nombres nous donne
   $\delta_{ K(i)/ \mathbb{Q}}=\delta_{ K/ \mathbb{Q}}^2 N_{K/\mathbb{Q}}(\delta_{ K(i)/ K})$, où $N_{K/\mathbb{Q}}$ est la norme. Sachant que $\delta_{ K/ \mathbb{Q}}=4d\mathbb{Z}$, on obtient $N_{K/\mathbb{Q}}(\delta_{ K(i)/ K})=4\mathbb{Z}$. Par suite $K(i)/ K$ est ramifiée et donc $K$ est de type classe.
\end{proof}

\begin{Prop}\label{fr-reel}
Soit $d$ un entier naturel   sans facteur carré et $K=\mathbb{Q}(\sqrt{d})$. Soit $\epsilon$ l'unité fondamental de $K$.

(i)  Le corps $K$ est de type unité si et seulement si  $N_{K/\mathbb{Q}}(\epsilon)=1$.

(ii)  Le corps $K$ est de type classe si et seulement si  $N_{K/\mathbb{Q}}(\epsilon)=-1$.

\end{Prop}

\begin{proof}
 Notons $s_{E_K}$ la restriction à $E_K$  du morphisme de signature $s$ de  $K^\times $ dans $\{1, -1\}^2$.
 Il est immédiat que si $N_{K/\mathbb{Q}}(\epsilon)=-1$, alors $Im(s_{E_K})=\{1, -1\}^2$,  et si $N_{K/\mathbb{Q}}(\epsilon)=1$, alors $Im(s_{E_K})=\{(1,1), (-1,-1)\}$ ; comme $K$ est totalement réel, la Proposition \ref{fr-totalement reel}  nous permet de conclure.
 \end{proof}

\begin{Prop}\label{fr-unité}
Soit $L=K(\sqrt{\alpha})$ une extension quadratique de $K$ non ramifiée. Alors,
  il existe $\varepsilon$ une unité de $K$ tel que $L=K(\sqrt{\varepsilon})$ si et seulement si $cl_{K}(\mathcal{O}_{L})=1$.

 \end{Prop}
 \begin{proof}
 Supposons $L=K(\sqrt{\varepsilon})$, où $\varepsilon$ une unité de $K$ ;  d'après le théorème d'Artin $cl_{K}(\mathcal{O}_{L})=\left[(\dfrac{\mathcal{O}_{K}}{4\varepsilon})^{1/2}\right]_{K}=1.$
Réciproquement,   supposons $cl_{K}(\mathcal{O}_{L})=1$ ;  d'après le Corollaire \ref{fr-cor1}, il existe un idéal fractionnaire $I$ de $K$ tel que $I^{2}=\alpha\mathcal{O}_{K}$ et $[I]_K=cl_{K}(\mathcal{O}_{L})$, d'où $[I]_K=1$ et donc il existe  $\beta \in K^{*}$ et une unité $\varepsilon$ de $K$ vérifiant $\beta^{2}\varepsilon=\alpha $. On en déduit $K(\sqrt{\alpha})=K(\sqrt{\varepsilon})$.
 \end{proof}

 \noindent \textbf{Remarque.} Le résultat de la proposition précédente est le même que celui de \cite[Proposition 4]{fr-Dummitquebec} (dans cette dernière on donne une base assez explicite du $\mathcal{O}_{K}$-module $\mathcal{O}_{L}$), mais ici la preuve est différente. $\blacksquare$
\smallskip

 Immédiatement on a la caractérisation suivante d'un corps de type classe à l'aide de la notion de la classe de Steinitz.
 \begin{Cor}\label{fr-cor2}
 Soit $K$ un corps de nombres. Alors, pour que  $K$ soit    de type classe, il faut et il suffit que,   pour toute  extension quadratique non ramifiée $L/K$, on ait $cl_{K}(\mathcal{O}_{L})\not = 1$ ($\Leftrightarrow$ $\mathcal{O}_L$ n'est pas un $\mathcal{O}_K$-module libre $\Leftrightarrow$ $\vert cl_{K}(\mathcal{O}_{L})\vert=2)$).
 \end{Cor}

 \noindent \textbf{Remarque.}  Dans le corollaire précédent, la condition $cl_{K}(\mathcal{O}_{L})\not = 1$ est équivalente d'une part à $\mathcal{O}_L$ n'est pas un $\mathcal{O}_K$-module libre par définition de la classes de Steinitz, et d'autre part à $cl_{K}(\mathcal{O}_{L})$ est d'ordre $2$ dans $Cl_{K}$ puisque $cl_{K}(\mathcal{O}_{L})^2=1$. $\blacksquare$

 \begin{proof}[\textbf{Démonstration du Théorème \ref{fr-t1}}]
 (1) D'après la première remarque qui suit la  Proposition \ref{fr-Artin},  $R'_{nr}(K,\mathbb{Z}/2\mathbb{Z})\subset Cl_{K,2}$.
  Comme  $R'_{nr}(K,\mathbb{Z}/2\mathbb{Z})$ est non vide et fini,    $R'_{nr}(K,\mathbb{Z}/2\mathbb{Z})$ est un sous-groupe de $Cl_{K,2}$ si et seulement s'il    est stable pour la loi de $Cl_{K}$. Soient $c_1$ et $c_2$ deux éléments de $R'_{nr}(K,\mathbb{Z}/2\mathbb{Z})$. Si  $c_1=c_2$, alors  $c_1c_2=1 \in R'_{nr}(K,\mathbb{Z}/2\mathbb{Z}) $. Si  $c_1\not=c_2$, Soient  $L_i/K$, $1\leq i \leq 2$, deux extensions quadratiques non ramifiées telles que $cl_{K}(\mathcal{O}_{L_1})=c_1$ et $cl_{K}(\mathcal{O}_{L_2})=c_2$. La composée $L_1L_2/K$ est une extension biquadratique non ramifiée ; elle contient une troisième extension quadratique non ramifiée $L_3/K$. Le Corollaire \ref{fr-cor4} nous donne $cl_{K}(\mathcal{O}_{L_3})=c_1c_2$. Donc
  $c_1c_2\in R'_{nr}(K,\mathbb{Z}/2\mathbb{Z})$.

(2) On sait déjà que si $h^+_K$ est impair, alors $R_{nr}(K,\mathbb{Z}/2\mathbb{Z})=\emptyset$ et\break  $R'_{nr}(K,\mathbb{Z}/2\mathbb{Z})=\{1\}$.

 (3)  Maintenant  supposons que $h^+_K$ est pair et $K$ est de type classe.  Rappelons que $\vert Cl_{K,2} \vert =2^{r_{2}(K)}$. Montrons que
  $R'_{nr}(K,\mathbb{Z}/2\mathbb{Z})=Cl_{K,2}$. Pour ce faire,
puisque $R'_{nr}(K,\mathbb{Z}/2\mathbb{Z})\subset Cl_{K,2}$, il   suffit de montrer que  $\vert R'_{nr}(K,\mathbb{Z}/2\mathbb{Z}) \vert$\break $=2^{r_{2}(K)}$.

 Par la  Proposition  \ref{fr-2rang} (2) et la théorie de Galois, il existe $n=:2^{r^+_{2}(K)}-1$ extensions quadratiques non ramifiées sur $K$ ; on les note $(K_{i})_{1\leq i\leq n}$.
Soient $i$ et $j$ deux éléments distincts  de  $ \lbrace 1,...,n\rbrace$. La composée de $K_{i}/K$ et $K_{j}/K$ est une extension biquadratique non ramifiée de $K$ ; il existe donc  $m\in \lbrace 1,...,n\rbrace \setminus \lbrace i,j\rbrace$ tel que sa troisième sous-extension quadratique est $K_{m}/K$.
   D'après le Corollaire \ref{fr-cor4},  $cl_{K}(\mathcal{O}_{K_{i}})cl_{K}(\mathcal{O}_{K_{j}})=cl_{K}(\mathcal{O}_{K_{m}})$.
 Grâce au Corollaire \ref{fr-cor2} on a $\vert cl_{K}(\mathcal{O}_{K_{m}}) \vert = 2$ , d'où   $cl_{K}(\mathcal{O}_{K_{i}})\not = cl_{K}(\mathcal{O}_{K_{j}})$. On en déduit que le cardinal de  $  R_{nr}(K,\mathbb{Z}/2\mathbb{Z})$ est égal à $ 2^{r^+_{2}(K)}-1$ et
$ \vert R'_{nr}(K,\mathbb{Z}/2\mathbb{Z}) \vert =2^{r^+_{2}(K)}$ ; la dernière égalité  implique $r^+_{2}(K)\leq r_{2}(K)$  puisque $R'_{nr}(K,\mathbb{Z}/2\mathbb{Z})\subset Cl_{K,2}$.
Le groupe $ Cl_{K}$ étant isomorphe à un groupe quotient de  $Cl^{+}_{K}$, il vient de  la Proposition \ref{fr-2rang} (5) que
  $r_{2}(K)\leq r^+_{2}(K)$. Par conséquent $r_{2}(K)=r^+_{2}(K)$, d'où  $\vert R'_{nr}(K,\mathbb{Z}/2\mathbb{Z})\vert=2^{r_{2}(K)}$ et  $R_{nr}(K,\mathbb{Z}/2\mathbb{Z}) =Cl_{K,2}\setminus\{1\}$.

  (4) Supposons $K$ de type unité. Par définition de $m$, il y a $2^m-1$ extensions quadratiques non ramifiées de $K$, $M_i=K(\sqrt{\varepsilon_i})$,  de classes de Steinitz  triviales. D'où
 $ (2^{r^+_{2}(K)}-1)-(2^m-1)=2^{r^+_{2}(K)}-2^m:=a$ extensions quadratiques non ramifiées de $K$ de classes de Steinitz  non triviales. Si $L=K(\sqrt{\alpha})$ est l'une d'entre elles, il y a   $2^m$ extensions quadratiques  ayant la même classe que $L$ : $L$ et les $2^m-1$   uniques sous-extensions quadratiques $K(\sqrt{\alpha \varepsilon_i})$ des extensions biquadratiques $LM_i$ (la composée de $L$ et $M_i$). Par suite Il y a $\frac{a}{2^m}=2^{r^+_{2}(K)-m}-1$ classes de Steinitz non triviales. On ajoute la classe triviale pour terminer la preuve.
  \end{proof}

\noindent \textbf {Remarques.}

\textbf{1.} On vient de voir que si $K$ est de type classe et $h^+_K$ est pair, alors $r_{2}(K)=r^+_{2}(K)$. Donc, dans ce cas, toute extension quadratique de $K$ non ramifiée aux places finies est aussi non ramifiée aux places infinies.

\textbf{2.} Sachant que si  $K$ est  de type unité, alors $h^+_K$ est pair,  on déduit de la remarque \textbf{1} : si $r_{2}(K)\not=r^+_{2}(K)$, alors $K$ est de type unité. $\blacksquare$

\medskip

\noindent \textbf{Exemples assez  explicites du Théorème \ref{fr-t1} dans le cas des extensions quadratiques  $K/\mathbb{Q}$.}

Nous utilisons les notations de la Proposition \ref{fr-K/kquadratique} et du Corollaire \ref{fr-corK/kquadratique}.

Rappelons la définition suivante : si $L/\mathbb{Q}$ est une extension galoisienne abélienne, le corps de genre de $L$ au sens restreint, noté $L^+_{gen}$, est la plus grande  extension  abélienne de $\mathbb{ Q}$, contenant $L$ et non ramifiée sur $L$ (aux places finies).
Rappelons aussi un résultat bien connu de la théorie du corps de classes : si $L'/L$ est une extension quadratique de corps de nombres et $L''/L'$ est une extension abélienne non ramifiée de degré premier avec $h^+_L$, alors $L''/L$ est galoisienne.

La proposition suivante nous dit, avec précision,  où se trouvent  les extensions quadratiques  de $K$  non ramifiées sur $K$.
\begin{Prop}\label{fr-genre}
  Soit  $K/\mathbb{Q}$ une extension quadratique. Soit  $M$ la composée de toutes les extensions quadratiques  de $K$  non ramifiées sur $K$. Alors $M=K^+_{gen}$.
\end{Prop}

\begin{proof}
 Le corps $K^+_{gen}$ est bien connu (voir par exemple \cite[(1.4), p. 4]{fr-Ichida}) ; en particulier on sait que $[K^+_{gen}:K]=2^{u-1}$ (\cite[(1.5), p. 5]{fr-Ichida}). Soit $L/K$ une extension quadratique non ramifiée. Puisque
$ [K:\mathbb{Q}]=2$, $PGCD(h^+_\mathbb{Q},[L:K])=1 $ et $L/K$ est abélienne non ramifiée, on a $L/\mathbb{Q}$ est galoisienne par le rappel précédent ; elle est abélienne car son degré est $4$. On en déduit que $M\subset K^+_{gen}$. Or $[M:K]=2^{r^+_2(K)}$ et $r^+_2(K)=u-1$ (Corollaire \ref{fr-corK/kquadratique} (i)). Par conséquent
$M=K^+_{gen}$.
  \end{proof}

  Rappelons la description explicite de $K^+_{gen}$, où  $K=\mathbb{Q}(\sqrt{d})$ et $d$ est un élément de $\mathbb{Z }$  sans facteur carré (voir \cite{fr-Ichida}).
Soit $p_i$, $1\leq i \leq u$,  un nombre premier divisant $\delta_{K/\mathbb{Q} }$.
 Si $d \equiv 1 \mod 4$ : on pose $p_i^*= (\dfrac{-1}{p})p_i=    (-1)^{\frac{p-1}{2}}p_i$.
Si $d \equiv 2, 3  \mod 4$ : on  note  $p_1=2$ ; on  pose $p_1^*=-4$ (resp. 8, -8) si $d \equiv 3 \mod 4$ ( resp. $d \equiv 2 \mod 8$, $d \equiv -2 \mod 8$) ;   pour $2\leq i \leq u$, on pose $p_i^*=    (-1)^{\frac{p-1}{2}}p_i$ .
On a
 $K^+_{gen}=K(\sqrt{p_i^*}, 2\leq i \leq u)=\mathbb{Q}(\sqrt{d}, \sqrt{p_i^*}, 2\leq i \leq u).$

 Notons $\mathfrak{p}_i$, $1\leq i \leq u$, l'idéal premier de $\mathcal{O}_K$ au dessus de $p_i$.

\begin{Cor}\label{fr-corgenre} Sous les notations précédentes,
    $R'_{nr}(K,\mathbb{Z}/2\mathbb{Z})$ est un sous-groupe de         $ Am_0(Cl_K)$  engendré par les classes $[\mathfrak{p}_i]_{K}$, $2\leq i \leq u$. De plus,  $R'_{nr}(K,\mathbb{Z}/2\mathbb{Z})= Am_0(K/\mathbb{Q})$  si
   $d\equiv 1 $ ou $2  \mod 4$.
\end{Cor}
\begin{proof} La première partie découle de la Proposition \ref{fr-genre} et la description de $K^+_{gen}$. Reste ce qui est en plus.

 Si $d\equiv 1 $ ou $2  \mod 4$, alors $\prod_{i=1}^{u}\mathfrak{p}_i=\sqrt{d}\mathcal{O}_K$, d'où $[\mathfrak{p}_1]_{K}\in R'_{nr}(K,\mathbb{Z}/2\mathbb{Z})$ et donc
  $R'_{nr}(K,\mathbb{Z}/2\mathbb{Z})= Am_0(K/\mathbb{Q})$.
\end{proof}

\textbf{1. Cas $h^+_K$ impair.}

Dans ce cas, $K$ est de type classe,   $R_{nr}(K,\mathbb{Z}/2\mathbb{Z})=\emptyset$ et  $R'_{nr}(K,\mathbb{Z}/2\mathbb{Z})=\{1\}.$

\begin{Prop}
  Le nombre $h^+_K$ est impair si et seulement si $K$ a l'une des formes suivantes :  $\mathbb{Q}(\sqrt{-1})$,  $\mathbb{Q}(\sqrt{\pm 2})$, $\mathbb{Q}(\sqrt{p})$ avec $p  \equiv 1 \mod{4}$, ou $\mathbb{Q}(\sqrt{-p})$ avec $p \equiv 3 \mod 4$, où $p$ est un nombre premier.
\end{Prop}

\begin{proof} Le cas $u=0$ correspond à $\mathbb{Q}(\sqrt{-1})$. Le cas $u\not =0$ est immédiat en utilisant l'assertion suivante :
  $h^+_K$ est impair si et seulement si $r^+_2(K)=0$ ;  ce qui est équivalent à $u=1$,  d'après le Corollaire \ref{fr-corK/kquadratique} (i).
\end{proof}

\textbf{2. Cas $h^+_K$ pair.}

\textbf{2.1. $K$ de type classe.}

 On a $R'_{nr}(K,\mathbb{Z}/2\mathbb{Z})=Cl_{K,2}.$ De plus
d'après la  Proposition \ref{fr-imaginaire} (ii),    la Proposition \ref{fr-reel} (ii) et le Corollaire \ref{fr-descclasseK2}, $R'_{nr}(K,\mathbb{Z}/2\mathbb{Z})=Am_0(K/\mathbb{Q}).$

\textbf{2.2. $K$ de type unité.}

 Rappelons que $| R_{nr}(K,\mathbb{Z}/2\mathbb{Z})| =2^{r^+_2(K)-m}$, où $m\geq 1$.

  \textbf{Le cas $r^+_2(K)=r_2(K)$ :}
   $K=\mathbb{Q}(\sqrt{d})$ imaginaire, avec    $d\not=-1$ et   $d  \equiv 3 \mod 4$, ou $K=\mathbb{Q}(\sqrt{d})$ est réel avec  $N_{K/\mathbb{Q}}(\epsilon)=1$ et $-1$ est  une norme dans  $K/\mathbb{Q}$.

    L'inégalité  $m\geq 1$ implique $R_{nr}(K,\mathbb{Z}/2\mathbb{Z})$ est un sous-groupe   de $Cl_{K,2}$ différent de $Cl_{K,2}$.

 Dans le cas imaginaire  $| R_{nr}(K,\mathbb{Z}/2\mathbb{Z})| =2^{r_2(K)-1}$.

  Dans le cas réel, en supposant $d$ impair, on a $d\equiv 1 \mod 4$ d'après le paragraphe juste après le Corollaire \ref{fr-corK/kquadratique} (Théorème de Hasse) ;   la non ramification de $2$ dans $K/\mathbb{Q}$ implique
$K(\sqrt{-1})/K$ est ramifiée (par exemple, car $K/\mathbb{Q}$ et $\mathbb{Q}(\sqrt{-1})/\mathbb{Q}$ sont arithmétiquement disjointes), d'où $m=1$ et  $| R_{nr}(K,\mathbb{Z}/2\mathbb{Z})| =2^{r_2(K)-1}$ ; notons que si $d\equiv 1 \mod 8$, d'après \cite[Théorème 1.1 (i)]{fr-LA}, $\epsilon$ est totalement positif implique  $\epsilon \equiv 1 \mod 4\mathcal{O}_K$ et donc $K(\sqrt{\epsilon})/K$ est non ramifiée (et par suite $K(\sqrt{-\epsilon})/K$ est ramifiée).

\textbf{Le cas $r^+_2(K)=r_2(K)+1$ :}
  $K$  réel, avec  $N_{K/\mathbb{Q}}(\epsilon)=1$ et $-1$ n'est pas une norme dans  $K/\mathbb{Q}$ (voir le Corollaire\ref{fr-r^+_2(K)= r_2(K)+1}).

   Dans ce cas $Cl_{K,2}=Am_0(K/\mathbb{Q})$ (d'ordre
  $2^{u-2}$) et $|R_{nr}(K,\mathbb{Z}/2\mathbb{Z})|=2^{r_2(K)+1 -m}$. La proposition suivante nous donne une infinité d'exemples d'un tel $K$ avec $R_{nr}(K,\mathbb{Z}/2\mathbb{Z})=Cl_{K,2}$ (équivalent à $m=1$).

 \begin{Prop}\label{fr-exempletheoreme1.1unite} Soit $K=\mathbb{Q}(\sqrt{d})$ un corps quadratique réel, où $d$ est un entier naturel  sans facteur carré. Supposons que les diviseurs $p$ de $d$ soient tous $\equiv 3 \mod 4$ et leur nombre soit pair non nul. Alors $K$ est de type unité
 et   $R_{nr}(K,\mathbb{Z}/2\mathbb{Z})=Cl_{K,2}=Am_0(K/\mathbb{Q})$.
 \end{Prop}

\begin{proof} Le discriminant de $K/\mathbb{Q}$ est divisible par un premier $p \equiv 3 \mod 4$, donc $N_{K/\mathbb{Q}}(\epsilon)=1$ (voir par exemple \cite[(1.7), p. 179]{fr-FT}) ; par suite $K$ est de type unité par la Proposition \ref{fr-reel} (ii). On a $(\dfrac{-1}{p})=-1$,  par conséquent $-1$ n'est pas une norme dans  $K/\mathbb{Q}$. On en déduit : $r^+_2(K)=r_2(K)+1$.  Les deux conditions sur les diviseurs premiers de $d$  nous donnent $d \equiv 1 \mod 4$, d'où $\delta_{K/\mathbb{Q}}=d\mathbb{Z}$. La non ramification de $2$ dans $K/\mathbb{Q}$ implique
$K(\sqrt{-1})/K$ est ramifiée.  La remarque 2 qui suit le Théorème \ref{fr-t1} nous permet d'obtenir $m=1$. On conclut que $|R_{nr}(K,\mathbb{Z}/2\mathbb{Z})|=2^{r_2(K)}$ et donc $R_{nr}(K,\mathbb{Z}/2\mathbb{Z})=Cl_{K,2}$
 \end{proof}

\section{Démonstration du Théorème 1.2}

  Dans cette section on démontre le Théorème \ref{fr-t2}.
  \begin{proof}[\textbf{Démonstration du Théorème \ref{fr-t2}}] Soit $\Gamma$ un groupe fini.

   Si l'ordre de $\Gamma$ est impair ou un $2$-sous-groupe de Sylow de $\Gamma$ n'est pas cyclique, alors $ R'_{nr}(K,\Gamma)=\{1\}$ d'après le Corollaire \ref{fr-cor3} (i) ; dans ce cas on a donc le Théorème \ref{fr-t2}.

    Dorénavant nous supposons que l'ordre de $\Gamma$ est pair et  un $2$-sous-groupe de Sylow de $\Gamma$ est cyclique. Dans la suite, et en vue d'avoir le résultat le plus général possible,  nous ne ferons intervenir la normalité (et donc l'unicité) du $2$-sous-groupe de Sylow qu'en cas de besoin.

     Tout d'abord, le  Corollaire \ref{fr-cor3} $(ii)$ nous dit que $ R'_{nr}(K,\Gamma) \subset R'_{nr}(K,\mathbb{Z}/2\mathbb{Z})$. Ensuite,  comme $R'_{nr}(K,\Gamma)$ est non vide et fini,  pour prouver le  Théorème \ref{fr-t2} il suffit de  montrer que $R'_{nr}(K,\Gamma)$ est stable pour la loi de $Cl_{K}$.

 Soient $c_1$ et $c_2$ des éléments de $R'_{nr}(K,\Gamma)$. Si $c_1=1$, ou $c_2=1$ ou  $c_1=c_2$, il est immédiat que  $c_1c_2\in R'_{nr}(K,\Gamma) $. Désormais on suppose  $c_1\not=1$,  $c_2\not= 1$ et  $c_1\not =c_2$.

 Soient $N_{1}/K$ et $N_{2}/K$ deux $\Gamma$-extensions non ramifiées telles que\break
   $cl_{K}(\mathcal{O}_{N_{1}})=c_1$ et  $cl_{K}(\mathcal{O}_{N_{2}})=c_2$.
   Grâce à la  Proposition  \ref{fr-prop1} $(ii)$,   l'extension $N_{i}/K$, $1\leq i \leq 2$, admet une seule sous-extension quadratique sur $K$ que l'on note  $L_{i}/K$.  Par   le Corollaire \ref{fr-cor3} $(ii)$, on a $c_i = cl_{K}(\mathcal{O}_{L_{i}})$, d'où
  $ c_1c_2 = cl_{K}(\mathcal{O}_{L_{1}}) cl_{K}(\mathcal{O}_{L_{2}}).$

  Considérons le diagramme suivant :

 \begin{center}

  \begin{tikzpicture}[node distance = 1cm, auto]
  \node (Q) {$K$};
  \node (E) [above of=Q] {$N_{1}\cap N_{2}$};
  \node (E1) [above of=E, left of=E] {$N_{1}$};
  \node (F) [above of=E, right of=E] {$N_{2}$};
  \node (K) [above of=E, node distance =2cm] {$N_{1}N_{2}$};
   \draw (Q) --(E);
  \draw (E) -- (F);
  \draw (E) --  (E1);
  \draw (F) --  (K);
  \draw (E1) -- (K);
  \end{tikzpicture}

  \end{center}

  Posons $Gal(N_{1}N_{2}/K)=G, Gal(N_{1}N_{2}/N_{1}\cap N_{2})=G', Gal(N_{i}/K)=G_i\simeq \Gamma, Gal(N_{i}/N_{1}\cap N_{2})=H_i, H'_{1}=Gal(N_1N_2/N_2), H'_{2}=Gal(N_1N_2/N_1).$

Raisonnons par disjonction des cas selon la parité de $[N_{1}\cap N_{2}: K]$.

  \textbf{Premier cas :}   $[N_{1}\cap N_{2}: K]$ est pair.

 Par la théorie de Galois, $Gal(N_{1}\cap N_{2}/K)\simeq G_1/H_1$. Comme $syl_2(G_1/H_1)= syl_2(G_1)H_1/H_1$ (Proposition \ref{fr-prp3} (1.2)) et $syl_2(G_1)$  est cyclique, on a\break $ syl_2(Gal(N_{1}\cap N_{2}/K))$ est cyclique. D'après la  Proposition  \ref{fr-prop1} $(ii)$, $N_{1}\cap N_{2}/K$ admet une seule sous-extension quadratique sur $K$, d'où $L_{1}=L_{2}$ et donc $c_{1}c_{2}=1$ appartient à $R'_{nr}(K,\Gamma)$. \smallskip

\noindent {\bf Remarque.} Dans ce cas nous n'avons pas besoin de l'hypothèse $syl_2(\Gamma)$ est normal. $\blacksquare$
\smallskip

   \textbf{Deuxième cas :}   $[N_{1}\cap N_{2}: K]$ est impair.

   Soit $L_{3}/K$ la $3$ème  sous-extension quadratique de l'extension biquadratique $L_{1}L_{2}/K$. D'après le Corollaire \ref{fr-cor4} on a $c_{1}c_{2}=cl_{K}(\mathcal{O}_{L_{3}}).$
    Montrons que $L_{3}/K$ est plongeable dans une $\Gamma$-extension non ramifiée sur $K$, ce qui nous donne $c_{1}c_{2} \in R'_{nr}(K,\Gamma)$ et termine la démonstration du Théorème \ref{fr-t2}.

Posons  $\vert\Gamma\vert=2^{n}k$, où $n$ est un entier naturel non nul et $k$ un entier naturel impair.

Commençons par la remarque suivante : si $M/K$ est une sous-extension galoisienne de $N_{1}N_{2}/K$ contenant $N_{1}\cap N_{2}$ et $H=Gal(M/N_{1}\cap N_{2})$,  alors
$syl_{2}(Gal(M/K))=syl_{2}(H)$. En effet, $syl_{2}(H)$ est un sous-groupe de\break  $Gal(M/K)$ et  $[N_{1}\cap N_{2}: K]$ est impair implique $\vert syl_{2}(H) \vert=$
 $\vert syl_{2}(Gal(M/K))\vert$.

En appliquant la remarque précédente à $Gal(N_{1}N_{2}/K)=G$ (resp. \break $Gal(N_{i}/K)=G_i$) et $Gal(N_{1}N_{2}/N_{1}\cap N_{2})=G'$ (resp. $Gal(N_{i}/N_{1}\cap N_{2})=H_i$) on obtient $syl_{2}(G)=syl_{2}(G')$ (resp. $syl_{2}(G_i)=syl_{2}(H_i))$.
Puisque $N_{1}/N_{1}\cap N_{2}$ et $N_{2}/N_{1}\cap N_{2}$ sont linéairement disjointes, on a   $G'\simeq H_1 \times H_2$, d'où $syl_{2}(G')\simeq syl_{2}(H_1) \times syl_{2}(H_2)$ par la Proposition \ref{fr-prp3} (2), et donc
$syl_{2}(G)\simeq syl_{2}(\Gamma) \times syl_{2}(\Gamma)$ (car  $syl_{2}(G_i)\simeq syl_{2}(\Gamma)$).

Voici une deuxième remarque :  par la théorie de Galois, $H'_{1}\simeq H_1$ et $H'_{2}\simeq H_2$ (isomorphisme de restriction) et par suite  $syl_{2}(H'_{i})\simeq syl_{2}(H_{i})  \simeq syl_{2}(\Gamma)$.

Notons que dans les deux remarques précédentes on n'avait pas besoin de l'hypothèse  $syl_{2}(\Gamma)$ cyclique.

 Comme   $syl_{2}(\Gamma)$ est  cyclique, alors $syl_{2}(G) \simeq (\mathbb{Z}/2^{n}\mathbb{Z})^{2}$ ; par conséquent
         $r_2(syl_{2}(G))=2$.

  Choisissons un générateur $\sigma$ (resp. $\tau$) de $syl_{2}(H'_{2})$ (resp. $syl_{2}(H'_{1})$). Immédiatement $syl_{2}(G)=\langle\sigma,\tau\rangle $
  (car $syl_{2}(G)\simeq syl_{2}(H'_{2})\times syl_{2}(H'_{1})$).

 Dans toute la suite on suppose que $\Gamma$ possède un  $2$-sous-groupe de Sylow normal ; $syl_{2}(\Gamma)$ est donc unique. Notons $I$ l'unique sous-groupe normal de $\Gamma$ d'ordre $k$ (voir Proposition \ref{fr-Scott}). D'après la Proposition  \ref{fr-prp3} (4), le groupe $\Gamma$ est un produit direct interne de ses sous-groupes $syl_{2}(\Gamma)$ et $I$, car ses derniers sont normaux et d'ordres premiers entre eux (on peut aussi utiliser la remarque qui suit la Proposition \ref{fr-Scott}). Donc
 $\Gamma \simeq  syl_{2}(\Gamma) \times I.$

   La théorie de Galois nous dit que $G$ est isomorphe à un sous-groupe de $$\Gamma \times \Gamma \simeq syl_{2}(\Gamma) \times syl_{2}(\Gamma) \times I \times I .$$ Les ordres des  groupes
    $syl_{2}(\Gamma) \times syl_{2}(\Gamma)$ et  $ I \times I$ étant premiers entre eux, $G$ est isomorphe au produit d'un sous-groupe $R$ du premier et d'un sous-groupe $S$ du second par la Proposition \ref{fr-prp3} (3). Puisque
   $syl_{2}(G)\simeq syl_{2}(\Gamma)\times syl_{2}(\Gamma)$, on a $R\simeq syl_{2}(\Gamma)\times syl_{2}(\Gamma)$, d'où
 $G\simeq syl_{2}(G) \times S.$ Dans la démonstration de (3) de la Proposition \ref{fr-2rang}, on a vu que  $\langle(syl_{2}(G) \times S)^2\rangle =\langle syl_{2}(G)^2\rangle  \times \langle S^2\rangle .$
Par conséquent  $G/\langle G^2\rangle \simeq (syl_{2}(G)/\langle syl_{2}(G)^2\rangle ) \times (S/\langle S^2\rangle ).$
 L'ordre de $S$ est impair implique $S^2=S$ (utiliser l'identité de Bezout), d'où $ r_2(G)=r_2(syl_{2}(G))$ ; par conséquent $ r_2(G)=2$.
  De la Proposition \ref {fr-2rang} (2) découle que $G$ admet exactement trois sous-groupes d'indice $2$.   On peut les définir  par $\Gamma_1:=\langle\sigma,\tau^{2}\rangle  S$, $\Gamma_2:=\langle\tau, \sigma^{2}\rangle  S$, $\Gamma_3:=\langle\sigma\tau,\tau^{2}\rangle S,$
  où  l'on a identifié $S $ avec un sous-groupe de $G$.

   Comme $N{_1}^{\langle\sigma\rangle }=N{_1}$ et  $N{_2}^{\langle\tau\rangle }=N{_2}$,  on peut choisir
 $Gal(N_{1}N_{2}/L_{1})=\Gamma_1, Gal(N_{1}N_{2}/L_{2})=\Gamma_2, Gal(N_{1}N_{2}/L_{3})=\Gamma_3$.

 Puisque $ PGCD(\vert syl_{2}(G)\vert, \vert S \vert)=1$, d'après la Proposition \ref{fr-prp3} (3),  le sous-groupe $H'_1$ de $G$ est un produit direct d'un sous-groupe de $syl_{2}(G)$ et d'un sous-groupe $S_1$ de $S$  ; ces deux sous-groupes sont distingués respectivement dans  $syl_{2}(G)$ et  $S$  car $H'_1$ est distingué dans $G$. Comme  $syl_{2}(H'_1)=\langle\tau\rangle $, on a :
$H'_1=\langle\tau\rangle  S_1.$

 Considérons le sous-groupe de G suivant : $\Gamma_4 : =\langle\sigma\tau\rangle  S_{1}.$
   Dans la suite on va utiliser la Proposition \ref{fr-prp3} sans la mentionner. Le groupe $\Gamma_4$ est un produit direct de $\langle\sigma\tau\rangle $ et  $S_{1}$, par conséquent il est normal dans $G$ car $\langle\sigma\tau\rangle $ et  $S_{1}$ sont respectivement normaux dans $syl_{2}(G)$ et $S$.  Notons $N_3/K$ la sous-extension de $N_{1}N_{2}/K$ fixe par $\Gamma_4$.
Par la théorie de Galois $Gal(N_{3}/K)\simeq G/\Gamma_4.$ On a
   $$Gal(N_{3}/K) \simeq (\langle\sigma,\tau\rangle \times S/\langle\sigma\tau\rangle \times S_{1})\simeq ((\langle\sigma,\tau\rangle /\langle\sigma\tau\rangle ) \times S/ S_{1}).$$ Donc
    $Gal(N_{3}/K)\simeq  \mathbb{Z}/2^{n}\mathbb{Z} \times S/S_{1}.$
 Rappelons que $\Gamma \simeq Gal(N_2/K) \simeq G/H'_1$. On a $$\Gamma \simeq G/H'_1\simeq (\langle\sigma,\tau\rangle \times S)/(\langle\tau\rangle \times S_{1})\simeq \mathbb{Z}/2^{n}\mathbb{Z} \times S/ S_{1}.$$
  On en déduit que $Gal(N_{3}/K)\simeq \Gamma$.

Comme $\langle\sigma\tau\rangle  S_{1}$ est un sous-groupe de $\Gamma_3$, on a  $L_{3}\subset N_{3}$.
L'extension $L_{3}/K$ est donc  plongeable dans une $\Gamma$-extension $N_{3}/K$ ; $N_{3}/K$ est  non ramifiée car c'est une sous-extension de la composée de deux extensions non ramifiées. Ceci termine la démonstration du Théorème \ref{fr-t2}.
\end{proof}

\noindent \textbf{Remarque.}
 Comme pour $H'_1$, le sous-groupe $H'_2$   peut s'écrire comme un produit direct  $H'_2=\langle\sigma\rangle  S_2$ ($syl_{2}(H'_{2})=\langle\sigma\rangle $), où $S_2$ est un sous-groupe normal de $S$.
En considérant le sous-groupe produit de G : $\Gamma_5 : =\langle\sigma\tau\rangle  S_{2}$, et  la sous-extension $N_{4}/K$ de $N_{1}N_{2}/K$ fixe par $\Gamma_5$, on montre comme ci-dessus que   $L_{3}/K$ est   plongeable dans la  $\Gamma$-extension non ramifiée $N_{4}/K$. $\blacksquare$
\smallskip

Le but de la proposition suivante est de donner un exemple du Théorème \ref{fr-t2} lorsque $\Gamma\simeq \mathbb{Z}/4\mathbb{Z}$.

Notons  $H^{+}_{K, 2}$ le $2$-corps de classes de Hilbert de $K$ au sens restreint (rappelons que c'est la $2$-extension abélienne non ramifiée (aux places finies)  maximale de $K$ et que $Gal(H^{+}_{K, 2}/K)\simeq syl_2(Cl^+_K)$). Si $L/K$ est une $\mathbb{Z}/4\mathbb{Z}$-extension non ramifiée, alors $L\subset H^{+}_{K, 2}$. Dans la suite on note $H_0:=Gal(H^{+}_{K, 2}/K)$.

\begin{Prop}\label{fr-exempletheor1.2}
	Soit $K$ un corps de nombres de type classe dont le $2$-groupe de classes au sens restreint est de type $(2, 4)$ (i.e., $syl_2(Cl^+_K)\simeq \mathbb{Z}/2\mathbb{Z} \times \mathbb{Z}/4\mathbb{Z}$). Alors $\left|R'_{nr}(K, \mathbb{Z}/4\mathbb{Z})\right|=2$.
\end{Prop}

\begin{proof}

Il existe deux éléments $\tau$ et $\sigma$ de $H_0$  d'ordres respectifs $2$ et $4$ tels que
 $H_0=\langle\tau, \sigma\rangle =\{\tau^i\sigma^j, 0\leq i \leq 1, 0\leq j \leq 3 \}.$
 Le groupe $H_0$ admet 3 sous-groupes
  d'ordre $4$ d'indice 2 :
 $H_{1}=\langle\sigma\rangle $, $H_{2}=\langle\tau \sigma\rangle $,
	$H_{3}=\langle\tau, \sigma^2\rangle $,
et  trois sous-groupes d'ordre 2  d'indice
4 :
$H_{4}=\langle\tau\rangle $, $H_{5}=\langle\tau \sigma^2\rangle $ et
	$H_{6}=\langle\sigma^2\rangle .$

Posons  $H_7=\{1\}$, $L_0=K$ et $L_7=H^{+}_{K, 2}$.  La théorie de Galois nous permet  d'associer à chaque sous-groupe   $H_{i}$, $0\leq i \leq 7$, une unique extension ${L}_{i}/K$ abélienne non ramifiée
telle que $H_0/H_i\simeq Gal(L_i/K)$.

 Les seules extensions cycliques de degré $4$ non ramifiées sur $K$ sont $L_4$ et $L_5$.   D'après le Corollaire \ref{fr-cor3} (ii),
$ cl_{K}(\mathcal{O}_{L_4})=cl_{K}(\mathcal{O}_{L_5})=cl_{K}(\mathcal{O}_{L_3}).$
On a donc $R'_{nr}(K, \mathbb{Z}/4\mathbb{Z})=\{1, cl_{K}(\mathcal{O}_{L_3})\}$.
Comme $K$ est de type classe,  $cl_{K}(\mathcal{O}_{L_3})\not = 1$, par conséquent  $\left|R'_{nr}(K, \mathbb{Z}/4\mathbb{Z})\right|=2$.
\end{proof}

\noindent \textbf{Remarque} Sous les hypothèses de la proposition précédente, on a \break  $R'_{nr}(K, \mathbb{Z}/2\mathbb{Z})=\{1, cl_{K}(\mathcal{O}_{L_i}), 1\leq i \leq 3 \}=Cl_{K,2}$. $\blacksquare$
\medskip

\textbf{Exemple.}
Soit $K=\mathbb{Q}(\sqrt{p_0p_1p_2})$ un corps quadratique réel, où les $p_i$, $0\leq i \leq 2$,  sont des nombres premiers congrus à 1 modulo 4  satisfaisant :
$(\dfrac{p_1}{p_0})=(\dfrac{p_2}{p_0})=1$, $(\dfrac{p_1}{p_2})=-1$, $(\dfrac{p_0}{p_1p_2})_4=-1$ ou $(\dfrac{p_1p_2}{p_0})_4=-1$, où $(\dfrac{a}{b})_4$ est le symbole biquadratique.
D'après \cite{fr-Ka76}[Proposition B1, p. 326] (on pourrait voir aussi \cite{fr-Ka76}[Proposition B, p. 327]),
$syl_2(Cl^+_K)$ est de type $(2, 4)$. Il suit de la Proposition \ref{fr-reel} que $K$ est de type classe si et seulement si la norme de son unité fondamental  $\epsilon$ est $-1$.

  Soient $p_0=5$, $p_1=13$, $p_2=101$ et $K=\mathbb{Q}(\sqrt{5\cdot 13\cdot 101})=\mathbb{Q}(\sqrt{6565})$. En utilisant Pari/GP (voir \cite{fr-Pari}), on obtient : $syl_2(Cl^+_K)$ est de type $(2, 4)$ et $\epsilon=\frac{81+\sqrt{6565}}{2}$ ; la   norme de $\epsilon $ est  \ $-1$. Ce corps $K$ vérifie donc les hypothèses de la proposition précédente. $\blacksquare$
\smallskip

Nous terminons cet article par dire qu'au début de nos investigations nous avons prouvé le  Théorème \ref{fr-t2} pour les groupes abéliens, ensuite pour les groupes nilpotents. Nous avons constaté que la clé de la preuve dans  ces cas est la normalité du $2$-sous-groupe de Sylow, ce qui a abouti au Théorème \ref{fr-t2} du présent article.

Naturellement on s'est posé la question si notre méthode de démonstration dans le deuxième cas ($[N_{1}\cap N_{2}: K]$  impair) marche si  le $2$-sous-groupe de Sylow n'est pas normal (rappelons que dans le premier cas  nous n'avons pas besoin de  $syl_2(\Gamma)$  normal). Pour cela nous avons considéré le groupe diédral $\Gamma=D_n$ d'ordre $2n$, où $n$ est un entier impair $\geq 3$, et supposé $[N_{1}\cap N_{2}: K]$ impair ; il est immédiat que $[N_{1}\cap N_{2}: K]=1$ et donc $N_1/K$ et $N_2/K$ sont linéairement disjointes de sorte que $Gal(N_1N_2/K)\simeq  D_n\times D_n$.

Le groupe $D_n$
 a la présentation suivante : $D_n=\langle\sigma, \tau :  \sigma^n=\tau^2=1; \tau \sigma \tau=\sigma^{-1}\rangle $. Il est facile de montrer que $r_2(D_n\times D_n)=2$ et que les trois sous-groupes d'indice $2$ de  $D_n\times D_n$ sont :  $D_n\times \langle\sigma\rangle $, $\langle\sigma\rangle \times D_n$ et $(\langle\sigma\rangle \times \langle\sigma\rangle )\cup (\langle\sigma\rangle \tau \times \tau\langle\sigma\rangle )$.
 On a utilisé SageMath (\cite{fr-Sage})  pour voir que, lorsque $n\leq 100$,  $D_n\times \{1\}$ et $\{1\} \times D_n$ sont les seuls sous-groupes normaux  $J$ de $D_n\times D_n$ tels que $((D_n\times D_n)/J) \simeq D_n$ ; par la théorie de Galois, ils correspondent à $N_1/K$ et $N_2/K$.  On en déduit que $L_3/K$ ne peut pas être plongeable dans une sous-extension de $N_1N_2/K$ à groupe de Galois isomorphe à $D_n$ (c'est là la limite de notre méthode). Peut-être est-elle plongeable dans une $D_n$-extension non ramifiée en dehors de $N_1N_2/K$, ce qui pourrait  nous donner une piste pour  généraliser notre Théorème \ref{fr-t2}.

\end{document}